\input amstex
\magnification=1200
\loadmsbm
\loadeufm
\loadeusm
\UseAMSsymbols
\input amssym.def

\font\BIGtitle=cmr10 scaled\magstep3
\font\bigtitle=cmr10 scaled\magstep1
\font\boldsectionfont=cmb10 scaled\magstep1
\font\section=cmsy10 scaled\magstep1

\def\scr#1{{\fam\eusmfam\relax#1}}

\def\scrA{{\scr A}}
\def\scrB{{\scr B}}
\def\scrC{{\scr C}}
\def\scrD{{\scr D}}
\def\scrE{{\scr E}}
\def\scrF{{\scr F}}
\def\scrG{{\scr G}}
\def\scrH{{\scr H}}
\def\scrI{{\scr I}}
\def\scrL{{\scr L}}
\def\scrK{{\scr K}}
\def\scrJ{{\scr J}}
\def\scrM{{\scr M}}
\def\scrN{{\scr N}}

\def\scrP{{\scr P}}
\def\scrQ{{\scr Q}}
\def\scrS{{\scr S}}

\def\scrR{{\scr R}}
\def\scrT{{\scr T}}

\def\scrV{{\scr V}}
\def\scrX{{\scr X}}

\def\scrW{{\scr W}}

\def\gr#1{{\fam\eufmfam\relax#1}}

\def\grA{{\gr A}}

\def\grD{{\gr D}}	
\def\grE{{\gr E}}	
	
\def\grG{{\gr G}}	\def\grg{{\gr g}}

\def\grK{{\gr K}}

\def\grN{{\gr N}}	
	
	\def\grp{{\gr p}} 
\def\grQ{{\gr Q}}	
	
\def\grS{{\gr S}}	
\def\grT{{\gr T}}

\def\db#1{{\fam\msbfam\relax#1}}

\def\dbA{{\db A}} \def\dbB{{\db B}}
\def\dbC{{\db C}} \def\dbD{{\db D}}
 \def\dbF{{\db F}}
\def\dbG{{\db G}} 
\def\dbI{{\db I}} 
 
 \def\dbN{{\db N}}
 \def\dbP{{\db P}}
\def\dbQ{{\db Q}}

 \def\dbZ{{\db Z}}

\def\eps{{\varepsilon}}

\def\Itil{\widetilde{I}}
\def\Jtil{\widetilde{J}}

\def\Rtil{\widetilde{R}}

\def\Ker{\text{Ker}}
\def\Coker{\text{Coker}}
\def\der{\text{der}}

\def\Gal{\text{Gal}}
\def\Hom{\text{Hom}}
\def\End{\text{End}}
\def\Spec{\text{Spec}}
\def\Spf{\text{Spf}}

\def\Lie{\text{Lie}}

\def\leaderfill{\leaders\hbox to 1em
     {\hss.\hss}\hfill}
\def\nspace{\lineskip=1pt\baselineskip=12pt\lineskiplimit=0pt}

\def\finishproclaim{\par\rm
     \ifdim\lastskip<\medskipamount\removelastskip
     \penalty55\medskip\fi}
\def\endproof{$\hfill \square$}     
\def\proof{\par\noindent {\it Proof:}\enspace}
\def\references#1{\par
  \centerline{\boldsectionfont References}\bigskip
     \parindent=#1pt\nspace}
\def\Ref[#1]{\par\hang\indent\llap{\hbox to\parindent
     {[#1]\hfil\enspace}}\ignorespaces}
\def\Item#1{\par\hang\indent\llap{\hbox to\parindent
     {#1\hfill$\,\,$}}\ignorespaces}
\def\ItemItem#1{\par\indent\hangindent2\parindent
     \hbox to \parindent{#1\hfill\enspace}\ignorespaces}

\def\Le{{\mathchoice{\,{\scriptstyle\le}\,}
  {\,{\scriptstyle\le}\,}
  {\,{\scriptscriptstyle\le}\,}{\,{\scriptscriptstyle\le}\,}}}
\def\Ge{{\mathchoice{\,{\scriptstyle\ge}\,}
  {\,{\scriptstyle\ge}\,}
  {\,{\scriptscriptstyle\ge}\,}{\,{\scriptscriptstyle\ge}\,}}}

\def\arrowsim{\,\smash{\mathop{\to}\limits^{\lower1.5pt
  \hbox{$\scriptstyle\sim$}}}\,}

\def\doublemaprights#1#2#3#4{\raise3pt\hbox{$\mathop{\,\,\hbox to
     #1pt{\rightarrowfill}\kern-30pt\lower3.95pt\hbox to
     #2pt{\rightarrowfill}\,\,}\limits_{#3}^{#4}$}}

\def\rightcapdownarrow{\raise9pt\hbox{$\ssize\cap$}\kern-7.75pt
     \Big\downarrow}

\def\rcapmapdown#1{\rightcapdownarrow\kern-1.0pt\vcenter{
     \hbox{$\scriptstyle#1$}}}

\def\rmapdown#1{\Big\downarrow\kern-1.0pt\vcenter{
     \hbox{$\scriptstyle#1$}}}
\def\rightsubsetarrow#1{{\ssize\subset}\kern-4.5pt\lower2.85pt
     \hbox to #1pt{\rightarrowfill}}
\def\longtwoheadedrightarrow#1{\raise2.2pt\hbox to #1pt{\hrulefill}
     \!\!\!\twoheadrightarrow}

\def\Gal{\operatorname{\hbox{Gal}}}

\def\Hom{\operatorname{\hbox{Hom}}}

\def\im{\hbox{Im}}

\NoBlackBoxes
\parindent=25pt
\document
\footline={\hfil}

\null
\centerline{\BIGtitle A motivic conjecture of Milne}
\vskip 0.2in
\centerline{\bigtitle Adrian Vasiu, Binghamton University}
\vskip 0.2in
\centerline{January 18, 2011}
\footline={\hfill}
\null

\noindent
{\bf ABSTRACT.} Let $k$ be an algebraically closed field of characteristic $p>0$. Let $W(k)$ be the ring of Witt vectors with coefficients in $k$. We prove a motivic conjecture of Milne that relates, in the case of abelian schemes, the \'etale cohomology with $\dbZ_p$ coefficients to the crystalline cohomology with integral coefficients, in the more general context of $p$-divisible groups endowed with {\it arbitrary} families of crystalline tensors over a finite, discrete valuation ring extension of $W(k)$. This extends a result  of Faltings in [Fa2]. As a main new tool we construct global deformations of $p$-divisible groups endowed with crystalline tensors over certain regular, formally smooth schemes over $W(k)$ whose special fibers over $k$ have a Zariski dense set of $k$-valued points.



\bigskip\medskip

\footline={\hss\tenrm \folio\hss}
\pageno=1

\bigskip
\noindent
{\boldsectionfont 1. Introduction}
\bigskip\smallskip

Let $p$ be a prime. Let $k$ be a perfect field of characteristic $p$. Let $W(k)$ be the ring of Witt vectors with coefficients in $k$. Let $B(k)$ be the field of fractions of $W(k)$. In this paper we study {\it $p$-divisible groups} endowed with families of {\it \'etale} and {\it crystalline tensors} over (finite discrete valuation ring extensions of) $W(k)$. We begin by introducing families of tensors. 

Let $\Spec\,R$ be an affine scheme. For a finitely generated, projective $R$-module $P$, let $P^*:=\Hom_R(P,R)$ and let $\pmb{GL}_P$ be the reductive group scheme over $\Spec\,R$ of linear automorphisms of $P$. For $s\in\dbN$, let $P^{\otimes s}:=P\otimes_R \cdots\otimes_R P$, the number of copies of $P$ being $s$ (here $P^{\otimes 0}:=R$). By the {\it essential tensor algebra} of $P$ we mean 
$$\scrT(P):=\oplus_{s,t\in\dbN} P^{\otimes s}\otimes_R P^{*\otimes t}.$$
Let $\iota\in R$ be a non-divisor of $0$. A family of tensors of $\scrT(P[{1\over \iota}])=\scrT(P)[{1\over \iota}]$ is denoted $(w_{\alpha})_{\alpha\in\scrJ}$, with $\scrJ$ as the set of indexes. The scalar extensions of $w_{\alpha}$ are also denoted by $w_{\alpha}$. Let $P_1$ be another finitely generated, projective $R$-module. Let $(w_{1,\alpha})_{\alpha\in\scrJ}$ be a family of tensors of $\scrT(P_1[{1\over \iota}])$ indexed by the same set $\scrJ$. By an isomorphism 
$$(P,(w_{\alpha})_{\alpha\in\scrJ})\arrowsim (P_1,(w_{1,\alpha})_{\alpha\in\scrJ})$$ 
we mean an $R$-linear isomorphism $P\arrowsim P_1$ that extends naturally to an $R[{1\over {\iota}}]$-linear isomorphism $\scrT(P[{1\over \iota}])\arrowsim\scrT(P_1[{1\over \iota}])$ that takes $w_{\alpha}$ to $w_{1,\alpha}$ for all $\alpha\in\scrJ$. If $E$ is a group scheme over $\Spec\,R$, let $\Lie(E)$ be its {\it Lie algebra} over $R$. 

If $\sharp$ is a $\Spec\,R$-scheme (or a morphism of $\Spec\,R$-schemes or a $p$-divisible group over $\Spec\,R$) and if $\Spec\,Q\to\Spec\,R$ is an affine morphism, we define $\sharp_{Q}:=\sharp\times_{\Spec\,R} \Spec\,Q$. Similarly we define $\sharp_{Q}$ (resp. $\sharp_{*,Q}$), starting from $\sharp_R$ (resp. $\sharp_*$ with $*$ an index).

Let $\sigma:=\sigma_k$ be the Frobenius automorphism of $W(k)$ and $B(k)$ induced from $k$. We fix an algebraic closure $\overline{B(k)}$ of $B(k)$. For $K$ a subfield of $\overline{B(k)}$ that contains $B(k)$, let $\Gal(K):=\Gal(\overline{B(k)}/K)$. 

\medskip\smallskip\noindent
{\bf 1.1. On $p$-divisible groups over $\Spec\,W(k)$ endowed with tensors.} Let $D$ be a $p$-divisible group over $\Spec\,W(k)$. Let $(M,\phi)$ be the contravariant {\it Dieudonn\'e module} of $D_k$. Thus $M$ is a free $W(k)$-module whose rank equals to the height of $D$ and $\phi$ is a $\sigma$-linear endomorphism of $M$ such that we have $pM\subseteq\phi(M)$. Let $F^1$ be the direct summand of $M$ that is the {\it Hodge filtration} defined by $D$. We have $\phi(M+{1\over p}F^1)=M$. The rank of $F^1$ is the dimension of $D_k$. We call $(M,F^1,\phi)$ the {\it filtered Dieudonn\'e module} of $D$. Let $(F^i(M))_{i\in\dbZ}$ be the decreasing, exhaustive, and separated filtration of $M$ defined uniquely by $F^2(M):=0$, $F^1(M):=F^1$, and $F^0(M):=M$. Let $(F^i(M^*))_{i\in\dbZ}$ be the decreasing, exhaustive, and separated  filtration of $M^*$ defined uniquely by $F^{-1}(M^*):=M^*$, $F^0(M^*):=\{x\in M^*|x(F^1)=0\}$, and $F^1(M^*):=0$. We endow $\scrT(M)$ with the {\it tensor product filtration} $F^i(\scrT(M))_{i\in\dbZ}$ defined by $(F^i(M))_{i\in\dbZ}$ and $(F^i(M^*))_{i\in\dbZ}$. As $(F^i(\scrT(M)))_{i\in\dbZ}$ depends only on $F^1$, we call it the {\it filtration of $\scrT(M)$ defined by} $F^1$. 

For $\flat\in M^*[{1\over p}]$ let $\phi(\flat):=\sigma\circ \flat\circ\phi^{-1}\in M^*[{1\over p}]$. Thus $\phi$ acts in the usual tensor product way on $\scrT(M[{1\over p}])$. Under the identification $\End_{W(k)}(M)=M\otimes_{W(k)} M^*$, $\phi$ maps $\vartriangle\in\End_{B(k)}(M[{1\over p}])$ to $\phi\circ \vartriangle\circ\phi^{-1}\in\End_{B(k)}(M[{1\over p}])$. Let $(t_{\alpha})_{\alpha\in\scrJ}$ be a family of tensors of $F^0(\scrT(M))[{1\over p}]\subseteq\scrT(M[{1\over p}])$ which are fixed by $\phi$; thus $t_{\alpha}\in\{x\in F^0(\scrT(M))[{1\over p}]|\phi(x)=x\}$. Let $G$ be the schematic closure in $\pmb{GL}_M$ of the subgroup of $\pmb{GL}_{M[{1\over p}]}$ that fixes $t_{\alpha}$ for all $\alpha\in\scrJ$. It is a flat, closed subgroup scheme of $\pmb{GL}_M$ such that we have $\phi(\Lie(G_{B(k)}))=\Lie(G_{B(k)})$. Thus the pair $(\Lie(G_{B(k)}),\phi)$ is a Lie $F$-subisocrystal of $(\End_{B(k)}(M[{1\over p}]),\phi)$. 

Quadruples of the form $(M,\phi,G,(t_{\alpha})_{\alpha\in\scrJ})$ play key roles in the study of special fibres of good integral models of {\it Shimura varieties} of {\it Hodge type} in mixed characteristic $(0,p)$ (see [LR], [Ko], [Va1, Sect. 5], etc., for concrete situations with $G$ a reductive group scheme). 

Let $D^{\text{t}}$ be the {\it Cartier dual} of $D$. Let $H^1(D):=T_p(D^{\text{t}}_{B(k)})(-1)$ be the dual of the {\it Tate-module} $T_p(D_{B(k)})$ of $D_{B(k)}$. Thus $H^1(D)$ is a free $\dbZ_p$-module of the same rank as $M$ and $\Gal(B(k))$ acts on it. We suppress the Galois action whenever we consider $\pmb{GL}_{H^1(D)\otimes_{\dbZ_p} R}$ or pairs of the form $(H^1(D)\otimes_{\dbZ_p} R,(w_{\alpha})_{\alpha\in\scrJ})$, with $R$ as a commutative $\dbZ_p$-algebra. Let $v_{\alpha}\in\scrT(H^1(D)[{1\over p}])$ be the $\Gal(B(k))$-invariant tensor that corresponds to $t_{\alpha}$ via {\it Fontaine comparison theory} (see [Fo4, Subsect. 5.5], [Fa2, Sect. 6], Subsubsection 2.2.4, etc.). Let $G^{\acute{et}}$ be the schematic closure in $\pmb{GL}_{H^1(D)}$ of the subgroup of $\pmb{GL}_{H^1(D)[{1\over p}]}$ that fixes $v_{\alpha}$ for all $\alpha\in\scrJ$. The next definition plays a key role when $p=2$.

\medskip\noindent
{\bf 1.1.1. Definition.} We say that the {\it property (EC)} holds for $D$ (or that $D$ is {\it essentially of connected type}) if $D$  is a direct sum $D_1\oplus D_2$ such that $D_1$ and $D_2^{\text{t}}$ are connected (e.g., this automatically holds if $G^{\acute{et}}$ is a torus). 

\medskip
The first main goal of this paper is to construct {\it global deformations} of the pairs $(D,(t_{\alpha})_{\alpha\in\scrJ})$ and $(D,(v_{\alpha})_{\alpha\in\scrJ})$ over certain regular, affine, formally smooth $\Spec\,W(k)$-schemes $\Spec\,Q$ whose special fibres are geometrically connected and have a Zariski dense set of $\bar k$-valued points (see Theorem 3.4.1 and Subsubsection 3.4.2). The deformations of $(D,(t_{\alpha})_{\alpha\in\scrJ})$ will be obtained from deformations of the quadruple $(M,F^1,\phi,(t_{\alpha})_{\alpha\in\scrJ})$ chosen in such a way that the induced deformations of the triple $(M,F^1,(t_{\alpha})_{\alpha\in\scrJ})$ are constant. If $p=2$, then for these constructions we will assume that either $D$ or $D^{\text{t}}$ is connected. Each such global deformation of $(D,(t_{\alpha})_{\alpha\in\scrJ})$ over $\Spec\,Q$ will consist in a $p$-divisible group over $\Spec\,Q$ endowed naturally with a family of crystalline tensors indexed by the set $\scrJ$, whose pull-back at a suitable $k$-valued point of $\Spec\,Q$ is $(D,(t_{\alpha})_{\alpha\in\scrJ})$ (cf. Subsection 3.4). Each $Q$ will be the $p$-adic completion of a particular type of ind-\'etale algebra over a polynomial $W(k)$-algebra, cf. Subsection 3.1 and Theorem 3.2.

The second main goal of this paper is to use the mentioned global deformations in order to prove the following Main Theorem.

\medskip\smallskip\noindent
{\bf 1.2. Main Theorem.} {\it Suppose that $k=\bar k$. If $p=2$ we also assume that the property (EC) holds for $D$. Then there exists an isomorphism 
$$\rho:(M,(t_{\alpha})_{\alpha\in\scrJ})\arrowsim (H^1(D)\otimes_{\dbZ_p} W(k),(v_{\alpha})_{\alpha\in\scrJ}).$$}
\noindent
{\bf 1.3. Corollary.} {\it If $p=2$, we assume that the property (EC) holds for $D$. 

\smallskip
{\bf (a)} Then the pull-backs to $\Spec\,W(\bar k)$ of $G$ and $G^{\acute{et}}$ are isomorphic. In particular, $G$ is smooth (resp. reductive) if and only if $G^{\acute{et}}$ is smooth (resp. reductive). 

\smallskip
{\bf (b)} Let $\alpha\in\scrJ$. Then the tensor $t_{\alpha}$ belongs to $\scrT(M)$ (resp. to $\scrT(M)\setminus p\scrT(M)$) if and only if the tensor $v_{\alpha}$ belongs to $\scrT(H^1(D))$ (resp. to $\scrT(H^1(D))\setminus p\scrT(H^1(D))$).}

\medskip
Let $i\in\dbZ$. If $t\in F^i(\scrT(M))[{1\over p}]$ is such that $\phi^i(t)=p^it$, then Fontaine comparison theory associates to $t$ an \'etale Tate-cycle $v\in\scrT(H^1(D)[{1\over p}])$ whose $\dbZ_p$-span is a $\Gal(B(k))$-module isomorphic to $\dbZ_p(-i)$. Strictly speaking $v$ is unique up to $i$-th powers of units of $\dbZ_p$ i.e., up to a choice of a generator $\beta$ of a suitable free $\dbZ_p$-submodule of $B^+(W(k))$ of rank $1$ (see Subsection 2.2 for the integral crystalline {\it Fontaine ring} $B^+(W(k))$ and for its $\dbZ_p$-submodule $\dbZ_p\beta$). We have the following equivalent form of the Main Theorem which involves also {\it Tate-twists}, in which all \'etale Tate-cycles are defined using a fixed generator $\beta$ (see Subsubsection 2.2.4), and which is proved in Subsection 4.4.

\medskip\smallskip\noindent
{\bf 1.4. Corollary.} {\it Suppose that $k=\bar k$. If $p=2$ we also assume that the property (EC) holds for $D$. Let $(t_{\alpha})_{\alpha\in\scrJ_{\text{twist}}}$ be the family of all tensors in the set 
$$\cup_{i\in\dbZ} \{x\in F^i(\scrT(M))[{1\over p}]|\phi(x)=p^ix\},$$ 
with $\scrJ_{\text{twist}}$ as an index set. Let $\scrW$ be the family of all direct summands $W$ of the $W(k)$-module $\scrT(M)$ which have finite ranks and for which we have $W=\sum_{i\in\dbZ} p^{-i}\phi(W\cap F^i(\scrT(M)))$. For $\alpha\in\scrJ_{\text{twist}}$ (resp. for $W\in\scrW)$ let $v_{\alpha}$ (resp. $W^{\acute{et}}$) be the tensor of $\scrT(H^1(D))[{1\over p}]$ (resp. be the direct summand of $\scrT(H^1(D))$) that corresponds to $t_{\alpha}$ (resp. to $W$) via Fontaine comparison theory. Then there exists an isomorphism 
$$\rho_{\text{twist}}:(M,(t_{\alpha})_{\alpha\in\scrJ_{\text{twist}}})\arrowsim (H^1(D)\otimes_{\dbZ_p} W(k),(v_{\alpha})_{\alpha\in\scrJ_{\text{twist}}})$$ 
such that the isomorphism $\scrT(M)\arrowsim \scrT(H^1(D)\otimes_{\dbZ_p} W(k))$ induced by it maps each $W\in\scrW$ onto $W^{\acute{et}}\otimes_{\dbZ_p} W(k)$.}

\medskip\smallskip\noindent
{\bf 1.5. On literature and motivation.} The existence of $\rho:(M,(t_{\alpha})_{\alpha\in\scrJ})\arrowsim (H^1(D)\otimes_{\dbZ_p} W(k),(v_{\alpha})_{\alpha\in\scrJ})$ was conjectured by Milne if $G^{\acute{et}}$ is a reductive group scheme over $\Spec\,\dbZ_p$, $D$ is the $p$-divisible group of an {\it abelian scheme} $\scrA$ over $\Spec\,W(k)$, and each $t_{\alpha}$ and $v_{\alpha}$ are the de Rham component and the $p$-component of the \'etale component (respectively) of a {\it Hodge cycle} on $\scrA_{B(k)}$ (see [De2, Sect. 2] for Hodge cycles). References to {\it Milne's conjecture} in the context of abelian schemes can be found in [Mi3, property (4.16.3) and Subsect. 4.23], an unpublished manuscript of Milne dated Aug. 1995 which worked with $p>>0$ (see also [Mi4]), [Va1, Rm. 5.6.5 and Conj. 5.6.6], [Va6], and [Va7]. In fact the conjecture [Va1, Conj. 5.6.6] is a slight restatement of Milne's conjecture of 1995. In [Va6] (and an earlier version of [Va7]) we proved [Va1, Conj. 5.6.6] in some cases. The methods used in [Va6] (and in an earlier version of [Va7]) in connection to Milne's conjecture are very short and simple and of purely reductive group scheme theoretical nature. These methods can be viewed as refinements of [Ko, Lem. 7.2] and can be used to regain the Main Theorem only if the following two strong conditions hold:
\medskip
-- $G^{\acute{et}}$ and $G$ are reductive group schemes and moreover $G$ is generated by cocharacters of $\pmb{GL}_M$ that have weights only $0$ or $1$;
\smallskip
-- the maximal torus of the center of $G$ is the maximal torus of the center of the double centralizer of $G$ in either $\pmb{GL}_M$ or $\pmb{GSp}(M,\lambda_M)$ (resp. or $\pmb{GSO}(M,\lambda_M)$), where $\lambda_M:M\times M\to W(k)$ is an alternating bilinear form (resp. a symmetric bilinear form which modulo $2$ is alternating) normalized by $G$ and defined by an isomorphism $D\arrowsim D^{\text{t}}$.
\medskip
The recent paper [Ki3] (written after this paper was) regains the Main Theorem only in the following particular case (cf. [Ki3, Cor. (1.4.3)]):

\medskip
-- $G^{\acute{et}}$ is a reductive group scheme.

\medskip
Corollary 1.3 (b) is implied by Faltings' results [Fa2, Thm. 5 and Cor. 9] if and only if we have $t_{\alpha}\in \oplus_{s,t\in\{0,\ldots,p-2\}} (M^{\otimes s}\otimes_{W(k)} M^{*\otimes t})[{1\over p}]$.

In Milne's motivic approach to the proof of the {\it Langlands--Rapoport conjecture} (of [Mi2, Conj. 4.4]) for Shimura varieties of Hodge type, Milne's conjecture plays a key role (see [Mi4], [Va7], and [Va9]). We recall that Shimura varieties of Hodge type are {\it moduli spaces} of polarized abelian schemes endowed with level structures and with (specializations of) Hodge cycles and, due to this, they are also the main testing ground for many parts of the {\it Langlands program} (like zeta functions, local correspondences). Thus the importance of Milne's conjecture stems mainly from its meaningful applications to the Langlands program in general and to Shimura varieties in particular. Though for $p\Ge 5$, [Va1] and [Va7] combined can suffice for applications to the Langlands--Rapoport conjecture for Shimura varieties of Hodge type, for refined applications one needs the much stronger result of the Main Theorem. In future works we will use the Main Theorem and its proof as follows:

\medskip
{\bf (i)} To prove the existence of {\it integral canonical models} in mixed characteristics $(0,3)$ and $(0,2)$ of Shimura varieties of abelian type with respect to hyperspecial subgroups (see [Va1] for definitions and for the analogue result in mixed characteristic $(0,p)$ with $p\Ge 5$). See already [Va8] and [Ki3].  We emphasize that the methods used in [Ki3] work for $p=2$ only in exceptional cases while [Va8] also works for $p=2$. The key idea is that if $p=2$ and $G$ is a reductive group scheme, then one can always replace $D$ by another $2$-divisible group $D^{\prime}$ over $\Spec\,W(k)$ whose filtered $F$-crystal is the same as the one of $D$ and for which the Main Theorem holds (see [Va8, Part I, App. B] for details).  

\smallskip
{\bf (ii)} To study different {\it stratifications} of special fibres of good integral models of Shimura varieties of abelian type which are proved to exist in [Va8, Part I, Thm. 1.5] (or which can be proved to exist based on loc. cit.) and which are far more general than the integral canonical models mentioned in (i) (see [Va3, Sect. 8], [Va4, Sect. 4], [Va5], etc.). 

\smallskip
{\bf (iii)} To construct a comprehensive theory of {\it automorphic vector bundles} on integral good models of Shimura varieties of abelian type (presently this theory is not complete even for Siegel modular varieties, cf. end of [FC, Ch. VI, p. 238]).

\smallskip
{\bf (iv)} To get analogues of the Main Theorem and Corollary 1.3 for different {\it motives} associated to classes of polarized projective schemes whose moduli spaces are related to Shimura varieties (see [An] for such classes; we have in mind mainly the case of $H^2_{\text{crys}}$ and $H^2_{\acute{et}}$ groups of {\it hyperk\"ahler schemes}).

\medskip
The applications (ii) to (iv) are not in the reach of [Ki3, Cor. (1.4.3)] as they involve general group schemes $G^{\acute{et}}$ and $G$ which are far from being reductive (often they are not even smooth) and thus as they require the Main Theorem in its full form. 

Both the Main Theorem and the Corollary 1.3 do not hold in general for $p=2$ (like when $D$ is isogenous but not isomorphic to $\dbQ_2/\dbZ_2\oplus\mu_{2^{\infty}}$ and $G$ is isomorphic to $\dbG_m$). 

\medskip\smallskip\noindent
{\bf 1.6. On contents.} The proof of the Main Theorem involves two steps. The first step shows that it is enough to prove the Main Theorem if $G$ is a torus (see Subsection 3.5) and the second step proves the Main Theorem if $G$ is a torus (see Subsections 4.1 and 4.2). The proof of the Main Theorem ends with Subsubsection 4.2.11.

The first step relies on the existence of global deformations of the pairs $(D,(t_{\alpha})_{\alpha\in\scrJ})$ and $(D,(v_{\alpha})_{\alpha\in\scrJ})$ over certain regular, affine, formally smooth $\Spec\,W(k)$-schemes $\Spec\,Q$ whose special fibres are geometrically connected and have a Zariski dense set of $\bar k$-valued points. Based on [Fa1, Thm. 7.1] and on a variant of it for $p=2$, the existence of such deformations boils down to proving the existence of certain integrable and nilpotent modulo $p$ {\it connections} on $M\otimes_{W(k)} Q$ (see  Subsection 3.4). See Subsections 3.2 and 3.3 for our basic results which pertain to moduli spaces of connections. See Theorem 3.2 (c) for a precise description of such $Q$. Though these global deformations are of interest in their own (for instance, they can be used to solve essentially any specialization problem which pertains to pairs of the form $(D_k,(t_{\alpha})_{\alpha\in\scrJ})$ and thus to the stratifications mentioned above in 1.5 (ii)), here we will use them only to accomplish the first step. 

The second step relies on the {\it integral version of Fontaine comparison theory} studied in [Fa2, Sect. 4] and it involves a non-trivial reduction to the {\it Lubin--Tate} case in which $k=\bar k$, $F^1$ has rank $1$, the $F$-isocrystal $(M[{1\over p}],\phi)$ over $k$ is simple, and $G$ is a torus. It is easy to prove Theorem 1.2 in the Lubin--Tate case, cf. Example 4.1.3 and Lemma 4.2.1.  

Section 2 gathers prerequisites for Sections 3 to 5. Corollary 1.4 is proved in Subsection 4.4. See Subsections 4.3 and 4.5 and Section 5 for a $\dbZ_p$ variant and complements to the Main Theorem and to Corollary 1.3. If $p\ge 3$, then each variant of the Main Theorem which is over the spectrum of a finite, totally ramified discrete valuation ring extension of $W(k)$ and which  involves a setting one expects to be deformable in a good way to $\Spec\,W(k)$, follows from the Main Theorem and  from a refinement of the deformation theory of [Fa2, Ch. 7] (see Subsections 5.3 and 5.4). See Subsection 5.5 for a crystalline variant (converse) to [Ki3, Cor. (1.4.3)] for $p>2$. Analogues of the Main Theorem do not exist in general for ramified settings that are not deformable in a good way to $\Spec\,W(k)$, cf. Subsection 5.6. At the first reading, one could assume everywhere that $k$ is algebraically closed.

\bigskip\smallskip
\noindent
{\boldsectionfont 2. Preliminaries}
\bigskip\smallskip

See Subsection 2.1 for conventions and notations to be used throughout the paper. In Subsections 2.2 and 2.3 we recall some facts which pertain to Fontaine comparison theories and to {\it Faltings--Fontaine categories}. Subsection 2.4 introduces Artin--Schreier systems of equations; they play a key role in Section 3. In Subsections 2.5 to 2.6 we include a few simple group scheme theoretical properties. Subsection 2.7 reviews and complements [Wi].

\medskip\smallskip\noindent
{\bf 2.1. Conventions and notations.} If $f_1$ and $f_2$ are endomorphisms of a group, we often denote $f_1\circ f_2$ by $f_1f_2$. For $n\in\dbN^{\ast}$ let $W_n(k):=W(k)/p^nW(k)$. Let $R^\wedge$ be the $p$-adic completion of a commutative, flat $W(k)$-algebra $R$. If $Y=\Spec\,R$ let $Y^\wedge:=\Spec\,R^\wedge$. By a {\it Frobenius lift} of $R^\wedge$ (resp. of $Y^\wedge$) we mean an endomorphism of $R^\wedge$ (resp. of $Y^\wedge$) that lifts the Frobenius endomorphism of $R/pR$ (resp. of $Y^\wedge_k=Y_k$). If $\bot:Z\to Y$ is a morphism of $\Spec\,W(k)$-schemes, then by its special fibre we mean $Z_k\to Y_k$, by its generic fibre we mean $Z_{B(k)}\to Y_{B(k)}$, and by $\bot$ modulo $p^n$ we mean $Z_{W_n(k)}\to Y_{W_n(k)}$. If $\Spec\,Q\to \Spec\,R$ is a morphism between affine, flat $\Spec\,W(k)$-schemes, then whenever we say that it (or the induced $W(k)$-homomorphism $R\to Q$) is formally smooth (or \'etale), both the $W(k)$-algebras $R$ and $Q$ are endowed with the $p$-adic topology. If $\hat R$ is the completion of $R$ with respect to an ideal $I$ of $R$ that contains $p$ and if $\Spf\,\hat R$ is the corresponding formal scheme, then we identify canonically the $\dbZ_p$-linear categories of $p$-divisible groups over $\Spec\,\hat R$ and respectively over $\Spf\,\hat R$ (cf. [Me, Ch. II, Lem. 4.16]); thus we will often use the Grothendieck--Messing deformation theory of [Me, Ch. 5] to lift $p$-divisible groups over $\Spec\,R/I$ to $p$-divisible groups over $\Spec\,\hat R$. 

Let $\delta(p)$ denote the unique divided power structure on the ideal $(p)$ of $R$ and the induced divided power structure on the ideal $(p)$ of $R/p^nR$.  Let $CRIS(Y_k/\Spec\,W(k))$ be the Berthelot crystalline site of [Be, Ch. III, Sect. 4]. It is convenient to call both pairs $(Y_k\hookrightarrow Y_{W_n(k)},\delta(p))$ and $(Y_k\hookrightarrow Y^\wedge,\delta(p))$ as thickenings, even though only the first pair is a {\it thickening} (i.e., an {\it \'epaississement \`a puissances divis\'ees}) in the sense of [Be, Ch. III, Def. 1.1.1]. For {\it Dieudonn\'e crystals} on $CRIS(Y_k/\Spec\,W(k))$ of either $p$-divisible groups or finite, flat, commutative  group schemes annihilated by a power of $p$ over $Y_k$, we refer to [BBM, Ch. 3] and [BM, Chs. 2 and 3]. If $\scrD_{Y_k}$ is a $p$-divisible group over $Y_k$ and if $\Phi_R$ is a fixed Frobenius lift of $R^\wedge$ compatible with $\sigma$, then by the {\it evaluation} $(N,\phi_N,\nabla_N)$ of the Dieudonn\'e crystal $\dbD(\scrD_{Y_k})$ on $CRIS(Y_k/\Spec\,W(k))$ at the thickening $(Y_k\hookrightarrow Y^\wedge,\delta(p))$ we mean the projective limit indexed by $n$ of the evaluation of $\dbD(\scrD_{Y_k})$ at the thickening $(Y_k\hookrightarrow Y_{W_n(k)},\delta(p))$ of $CRIS(Y_k/\Spec\,W(k))$, the Verschiebung maps being suppressed. Thus $N$ is a projective $R^\wedge$-module, $\phi_N$ is a $\Phi_R$-linear endomorphism of $N$, and $\nabla_N$ is an integrable, nilpotent modulo $p$ connection on $N$. By the {\it filtered Dieudonn\'e crystal} of a $p$-divisible group $\scrD$ over $Y$ we mean $\dbD(\scrD_{Y_k})$ endowed with the Hodge filtration defined by $\scrD$; thus the evaluation of the mentioned filtered Dieudonn\'e crystal at the thickening $(Y_k\hookrightarrow Y^\wedge,\delta(p))$ is the quadruple $(N,F^1_N,\phi_N,\nabla_N)$, where $F^1_N\subseteq N$ is the Hodge filtration defined by $\scrD$. If $p>2$, then $N/p^nN$ is the Lie algebra of the universal vector extension of the Cartier dual $\scrD^{\text{t}}_{Y_{W_n(k)}}$  of $\scrD_{Y_{W_n(k)}}$ (cf. [BM, Cor. 3.2.11] and [Me, Ch. 4]). 

A morphism of schemes will be called  an \'etale cover (resp. a pro-\'etale cover) if it is finite \'etale (resp. if it is a projective limit of finite \'etale morphisms). Let $K_A$ be the field of fractions of a $\dbZ$-algebra $A$ that is an integral domain. Let $\overline{W(k)}$ be the normalization of $W(k)$ in $\overline{B(k)}$. Let $V(k):=\overline{W(k)}^\wedge$. Let $K(k):=V(k)[{1\over p}]$.

The notations $D$, $(M,\phi)$, $F^1$, $(t_{\alpha})_{\alpha\in\scrJ}$, $G$, $(v_{\alpha})_{\alpha\in\scrJ}$, and $G^{\acute{et}}$ will always be as in Subsection 1.1. Let $\mu:\dbG_m\to G$ be the factorization through $G$ of the inverse of the {\it canonical split cocharacter} $\mu_{\text{can}}:\dbG_m\to \pmb{GL}_M$ of $(M,F^1,\phi)$ defined in [Wi, p. 512] (the cocharacter $\mu_{\text{can}}$ fixes each $t_{\alpha}$, cf. the functorial aspects of [Wi, p. 513]). The cocharacter $\mu$ produces a direct sum decomposition $M=F^1\oplus F^0$ such that $\triangle\in\dbG_m(W(k))$ acts through $\mu$ on $F^i$ as the multiplication with $\triangle^{-i}$, $i\in\{0,1\}$. Note that $F^1=F^1(M)$ but $F^0$ is in general only a direct summand of $F^0(M)=M$. 
 
\medskip\smallskip\noindent
{\bf 2.2. Fontaine comparison theory.} 
For the below review we refer to [Fo3] and [Fa2, Sect. 2]. Let $A_k$ be the perfect integral domain of
sequences $(x_n)_{n\in\dbN}$ of $V(k)/pV(k)$ which satisfy the identity $x_n=x_{n+1}^p$ for all $n\in\dbN$. The Galois group $\Gal(B(k))$ acts on $V(k)$ and thus also on $A_k$. The $\Gal(B(k))$-module $\dbQ_p(1)$ can be identified with sequences $(\eta_n)_{n\in\dbN}$ in $V(k)$
 of $p$-power roots of unity
which satisfy the identity $\eta_n=\eta_{n+1}^p$ for all $n\in\dbN$. The inclusion $\dbZ_p(1)\subset \dbQ_p(1)$ is defined by the extra identity $\eta_0=1$. 
Taking such sequences modulo $p$, we get a group homomorphism $\eta_k\colon\dbQ_p(1)\to\dbG_m(A_k)$ 
that respects the Galois actions.
For an element $z\in V(k)$, we choose a sequence $(z_n)_{n\in\dbN}$ of elements of $V(k)$ with the properties that $z_0=z$ and that $z_n=z_{n+1}^p$ for all $n\in\dbN$. Taking this sequence modulo $p$ we obtain an element 
$\underline{z}\in A_k$, well
defined by $z$ up to multiplication with an element of $\eta_k(\dbZ_p(1))$. If $x\in V(k)/pV(k)$ let $\tilde x\in V(k)$ be a lift of it. The ring $W(A_k)$ of Witt vectors with coefficients in $A_k$ is an integral domain as $A_k$ is so. 

The rule $x\to (\sigma^{-n}(x))_{n\in\dbN}$ defines a $k$-algebra structure on $A_k$. Let $s_k\colon W(A_k)\twoheadrightarrow V(k)$ 
be
the $W(k)$-epimorphism defined by (to be compared with [Fo3, Subsubsect. 1.2.2]):
$$s_k((x_0,x_1,\ldots))=\sum_{n\ge 0} p^nx_n^{(n)},$$
where $(x_{n,m})_{m\in\dbN}$ is the sequence of elements of $V(k)/pV(k)$ that defines $x_n\in A_k$ and where $x_n^{(n)}$ is the $p$-adic limit in $V(k)$ of the sequence $(\tilde x_{n,n+m}^{p^m})_{m\in\dbN}$ (it does not depend on the choice of the lifts $\tilde x_{n,n+m}$'s). Let 
$$\xi:=(\underline{p},\underline{-1}^p,0,0,\ldots)\in W(A_k).$$ 
We have $s_k(\xi)=p+ p(-1)=0$. The kernel of $s_k$ modulo $p$ is generated by $\xi$ modulo $p$. From the last two sentences we get that $\xi$ generates $\Ker(s_k)$.

Let $W^+(A_k)$ be the $W(A_k)$-subalgebra of the field of fractions of $W(A_k)$ (equivalently of $W(A_k)[{1\over p}]$) generated by all ${{\xi^n}\over {n!}}$ with $n\in\dbN^{\ast}$. It is equipped with a $W(k)$-epimorphism $s_k^+\colon W^+(A_k)\twoheadrightarrow V(k)$ which extends $s_k$ and which maps each ${{\xi^n}\over {n!}}$ to $0$. The kernel $\Ker(s_k^+)$ has a canonical divided power structure. Working with $W(k)$-algebras, $W^+(A_k)$ is the divided power hull of $\Ker(s_k^+)$ (cf. [Fo3, Subsubsect. 2.2.3] and [Fo2, Thm. 1 (ii), p. 98]).  

Let $B^+(W(k)):=W^+(A_k)^\wedge$; it is a $p$-adically complete integral domain that is a $W^+(A_k)$-algebra and thus also a $W(k)$-algebra. Note down that $B^+(W(k))$ is the integral crystalline Fontaine ring $A_{\text{cris}}$ of [Fo3] and [Fo4], cf. [Fo3, Subsubsect. 2.3.3]. The $W(k)$-algebra $B^+(W(k))$ has a decreasing filtration $(F^n(B^+(W(k))))_{n\in\dbZ}$ by ideals, where $F^n(B^+(W(k)))$ is the $p$-adic completion of $\Ker(s_k^+)^{[n]}$ (here $\Ker(s_k^+)^{[n]}:=W^+(A_k)$ if $n\le 0$ and $\Ker(s_k^+)^{[n]}$ is the $n$-th divided power ideal associated to the divided power structure on $\Ker(s_k^+)$ if $n> 0$). It is easy to see that the $W(A_k)$-module $W^+(A_k)/\Ker(s_k^+)^{[n]}$ is finitely presented. Thus we have a canonical identification $B^+(W(k))/F^n(B^+(W(k)))=W^+(A_k)/\Ker(s_k^+)^{[n]}$ of $p$-adically complete $W(A_k)$-modules. 

As the canonical Frobenius endomorphism of $W(A_k)$ maps $\xi$ to $\xi^p+pW(A_k)\subset pB^+(W(k))$, it extends to a Frobenius lift $\Phi_k$ of $B^+(W(k))$. 

The group $\Gal(B(k))$ acts on the filtered $W(k)$-algebra $B^+(W(k))$. The $W(k)$-epimorphism $s_k^+:W^+(A_k)\twoheadrightarrow V(k)$ induces naturally a $W(k)$-epimorphism 
$$s_k^+:B^+(W(k))\twoheadrightarrow B^+(W(k))/F^1(B^+(W(k)))=V(k).$$ 
\indent
We have a Teichm\"uller monomorphism $\dbG_m(A_k)\hookrightarrow \dbG_m(W(A_k))$ defined by the rule $a\mapsto (a,0,0,\ldots)$. Composing the restriction of $\eta_k$ to $\dbZ_p(1)$ with this Teichm\"uller monomorphism, we get a homomorphism $t_k:\dbZ_p(1)\to\dbG_m(W(A_k))$. Let $v_k:\dbZ_p(1)\to\dbG_m(B^+(W(k)))$ be the composite of $t_k$ with the natural monomorphism $\dbG_m(W(A_k))\hookrightarrow\dbG_m(B^+(W(k)))$. The composite $s_k^+\circ v_k$ is the trivial homomorphism $\dbZ_p(1)\to\dbG_m(V(k))$. Thus let (to be compared with [Fo3, Subsubsect. 2.3.4])
$$\beta_k\colon\dbZ_p(1)\to F^1(B^+(W(k)))$$ 
be the homomorphism obtained by taking log of $v_k$ (i.e., by inserting $v_k$ into the power series $\log(X)=\sum_{n=1}^{\infty} (-1)^{n+1}{{(X-1)^n}\over n}$ centered at $1$). Throughout the paper we {\it fix} a generator of $\dbZ_p(1)$ and call $\beta$ its image in $F^1(B^+(W(k)))$ through $\beta_k$. As $\Ker(s_k^+)$ has a divided power structure, we have ${{\beta^p}\over p}\in B^+(W(k))$. Thus $B^+(W(k))[{1\over {\beta}}]$ is a $B(k)$-algebra. We have $\Phi_k\circ\beta_k=p\beta_k$.

As the $W(A_k)$-algebras 
$$B^{++}(W(k)):=\text{proj}.\text{lim}._{n\in\dbN^{\ast}} B^+(W(k))/F^n(B^+(W(k)))$$ and 
$\text{proj}.\text{lim}._{n\in\dbN^{\ast}} W^+(A_k)/F^n(B^+(W(k)))\cap W^+(A_k)=\text{proj}.\text{lim}._{n\in\dbN^{\ast}} W^+(A_k)/\Ker(s_k^+)^{[n]}$ are naturally identified, the $W(A_k)$-monomorphism $B^+(W(k))\hookrightarrow B^{++}(W(k))$ gives rise for each $n\in\dbN^{\ast}$ to an identification 
$$F^n(B^+(W(k)))=B^+(W(k))\cap\Ker(B^{++}(W(k))\twoheadrightarrow B^+(W(k))/F^n(B^+(W(k)))).$$

\noindent
{\bf 2.2.1. Key facts on $B^+(W(k))$.} Let $\bar\xi\in A_k$ be the reduction modulo $p$ of  $\xi\in W(A_k)$. We have isomorphisms
$$q_{k,p}:B^+(W(k))/(F^p(B^+(W(k)))+pB^+(W(k)))=W(A_k)/(\xi^p,p)=A_k/(\bar\xi^p)\arrowsim V(k)/pV(k),$$ 
where all except the last one are canonical identifications and where the last isomorphism is defined by the epimorphism $A_k\twoheadrightarrow V(k)/pV(k)$ that takes $(x_n)_{n\in\dbN}\in A_k$ to $x_1$ (see [Fa1, Sect. 2, p. 30] and [Fa2, Sect. 4, top of p. 126]). The isomorphism $q_{k,p}$ is defined naturally by a $\sigma^{-1}$-linear $\dbZ_p$-epimorphism
$$q_k:B^+(W(k))\twoheadrightarrow V(k)/pV(k).$$

We fix a $p$-th root $p^{1\over p}$ of $p$ such that $q_k(\xi)$ is $p^{1\over p}$ modulo $p$; thus, if $(x_n)_{n\in\dbN}\in A_k$ defines $\underline{p}$, then $x_1$ is $p^{1\over p}$ modulo $p$. For $i\in\{1,\ldots,p-1\}$, the ideal $q_k(F^i(B^+(W(k))))$ of $V(k)/pV(k)$ is generated by $p^{i\over p}$ modulo $p$. Let $\Phi_{1k}:F^1(B^+(W(k)))\to B^+(W(k))$ be the map defined by the rule: if $z\in F^1(B^+(W(k)))$, then we have an equality $\Phi_k(z)=p\Phi_{1k}(z)$. As $\Phi_{k}(\xi)=(\underline{p}^p,\underline{-1}^{p^2},0,\ldots)$, we easily get that $\Phi_{k}(\xi)-\xi^p-p(\underline{-1}^{p},0,\ldots)\in p^2W(A_k)$. As ${{\xi^p}\over p}\in F^p(B^+(W(k)))$ and as $(-1)^p$ and $-1$ are congruent modulo $p$, we get that $\Phi_{1k}(\xi)-(\underline{-1}^{p^2},0,0,\ldots)\in\Ker(q_k)$. Thus the Frobenius lift $\Phi_k$ of $B^+(W(k))$ and the $\Phi_k$-linear map $\Phi_{1k}$ become under $q_k$ the Frobenius endomorphism $\bar\Phi_k$ of $V(k)/pV(k)$ and respectively the map 
$\bar\Phi_{1k}:p^{1\over p}V(k)/pV(k)\to V(k)/pV(k)$
 that takes $p^{1\over p}z$ modulo $p$ to $(-z)^p$ (equivalently to $-z^p$) modulo $p$; here $z\in V(k)$. Let 
$$\dbD(B^+(W(k))/pB^+(W(k)))$$ 
be the reduction of the quadruple $(B^+(W(k)),F^1(B^+(W(k))),\Phi_k,\Phi_{1k})$ modulo $p$ and let 
$$\dbD(V(k)/pV(k)):=(V(k)/pV(k),p^{1\over p}V(k)/pV(k),\bar\Phi_k,\bar\Phi_{1k}).$$  

Let $\text{gr}^0:=V(k)=F^0(B^+(W(k)))/F^1(B^+(W(k)))$. The $V(k)$-module 
$$\text{gr}^1:=F^1(B^+(W(k)))/F^2(B^+(W(k)))$$ 
is generated by the image $\xi_1$ of $\xi$ in $\text{gr}^1$. We recall the standard argument that $\text{gr}^1$ is a free $V(k)$-module of rank $1$. We show that the assumption that the cyclic $V(k)$-module $\text{gr}^1$ is not free leads to a contradiction. As each non-zero ideal of $V(k)$ contains a power of $p$, our assumption implies that there exists $n\in\dbN^{\ast}$ such that $p^n\text{gr}^1=0$. Thus $p^n\xi$  is a linear combination of $\xi^2,\xi^3, \ldots$ with coefficients in $W(A_k)[{1\over p}]$; by enlarging $n$ and by clearing out the powers of $p$ from denominators we can assume that in fact all coefficients are in $W(A_k)$. As $W(A_k)$ is an integral domain we get that $p^n$ belongs to the ideal $(\xi)$ of $W(A_k)$ and this contradicts the fact that $W(A_k)/(\xi)=V(k)$ is torsion free. 

As $\text{gr}^1$ is a free $V(k)$-module of rank $1$, the image $\beta_1$ of $\beta$ in $\text{gr}^1$ is $\psi\xi_1$, where $\psi\in V(k)$ is a $\dbG_m(V(k))$-multiple of a $p-1$-th root of $p$ (cf. [Fo3, Subsubsects. 5.1.2 and 5.2.4] where the pair $(\beta,\xi)$ is denoted as $(t,\xi)$; $\psi$ corresponds to $\pi_{\eps}^\prime v_0^{-1}$ of loc. cit.). 

This implies that $q_k(\beta)$ is $q_k(\xi)$ times a $\dbG_m(V(k)/pV(k))$-multiple of the reduction modulo $p$ of a $p(p-1)$-th root of $p$. As ${1\over p}+{1\over {p(p-1)}}={1\over {p-1}}$, we get that $q_k(\beta)$ is the image in $V(k)/pV(k)$ of a $\dbG_m(V(k))$-multiple of a $p-1$-th root of $p$. 

\medskip\noindent
{\bf 2.2.2. First applications to $D$.} 
We use the notations of Subsection 1.1. We consider the decreasing, exhaustive, and separated filtration $(F^i(H^1(D)))_{i\in\dbZ}$ of $H^1(D)$ defined by $F^1(H^1(D)):=0$ and $F^0(H^1(D)):=H^1(D)$. We endow $M$ with the trivial $\Gal(B(k))$-action. 

The Fontaine comparison theory provides us with a $B^+(W(k))$-linear monomorphism
$$i_D:M\otimes_{W(k)} B^+(W(k))\hookrightarrow H^1(D)\otimes_{\dbZ_p} B^+(W(k))\leqno (1)$$
that respects the tensor product filtrations and the Galois actions. Following [Fa2, Sect. 6], we recall how $i_D$ is constructed. Let $j:\dbQ_p/\dbZ_p\to D_{V(k)}$ be a homomorphism of $p$-divisible groups over $\Spec\,V(k)$. For $n\in\dbN^{\ast}$ let $\delta^{\text{can}}_n$ be the natural divided power structure of the ideal $\text{Im}(F^1(B^+(W(k)))+pB^+(W(k))\to B^+(W(k))/p^nB^+(W(k)))$ of $B^+(W(k))/p^nB^+(W(k))$. By evaluating $\dbD(j):\dbD(D_{V(k)})\to \dbD(\dbQ_p/\dbZ_p)$ at the thickening 
$(\Spec\,V(k)/pV(k)\hookrightarrow \Spec\,B^+(W(k))/p^nB^+(W(k)),\delta_n^{\text{can}}),$
we get a map 
$$M/p^nM\otimes_{W_n(k)} B^+(W(k))/p^nB^+(W(k))\to B^+(W(k))/p^nB^+(W(k))$$
which is $B^+(W(k))/p^nB^+(W(k))$-linear. Passing to limit $n\to\infty$ we get a $B^+(W(k))$-linear map $M\otimes_{W(k)} B^+(W(k))\to B^+(W(k))$ i.e., an element of $M^*\otimes_{W(k)} B^+(W(k))$. Varying $j$ we get a $B^+(W(k))$-linear monomorphism 
$$i_D^*:T_p(D_{B(k)})\otimes_{\dbZ_p} B^+(W(k))\hookrightarrow M^*\otimes_{W(k)} B^+(W(k)).$$ 
The transpose (dual) of $i_D^*$ is $i_D$. 

\medskip\noindent
{\bf 2.2.3. Key facts on $i_D$.} We check that $\beta$ annihilates $\text{Coker}(i_D)$ and $\text{Coker}(i_D^*)$. Considering homomorphisms $\dbQ_p/\dbZ_p\to D_{V(k)}\to\mu_{p^{\infty}}$, it is enough to handle the case when $D=\mu_{p^{\infty}}^m$, with $m\in\dbN^{\ast}$. In this case we have $\im(i_D)=H^1(D)\otimes_{\dbZ_p} \beta B^+(W(k))$, cf. the definition of $\beta_k$ and [BM, Cor. 2.2.4]; thus $\beta$ annihilates $\text{Coker}(i_D)$ and there exists 
$$(\beta i_D^{-1}):H^1(D)\otimes_{\dbZ_p} B^+(W(k))\to M\otimes_{W(k)} B^+(W(k))$$ such that the endomorphisms $(\beta i_D^{-1})i_D$ and $i_D(\beta i_D^{-1})$ are the multiplication by $\beta$. By transposition (duality), we get that $\beta$ annihilates $\text{Coker}(i_D^*)$.

The $B^+(W(k))$-linear monomorphism $i_D$ is strict with respect to filtrations, cf. [Fa2, Thms. 5 and 7] and [Fa2, Ch. 8] (due to the paragraph before Subsubsection 2.2.1, to check this property it is irrelevant if we use $i_D$ or its tensorization over $B^+(W(k))$ with  $B^{++}(W(k))$). This strictness property implies that the $V(k)$-linear map
$$j_D:F^0\otimes_{W(k)} \text{gr}^1\oplus F^1\otimes_{W(k)} \text{gr}^0\hookrightarrow H^1(D)\otimes_{\dbZ_p} \text{gr}^1\leqno (2)$$ defined by $i_D$ at the level of one gradings, is injective. As $\text{Coker}(i_D)$ is annihilated by $\beta$, $\text{Coker}(j_D)$ is annihilated by the element $\psi\in V(k)$ of Subsubsection 2.2.1. 

\medskip\noindent
{\bf 2.2.4. \'Etale Tate-cycles and $\dbZ_p$-submodules.} The category of crystalline representations of $\Gal(B(k))$ over $\dbQ_p$ is stable under tensor products and duals, cf. [Fo4, Subsect. 5.5]. This implies that there exists a tensor $v_{\alpha}\in\scrT(H^1(D)[{1\over p}])$ that corresponds to $t_{\alpha}$ via the $B^+(W(k))[{1\over \beta}]$-linear isomorphism 
$$\scrT(i_D[{1\over \beta}]):\scrT(M[{1\over p}])\otimes_{B(k)} B^+(W(k))[{1\over \beta}]\arrowsim\scrT(H^1(D)[{1\over p}])\otimes_{\dbQ_p} B^+(W(k))[{1\over \beta}]$$
defined by $i_D[{1\over \beta}]$. The tensor $v_{\alpha}$ is fixed by $\Gal(B(k))$. 

Let $\mu(\beta)=\mu_{\text{can}}(\beta^{-1})$ be the $B^+(W(k))[{1\over \beta}]$-linear automorphism of $M\otimes_{W(k)} B^+(W(k))[{1\over \beta}]$ defined by the evaluation of the cocharacter $\mu:\dbG_m\to G$ at $\beta$. If $i\in\dbZ$ and $t\in\{x\in F^i(\scrT(M))[{1\over p}]|\phi(x)=p^ix\}$, then the $B^+(W(k))[{1\over \beta}]$-linear isomorphism $\scrT(M)\otimes_{W(k)} B^+(W(k))[{1\over \beta}]\arrowsim\scrT(H^1(D))\otimes_{\dbZ_p} B^+(W(k))[{1\over \beta}]$  defined by $i_D[{1\over \beta}]\circ\mu(\beta)$ takes $t$ to a tensor $v\in\scrT(H^1(D))[{1\over p}]$ whose $\dbZ_p$-span $<v>$ is a $\Gal(B(k))$-module isomorphic to $\dbZ_p(-i)$. We say that $v$ corresponds (or is associated) to $t$ via Fontaine comparison theory for $D$ (and via the choice of the generator $\beta$ of the free $\dbZ_p$-module $\im(\beta_k)$).

Let $W$ be a direct summand of the $W(k)$-module $\scrT(M)$ which has finite rank and for which we have an identity 
$$W=\sum_{i\in\dbZ} p^{-i}\phi(W\cap F^i(\scrT(M))).$$ 
Then Fontaine comparison theory shows that there exists a unique direct summand $W^{\acute{et}}$ of the $\dbZ_p$-module $\scrT(H^1(D))$ such that the above $B^+(W(k))[{1\over \beta}]$-linear isomorphism $\scrT(i_D[{1\over \beta}])$ maps $W\otimes_{W(k)} B^+(W(k))[{1\over \beta}]$ onto $W^{\acute{et}}\otimes_{\dbZ_p} B^+(W(k))[{1\over \beta}]$. 

\medskip
The next two results will be often used to reconstruct $D$ and $\End(D)$ from $(M,F^1,\phi)$ and $\End((M,F^1,\phi))$ (respectively).

\medskip\noindent
{\bf 2.2.5. Lemma.} {\it We recall that $D$ is a $p$-divisible group over $\Spec\,W(k)$. If $p=2$, we assume that either $D_k$ or $D_k^{\text{t}}$ is connected. Then $D$ is determined up to unique isomorphism by its filtered Dieudonn\'e module $(M,F^1,\phi)$ and moreover the natural homomorphism of $\dbZ_p$-algebras $\End(D)^{\text{opp}}=\End(D^{\text{t}})=\End(H^1(D))\to \End((M,F^1,\phi))$ is an isomorphism.}

\medskip
\proof
The classical Dieudonn\'e theory says that $D_k$ is determined up to unique isomorphism by $(M,\phi)$ and that the natural homomorphism $\End(D_k)^{\text{opp}}=\End(D_k^{\text{t}})\to\End((M,\phi))$ is an isomorphism (see [Dem, Ch. III, Sect. 8, Thm.], [Fo1, Ch. III, p. 128 or p. 153], etc.). Thus for $p>2$ the lemma follows from Grothendieck--Messing deformation theory (see also [Fa2, Thm. 5] or Theorem 2.3.4 below). For $p\Ge 2$ the lemma (as well as Proposition 2.2.6 below) can be deduced from [Fo1, Ch. IV, 1.6, p. 186]. As loc. cit. is stated in terms of Honda triples of the form $(M,\phi({1\over p}F^1),\phi)$, we include here a proof of the lemma for $p=2$ that appeals to Subsubsection 2.2.3. 

Suppose that $p=2$; thus $p-1=1$. As $\psi\in 2\dbG_m(V(k))$, $2$ annihilates $\text{Coker}(j_D)$. It suffices to prove the lemma under the assumption that $D_k^{\text{t}}$ is connected. Let $D_1$ be another $2$-divisible group over $\Spec\,W(k)$ that has $(M,F^1,\phi)$ as its filtered Dieudonn\'e module. As $2$ annihilates both $\text{Coker}(j_D)$ and $\text{Coker}(j_{D_1})$, via either $j_{D_1}$ or $i_{D_1}$ we can identify $H^1(D_1)$ with a $\dbZ_2$-submodule of ${1\over 2}H^1(D)$ that contains $2H^1(D)$.

We show that the assumption $H^1(D)\neq H^1(D_1)$ leads to a contradiction. We can assume $k=\bar k$. Let $D^{\acute{et}}$ and $D_1^{\acute{et}}$ be the maximal \'etale $2$-divisible groups that are quotients of $D$ and $D_1$ (respectively). Their filtered Dieudonn\'e modules are $(M^{\acute{et}},0,\phi)$, where $M^{\acute{et}}=\cap_{n\in\dbN^{\ast}} \phi^n(M)$ is the maximal direct summand of $M$ that is $W(k)$-generated by elements fixed by $\phi$. Let $r_{\acute{et}}\in\dbN$ be the rank of $M^{\acute{et}}$. Let $\{a_1,\ldots,a_r\}$ be a $W(k)$-basis for $M$ formed by elements of $F^1\cup F^0$ and such that we have $\phi(a_i)=a_i$ if $i\Le r_{\acute{et}}$; this makes sense as the canonical split cocharacter $\mu_{\text{can}}:\dbG_m\to G$ of $(M,F^1,\phi)$ fixes all elements of $M$ fixed by $\phi$ (cf. the functorial aspects of [Wi, p. 513]). As $i_{D^{\acute{et}}}$ and $i_{D_1^{\acute{et}}}$ are isomorphisms, $H^1(D^{\acute{et}})=H^1(D_1^{\acute{et}})$ is a direct summand of both $H^1(D)$ and $H^1(D_1)$. Let $x\in H^1(D)\cap H^1(D_1)$ be such that $\dbZ_2x\oplus H^1(D^{\acute{et}})$ is a direct summand of $H^1(D)\cap H^1(D_1)$ and $\dbZ_2{x\over 2}$ is a direct summand of precisely one of the $\dbZ_2$-modules $H^1(D)$ and $H^1(D_1)$. To fix the ideas we assume that ${x\over 2}\in H^1(D)\setminus H^1(D_1)$. As ${x\over 2}\beta\in i_D(M\otimes_{W(k)} B^+(W(k)))$, we can write ${x\over 2}\beta=i_D(\sum_{i=1}^r a_i\otimes \alpha_i)$, where $\alpha_i\in B^+(W(k))$. Let $\alpha_i^0\in\text{gr}^0$ be $\alpha_i$ modulo $F^1(B^+(W(k)))$. As $i_D$ is strict with respect to filtrations, we have $\alpha_i\in F^1(B^+(W(k)))$ for each $i\in\{1,\ldots,r\}$  such that $a_i\in F^0$. One can check that in fact we have $\beta\in 2B^+(W(k))$, cf. [Fo1, top of p. 79]. Either from this or from the fact that $\psi\in 2\dbG_m(V(k))$ we get that the image of ${x\over 2}\beta=x{{\beta}\over 2}\in H^1(D_1)\otimes_{\dbZ_2} B^+(W(k))[{1\over 2}]$ in $H^1(D_1)\otimes_{\dbZ_2} \text{gr}^1[{1\over 2}]$ belongs to $H^1(D_1)\otimes_{\dbZ_2} \text{gr}^1$ and generates a direct summand which when tensored with $k$ does not belong to $H^1(D^{\acute{et}})\otimes_{\dbZ_2} \text{gr}^1\otimes_{V(k)} k$. Thus there exists $i_0\in\{r_{\acute{et}}+1,\ldots,r\}$ such that:
\medskip
{\bf (i)} either $a_{i_0}\in F^1$ and $\alpha_{i_0}\in\dbG_m(B^+(W(k)))$
\smallskip
{\bf (ii)} or $a_{i_0}\in F^0$ and $j_{D_1}(a_{i_0}\otimes\alpha_{i_0}^0)$ generates a direct summand of $H^1(D_1)\otimes_{\dbZ_2} \text{gr}^1$ which when tensored with $k$ does not belong to $H^1(D_1^{\acute{et}})\otimes_{\dbZ_2} \text{gr}^1\otimes_{V(k)} k$. 

\medskip
We first assume that (i) holds. As $\Phi_k(\beta)=2\beta$, for $n\in\dbN^{\ast}$ we have $\sum_{i=1}^r \phi^n(a_i)\otimes\Phi_k^n(\alpha_i)=i_D^{-1}(2^n{x\over 2}\beta)\in 2^nM\otimes_{W(k)} B^+(W(k))$. As $\Phi_k^n(\alpha_{i_0})\in\dbG_m(B^+(W(k)))$, we get that the $W(k)$-submodule $\phi^n(M)$ of $M$ is $W(k)$-generated by $2^nM$ and $\{\phi^n(a_i)|i\in\{1,\ldots,r\}\setminus\{i_0\}\}$. Taking $n>>0$ we easily get that this implies that $(M,\phi)$ has Newton polygon slope $1$ with positive multiplicity. Thus $D_k^{\text{t}}$ is not connected. Contradiction.

We assume that (ii) holds. Due to the property enjoyed by $j_{D_1}(a_{i_0}\otimes\alpha_{i_0}^0)$, the image of $i_{D_1}(a_{i_0}\otimes 1)$ in $H^1(D_1)\otimes_{\dbZ_2} k$ does not belong to $H^1(D_1^{\acute{et}})\otimes_{\dbZ_2} k$. This statement holds also if $a_{i_0}$ is replaced by $\phi^n(a_{i_0})$. Taking $n>>0$, we get that $(M,\phi)$ has Newton polygon slope $0$ with multiplicity at least $r_{\acute{et}}+1$. This contradicts the definition of $r_{\acute{et}}$.

Thus $H^1(D)=H^1(D_1)$. Therefore $D=D_1$, cf. a classical theorem of Tate. The standard trick of replacing $D$ by $D\oplus D$ shows that to prove that $\End(H^1(D))$ surjects onto $\End((M,F^1,\phi))$, it suffices to show that each automorphism $a$ of $(M,F^1,\phi)$ is the image of an automorphism of $H^1(D)$. We identify $a$ via $i_D[{1\over {\beta}}]$ with a Galois automorphism $a^{\acute{et}}:H^1(D)[{1\over 2}]\arrowsim H^1(D)[{1\over 2}]$. An argument similar to the one above pertaining to the assumption $H^1(D)\neq H^1(D_1)$ shows that the assumption $H^1(D)\neq a^{\acute{et}}(H^1(D))$ leads to a contradiction. Thus $H^1(D)=a^{\acute{et}}(H^1(D))$ and therefore $a$ is the image of the automorphism $a^{\acute{et}}$ of $H^1(D)$.\endproof

\medskip\noindent
{\bf 2.2.6. Proposition.} {\it Suppose that $p=2$ and either $D_k$ or $D_k^{\text{t}}$ is connected. Then each direct summand of $M$ that lifts $F^1/2F^1$, is the Hodge filtration of a $2$-divisible group $D^{\prime}$ over $\Spec\,W(k)$ that lifts $D_k$. Moreover, $D^{\prime}$ is unique (up to a unique isomorphism).}

\medskip
\proof
The divided power structure of the ideal $(4)$ of $W(k)$ is nilpotent modulo $(2^n)$ for all $n\Ge 2$. Thus based on Grothendieck--Messing deformation theory, it suffices to show that there exists a unique lift $D^{\prime}_{W_2(k)}$ of $D_k$ to $\Spec\,W_2(k))$ whose Hodge filtration is a fixed direct summand of $M/4M$ that lifts $F^1/2F^1$. 

Let $d$ be the product of the dimensions of $D_k$ and $D_k^{\text{t}}$. Let $D^{\text{univ}}$ be the universal $2$-divisible group over the deformation space $\grD$ of $D_k$, cf. [Il, Cor. 4.8 (i)]. We view $\grD$ as a scheme naturally identified with $\Spec\,W(k)[[x_1,\ldots,x_d]]$ and not as a formal scheme of the form $\Spf\,W(k)[[x_1,\ldots,x_d]]$, cf. Subsection 2.1. Let $\scrR$ be the $2$-adic completion of the divided power hull of the maximal ideal of $W(k)[[x_1,\ldots,x_d]]$. Let $\delta_{2,d}$ be the natural divided power structure of the maximal ideal of $\scrR$. Let $\Phi_{\scrR}$ be the Frobenius lift of $\scrR$ that is compatible with $\sigma$ and that takes $x_i$ to $x_i^2$ for all $i\in\{1,\ldots,d\}$. The evaluation of the filtered Dieudonn\'e crystal of $D^{\text{univ}}_{\scrR}$ at the thickening $(\Spec\,k\hookrightarrow\Spec\,\scrR,\delta_{2,d})$ is defined as in Subsection 2.1 and it is $(M\otimes_{W(k)} \scrR,F^1_{\text{univ}},\phi\otimes\Phi_{\scrR},\Delta)$, where $\Delta$ is the flat connection on $M\otimes_{W(k)} \scrR$ that annihilates $M\otimes 1$ and where $F^1_{\text{univ}}$ is a direct summand of $M\otimes_{W(k)} \scrR$ that lifts $F^1/2F^1$. Let $\dbA^d_{\text{def}}$ (resp. $\dbA^d_{\text{lift}}$) be the $d$ dimensional affine space over $\Spec\,k$ that parametrizes lifts of $D_k$ to $\Spec\,W_2(k)$ (resp. lifts of $F^1/2F^1$ to direct summands of $M/4M$), the origin corresponding to $D_{W_2(k)}$ (resp. to $F^1/4F^1$). 

The existence of $F^1_{\text{univ}}$ implies that there exists a morphism $\scrN:\dbA^d_{\text{def}}\to\dbA^d_{\text{lift}}$ which on $k$-valued points takes a lift of $D_k$ to $\Spec\,W_2(k)$ to the direct summand of $M/4M$ that is the Hodge filtration of this lift. We know that $\scrN(k)$ is injective, cf. Lemma 2.2.5 and the beginning of this proof. As $m_{D_{\bar k}}$ is the pull-back of $\scrN$ to a morphism of $\Spec\,\bar k$-schemes, $\scrN(\bar k)$ is also injective. Thus the morphism $\scrN$ is quasi-finite, (due to reasons of dimensions) dominant, and generically purely inseparable. In particular, $\dbA^d_{\text{def}}$  is an open subscheme of the normalization $\dbB^d_{\text{def}}$ of $\dbA^d_{\text{lift}}$ in the field of fractions of $\dbA^d_{\text{def}}$, cf. Zariski Main Theorem (see [Ra1, Ch. IV]). We show that the assumption $\dbA^d_{\text{def}}\neq \dbB^d_{\text{def}}$ leads to a contradiction. From [Ma, Subsect. (17.H), Thm. 38] we easily get that $\dbB^d_{\text{def}}\setminus\dbA^d_{\text{def}}$ is purely of codimension 1 in $\dbB^d_{\text{def}}$ and thus contains an irreducible divisor of $\dbB^d_{\text{def}}$. As $\scrN$ is generically purely inseparable, such a divisor is the reduced scheme of the pull-back of an irreducible divisor $C$ of $\dbA^d_{\text{lift}}$. As a polynomial ring over $k$ is a unique factorization domain, $C$ is the principal divisor of a global function $f$ of $\dbA^d_{\text{lift}}$ that does not belong to $k$. We get that the polynomial ring over $k$ of global functions of $\dbA^d_{\text{def}}$ has an invertible element $f$ that does not belong to $k$. Contradiction. Therefore $\dbA^d_{\text{def}}=\dbB^d_{\text{def}}$. 

Thus the morphism $\scrN$ is finite. As $\scrN$ is also generically purely inseparable and as $\dbA_{\text{lift}}^d$ is normal, $\scrN$ is in fact a finite, purely inseparable morphism. As $k$ perfect, this implies that the map $\scrN(k)$ is surjective. As $\scrN(k)$ is also injective, the map $\scrN(k)$ is a bijection. Thus $D^{\prime}_{W_2(k)}$ exists and it is unique.\endproof

\medskip\smallskip\noindent
{\bf 2.3. Smooth Faltings--Fontaine theory.} Let $R_{00}$ be a smooth $W(k)$-algebra. Let $R_0$ be a regular, formally smooth $R_{00}$-algebra. The main examples of interest are ind-\'etale algebras over $R_{00}$ or different completions of them. We recall that an ind-\'etale algebra over $R_{00}$ is the direct limit of a functor from a filtered category (see [Mi1, App. A]) to the category of \'etale $R_{00}$-algebras. Let $R$ be either $R_0$ or $R_0^\wedge$. Let $\Phi_R$ be a Frobenius lift of $R^\wedge=R_0^\wedge$ that is compatible with $\sigma$. 

Until Subsubsection 2.3.2 we will assume that $R_0$ and $R_0/pR_0$ are integral domains; thus $R$ is also an integral domain. Let $K_{\bar R}$ be the maximal field extension of $K_R$ such that the normalization $\bar R$ of $R$ in $K_{\bar R}$ has the property that $\Spec\,\bar R[{1\over p}]$ is a pro-\'etale cover of $\Spec\,R[{1\over p}]$. The field 
extension $K_R\hookrightarrow K_{\bar R}$ is a pro-finite Galois extension (see [Gr, Exp. V, Prop. 8.2]) and thus the notations match (i.e., $K_{\bar R}$ is the field of fractions of $\bar R$). Let $A_{R/pR}:=\text{proj}.\text{lim}._{n\in\dbN^{\ast}} \bar R/p\bar R$, the transition homomorphisms being Frobenius endomorphisms. Using $A_{R/pR}$ as a substitute of $A_k$, one constructs $R$-homomorphisms $s_{R/pR}:W(A_{R/pR})\to\bar R^\wedge$ and  
$$s_{R/pR}^+:B^+(R)\to \bar R^\wedge$$ 
which are the analogous to $s_k:W(A_k)\to V(k)$ and $s_k^+:B^+(W(k))\to V(k)$ of Subsection 2.2 (if $R$ is a smooth $W(k)$-algebra, then this construction was first performed in [Fa1, Sect. 2]; see [Fa3] for other general forms of it). The $W(k)$-algebra $B^+(R)$ is equipped with a Frobenius lift $\Phi_{R/pR}$, a decreasing, exhaustive, and separated filtration $(F^n(B^+(R)))_{n\in\dbZ}$, and a $\Gal(K_{\bar R}/K_R)$-action that respects the filtration. Though $A_{R/pR}$, $s_{R/pR}$, and $\Phi_{R/pR}$ depend only on $\bar R/p\bar R$, they are not always uniquely determined but $R/pR$. However, to ease and uniformize notations, we will use the lower right index $R/pR$ and not $\bar R/p\bar R$. 

We identify $\overline{W(k)}$ with the normalization of $W(k)$ in $\bar R$. We have a canonical $W(k)$-monomorphism $B^+(W(k))\hookrightarrow B^+(R)$ that respects the Frobenius lifts and the filtrations. In particular, both $\beta$ and $\xi$ are naturally identified with elements of $F^1(B^+(R))$. 

\medskip\noindent
{\bf 2.3.1. Lemma.} {\it If $R=R_0^\wedge$ is $p$-adically complete, then the following two properties hold:

\medskip
{\bf (a)} The $R$-homomorphism $s_{R/pR}:W(A_{R/pR})\to\bar R^\wedge$ is onto.

\smallskip
{\bf (b)} The kernel of the $R$-epimorphism $s_{R/pR}:W(A_{R/pR})\twoheadrightarrow\bar R^\wedge$ is generated by $\xi$.}

\medskip
\proof
Part (a) is due in essence to Faltings. Here is a slight modification of his arguments. It is enough to show that $s_{R/pR}$ modulo $p$ is onto and thus that the Frobenius endomorphism of $\bar R/p\bar R$ is onto. We fix a $2p$-th root $p^{1\over {2p}}$ of $p$ in $V(k)$. Let $y\in \bar R$. The $\bar R$-algebra $\bar R_y:=\bar R[z]/(z^p-p^{1\over 2}z-y)$ is finite and after inverting $p$ becomes \'etale, cf. proof of [Fa3, p. 219, Lem. 5]. The $\bar R$-algebra $\bar R_y$ has a section $\bar R_y\twoheadrightarrow\bar R$, cf. the very definition of $\bar R$. Thus there exists $s(y)\in \bar R$ such that $s(y)^p-p^{1\over 2}s(y)-y=0$. A simple calculation shows that $s(y)^p-[p^{1\over {2p}}s(s(y))]^p-y\in p\bar R$. Therefore $y$ modulo $p$ has a $p$-th root in $\bar R/p\bar R$. Thus $s_{R/pR}$ is onto. The proof of (b) is the same as for the case $R=W(k)$.\endproof

\medskip\noindent
{\bf 2.3.2. Categories $\scrM\scrF_{[0,1]}^\nabla(R)$ and $\scrM\scrF_{[0,1]}(R)$.} In this subsubsection we do not assume that $R_0$ and $R_0/pR_0$ are integral domains. Let $\Omega_{R^\wedge/W(k)}$ be a suitable completion of the $R^\wedge$-module of relative differentials (e.g., if $R_0$ is a smooth $W(k)$-algebra, then it is the $p$-adic completion and if $R_0=W(k)[[z_1,\ldots,z_m]]$, then it is the $(z_1,\ldots,z_m)$-adic completion). Let $d\Phi_{R}/p:\Omega_{R^\wedge/W(k)}\to \Omega_{R^\wedge/W(k)}$ be the differential of $\Phi_R$ divided by $p$. Let $\scrM\scrF_{[0,1]}^\nabla(R)$ be the Faltings--Fontaine $\dbZ_p$-linear category defined as follows (see [FL, Sect. 1] and [Wi, Sect. 1] for $R=W(k)$, [Fa1, Sect. 2] for $R$ smooth, and [Va2, Sect. 2] in general). Its objects are quintuples $(N,F,\phi_0,\phi_1,\nabla)$,
where $N$ is a finitely generated $p$-torsion $R$-module, $F$ is a direct summand of $N$, $\phi_0:N\to N$ and $\phi_1:F\to N$ are $\Phi_R$-linear maps, and $\nabla:N\to N\otimes_R \Omega_{R^\wedge/W(k)}$ is an integrable, nilpotent modulo $p$ connection on $N$, such that the below five axioms hold:

\medskip
{\bf (a)} we have $\phi_0(x)=p\phi_1(x)$ for all $x\in F$;
\smallskip
{\bf (b)} the $R$-module $N$ is $R$-generated by $\phi_0(N)+\phi_1(F)$;
\smallskip
{\bf (c)} we have $\nabla\circ\phi_0(x)=p(\phi_0\otimes d\Phi_{R}/p)\circ\nabla(x)$ for all $x\in N$;
\smallskip
{\bf (d)} we have $\nabla\circ\phi_1(x)=(\phi_0\otimes d\Phi_{R}/p)\circ\nabla(x)$ for all $x\in F$;
\smallskip
{\bf (e)} locally in the Zariski topology of $\Spec\,R/pR$, $N$ is a finite direct sum of $R$-modules of the form $R/p^sR$, where $s\in\dbN$.

\medskip
A morphism $f:(N,F,\phi_0,\phi_1,\nabla)\to (N^\prime,F^\prime,\phi_0^\prime,\phi_1^\prime,\nabla^\prime)$ between two such quintuples is an $R$-linear map $f:N\to N^\prime$ such that we have an inclusion $f(F)\subseteq F^\prime$ and identities $\phi_0^\prime\circ f=f\circ\phi_0$, $\phi_1^\prime\circ f=f\circ\phi_1$, and $\nabla^\prime\circ f=(f\otimes_R 1_{\Omega_{R/W(k)}})\circ\nabla$. Additions of $f$'s and multiplications of $f$'s by elements of $\dbZ_p$ define the $\dbZ_p$-linear structure of $\scrM\scrF_{[0,1]}^\nabla(R)$. Due to (e) we can identify $\scrM\scrF_{[0,1]}^\nabla(R^\wedge)=\scrM\scrF_{[0,1]}^\nabla(R)$. Suppressing the connections and thus axioms (c) and (d), we get the $\dbZ_p$-linear category $\scrM\scrF_{[0,1]}(R^\wedge)=\scrM\scrF_{[0,1]}(R)$.  

\medskip\noindent
{\bf 2.3.3. The functor $\dbD$.} Let $p-FF(\Spec\,R)$ be the $\dbZ_p$-linear category of finite, flat, commutative group schemes over $\Spec\,R$ of $p$-power rank (order). We recall from [Va2, Sect. 2] that there exists a contravariant $\dbZ_p$-linear functor 
$$\dbD:p-FF(\Spec\,R)\to\scrM\scrF_{[0,1]}^\nabla(R)=\scrM\scrF_{[0,1]}^\nabla(R^\wedge).$$ 
\indent
Let $H$ be an object of $p-FF(\Spec\,R)$. We review from loc. cit. the construction of $\dbD(H)=(N,F,\phi_0,\phi_1,\nabla)$. Let $n\in\dbN^{\ast}$ be such that $p^n$ annihilates $H$. The triple $(N,\phi_0,\nabla)$ is part of the evaluation $(N,\phi_0,\vartheta_0,\nabla)$ of the Dieudonn\'e crystal $\dbD(H_{R/pR})$ at the thickening $(\Spec\,R/pR\hookrightarrow\Spec\,R/p^nR,\delta(p))$ and this motivates the notation $\dbD(H)$. The direct summand $F$ of $N$ is the Hodge filtration of $N$ defined by the lift $H_{R/p^nR}$. To define $\phi_1$ we can work locally in the Zariski topology of $\Spec\,R/pR$ and thus we can assume  that $R=R^\wedge=R_0^\wedge$ and (cf. Raynaud theorem of [BBM, Thm. 3.1.1]) that $H$ is a closed subgroup scheme of an abelian scheme $\scrA$ over $\Spec\,R$. Let $N_\scrA:=H^1_{\text{dR}}(\scrA/R)=H^1_{\text{crys}}(\scrA/R)$. Let $F_\scrA$  be the direct summand of $N_\scrA$ that is the Hodge filtration defined by $\scrA$. Let $\nabla_\scrA$ be the Gauss--Manin connection on $N_\scrA$. Let $\phi_{0\scrA}$ be the $\Phi_R$-linear endomorphism of $N_\scrA$. The embedding $H\hookrightarrow \scrA$ defines a canonical $R$-epimorphism $(N_\scrA,F_\scrA,\phi_{0\scrA},\nabla_\scrA)\twoheadrightarrow (N,F,\phi_0,\nabla)$, cf. [BBM, Thm. 3.1.2 and pp. 132--133]. Thus as $\phi_1:F\to N$ we can take the natural quotient of the restriction of ${1\over p}\phi_{0\scrA}$ to $F_\scrA$. 

As the $\dbZ_p$-linear categories of finite, flat, commutative group schemes of $p$-power rank over $\Spec\,R_0^\wedge$ and over $\Spf\,R_0^\wedge$ are canonically isomorphic, from [Fa1, Thm. 7.1] we get that:

\medskip\noindent
{\bf 2.3.4. Theorem.} {\it Suppose that $p>2$, that the $W(k)$-algebra $R_0$ is smooth, and that $R=R_0^\wedge$. Then the $\dbZ_p$-linear functor $\dbD:p-FF(\Spec\,R)\to\scrM\scrF_{[0,1]}^\nabla(R)$ is an antiequivalence of $\dbZ_p$-linear categories.}

\medskip
One can use [Fa1, Sect. 2, pp. 31--33] and [Ra2, Thm. 3.3.3] to check that the above $\dbZ_p$-linear categories are abelian and that the functor $\dbD$ is in fact an antiequivalence of abelian categories. Theorem 2.3.4 will be used in the proof of Theorem 3.4.1 that translates Theorem 3.2 in terms of $p$-divisible groups.

\medskip\noindent 
{\bf 2.3.5. A more general form of the monomorphism $i_D$.} In this subsubsection we assume that $R_0$ and $R_0/pR_0$ are integral domains, that $R=R_0^\wedge$, and that there exist elements $z_1,\ldots,z_m\in R$ such that $R$ is a formally \'etale $W(k)[z_1,\ldots,z_m]$-algebra and the equality $\Phi_R(z_i)=z_i^p$ holds for all $i\in\{1,\ldots,m\}$. For instance, $R$ could be $W(k)[[z_1,\ldots,z_m]]$. For $n\in\dbN^{\ast}$, each element of $\Ker(B^+(R)/p^nB^+(R)\to \bar R^\wedge/p^n\bar R^\wedge)$ is nilpotent. Thus the formally \'etale assumption and the fact that $R$ is $p$-adically complete imply that there exists a unique $W(k)$-monomorphism $i_R:R\hookrightarrow B^+(R)$ that lifts the inclusion $R\hookrightarrow \bar R^\wedge$ and that takes $z_i$ to $(\underline{z_i},0,0,\ldots)\in W(A_{R/pR})\hookrightarrow B^+(R)$. The $p$-power roots of $z_i$ might not exist in $\bar R$ (as we have not assumed that $z_i$ is a unit of $R$) but as the Frobenius endomorphism of $\bar R/p\bar R$ is surjective (see proof of Lemma 2.3.1 (a)), one can similarly define an element $\underline{z_i}\in A_{R/pR}$ using $p$-power roots of $z_i$ modulo $p$. 

Let $\scrD$ be a $p$-divisible group over $\Spec\,R$. Let $(N,F^1_N,\phi_N,\nabla_N)$ be as in Subsection 2.1. As the $R$-homomorphism $s_{R/pR}^+:B^+(R)\to \bar R^\wedge$ is onto (cf. Lemma 2.3.1 (a)), there exists a $B^+(R)$-linear monomorphism
$$i_{\scrD}:N\otimes_{R} B^+(R)\hookrightarrow H^1({\scrD}_{K_{R}})\otimes_{\dbZ_p} B^+(R)\leqno (3)$$
that is constructed and has the same properties as the $B^+(W(k))$-linear monomorphism $i_D$ of Subsubsection 2.2.2. Here, similarly to $H^1(D)$, $H^1({\scrD}_{K_{R}})$ is the dual of the Tate-module of $\scrD_{K_R}$. As the $W(k)$-monomorphism $i_R:R\hookrightarrow W^+(R)$ respects Frobenius lifts, the Frobenius endomorphism of $N\otimes_{R} B^+(R)$ is $\phi_N\otimes \Phi_{R/pR}$. 

\medskip\smallskip\noindent 
{\bf 2.4. Artin--Schreier systems of equations.} \rm
For a matrix $\dag$, let $\dag^{[p]}$ be the matrix obtained by raising the entries of $\dag$ to their $p$-powers. Let $n\in\dbN^{\ast}$. Let $x_1,\ldots,x_n$ be variables. Let $x:=(x_1,\ldots,x_n)^{\text{t}}$. Let $\Spec\,A_0$ be an affine $\Spec\,k$-scheme. An {\it Artin--Schreier system of equations} in $n$ variables over $A_0$ is a system
$$x=Bx^{[p]}+C_0,\; \text{where}\; B=(b_{ij})_{1\le i,j\le n}\in M_{n\times n}(A_0)\;\text{and}\; C_0=(c_1,\ldots,c_n)^{\text{t}}\in M_{n\times 1}(A_0).\leqno (4)$$
Let $A$ be the finitely generated $k$-subalgebra of $A_0$ generated by the entries of $B$. Let $A_1:=A_0[x_1,\ldots,x_n]/(x_1-c_1-\sum_{j=1}^n b_{1j}x_j^p,\ldots,x_n-c_n-\sum_{j=1}^n b_{nj}x_j^p)$ be the $A_0$-algebra defined by (4) i.e., the $A_0$-algebra that parametrizes solutions of (4). For $l\in\dbN^{\ast}$ with $l\le 2$, we construct inductively an $A_{l-1}$-algebra $A_l$ that parametrizes solutions of the Artin--Schreier system of equations $x=Bx^{[p]}+C_{l-1}$, where $C_{l-1}\in M_{n\times 1}(A_{l-1})$. Let $A_{\infty}$ be the inductive limit of the $k$-algebras $A_l$ with $l\in\dbN$.  Each projective system of affine $k$-schemes of the form 
$$\cdots\to \Spec\,A_n\to \Spec\,A_{n-1}\to\cdots\to\Spec\,A_0$$ 
will be called an {\it Artin--Schreier tower}. 
Let $\Spec\,\Lambda$ be the group scheme over $\Spec\,A$ defined by the addition of solutions of the Artin--Schreier system of equations $x=Bx^{[p]}$; we call it an {\it Artin--Schreier group scheme}.

Next theorem plays a key role in the proof and applications of Theorem 3.2 below. 

\medskip\noindent
{\bf 2.4.1. Theorem.} {\it The following five properties hold:

\medskip
{\bf (a)} The affine morphism $\Spec\,A_1\to \Spec\,A_0$ is \'etale. 

\smallskip
{\bf (b)} Each geometric fibre of the morphism $\Spec\,A_1\to \Spec\,A_0$ has a number of geometric points equal to $p^m$, where $m\in\{0,\ldots,n\}$ depends on the fibre. 

\smallskip
{\bf (c)} There exists a stratification of $\Spec\,A$ in reduced, locally closed subschemes $U^{(1)}$, $\ldots,U^{(s)}$ defined inductively by the following property holds: for each $i\in\{1,2,\ldots,s\}$ the scheme $U^{(i)}$ is the maximal open subscheme of the reduced scheme of $\Spec\,A\setminus (\cup_{j=1}^{i-1} U^{(j)})$ which has the property that the morphism $\Spec\,\Lambda\times_{\Spec\,A} U^{(i)} \to U^{(i)}$ is an \'etale cover. Moreover, for all $l\in\dbN^{\ast}$ and $i\in\{1,\ldots,s\}$, if $U_l^{(i)}:=\Spec\,A_l\times_{\Spec\,A} U^{(i)}$ and $U_{l-1}^{(i)}:=\Spec\,A_{l-1}\times_{\Spec\,A} U^{(i)}$, the morphism $U_l^{(i)}\to U_{l-1}^{(i)}$ is an \'etale cover.

\smallskip
{\bf (d)} If $A_0$ is noetherian, then the image of each connected component $\scrC_{\infty}$ of $\Spec\,A_{\infty}$ in $\Spec\,A_0$ is an open subscheme of a connected component of $\Spec\,A_0$.

\smallskip
{\bf (e)} Suppose that $A_0$ is a finitely generated $k$-algebra. Then each point $\grp$ of $\Spec\,A_{\infty}$ specializes to a point of $\Spec\,A_{\infty}$ whose residue field is algebraic over $k$.}

\medskip
\proof
We prove (a). The $A_0$-algebra $A_1$ is of finite presentation.  We have $\Omega_{A_1/A_0}=0$, cf. the shape of the system of equations $(4)$. We check that the criterion of formal \'etaleness holds for $\Spec\,A_1\to\Spec\,A_0$. Let $S$ be a commutative $A_0$-algebra and let $I$ be an ideal of $S$ of such that $I^2=0$. Let $\bar w$ be a solution of (4) in $S/I$. If $w\in M_{n\times 1}(S)$ lifts $\bar w$, then $Bw^{[p]}+C_0$ does not depend on the choice of $w$ and it is the unique solution of (4) in $S$ that lifts $\bar w$. We conclude that the morphism $\Spec\,A_1\to\Spec\,A_0$ is \'etale, cf. [BLR, Ch. 2, Sect. 2.2, Props. 2 and 6]. 

To prove (b) we can assume that $A_0=k=\bar k$ and we have to show that $A_1=k^{p^m}$ (i.e., $(4)$ has $p^m$ solutions in $k$) with $m\in\{0,\ldots,n\}$. We will include three proofs of this. 

The first proof of (b) is due to Ofer Gabber. We identify $\dbG_a^n(k)=M_{n\times 1}(k)$. Let $e_B:\dbG_a^n\to \dbG_a^n$ be the endomorphism which at the level of $k$-valued points (viewed as columns) maps $x=(x_1,\ldots,x_n)^t$ to $x-Bx^{[p]}$. As $e_B$ is \'etale (cf. either proof of (a) or [Gr, Exp. II, Cor. 4.6]), the image of $e_B$ is an open subgroup of $\dbG_a^n$. As $\dbG_a^n$ is connected, we get that $e_B$ is an \'etale epimorphism. Thus for each $C_0\in \dbG_a^n(k)$, $(4)$ has as many solutions in $k$ as $\Ker(e_B)(k)$. But $\Ker(e_B)$ is an \'etale group over $k$ annihilated by $p$ and thus it is isomorphic to $(\dbZ/p\dbZ)^m$ with $m\ge 0$. Therefore the system of equations $(4)$ has $p^m$ solutions in $k$, where $m\ge 0$ does not depend on $C_0$. Each point of $\Spec\,A_1$ is an irreducible component of the intersection of those $n$ hypersurfaces in the projective space $\dbP^n_k$ over $\Spec\,k$ that are of degree $p$ and that are defined by the rows of the homogeneous system of equations $x_0^{p-1}x=Bx^{[p]}+x_0^pC_0$ in the variables $x_0,x_1,\ldots,x_n$. Here the multiplication by $x_0$ is a scalar multiplication. From the B\'ezout inequality we get that $p^m$ is at the most the product $p^n$ of the degrees of these $n$ hypersurfaces, cf. [Fu, Ex. 8.4.6]. Thus $m\le n$. 

The second proof of (b) was suggested by a referee and goes as follows. We consider the abelian category of pairs $(\bar M,\bar\phi)$, where $\bar M$ is a finite dimensional $k$-vector space and $\bar\phi:\bar M\to\bar M$ is a $\sigma$-linear endomorphism. Each object of this category has finite length and the simple objects are $1$-dimensional and isomorphic to either $(k,0)$ or $(k,\sigma)$ (see [GKHR, Subsect. 5.1]). Based on this applied to $\bar M=k^n$ and to its $\sigma$-linear endomorphism $\bar\phi$ whose matrix representation with respect to the standard $k$-basis for $k^n$ is $B$, it is easy to see that by using a substitution of the form $x=\scrB y$ with $\scrB\in \pmb{GL}_n(k)$ and $y=(y_1,\ldots,y_n)^t$ and thus by replacing $(B,C_0)$ with $(\scrB^{-1}B\scrB^{[p]},\scrB^{-1}C_0)$, to solve $(4)$ we can assume that the matrix $B\in M_{n\times n}(k)$ is upper triangular. Thus the system of equations $(4)$ becomes
$$x_1=b_{11}x_1^p+\cdots +b_{1n}x_n^p+c_1$$
$$\cdots=\cdots$$
$$x_{n-1}=b_{n-1 n-1}x_{n-1}^p+b_{n-1 n}x_n^p+c_{n-1}$$
$$x_n=b_{nn}x_n^p+c_n.$$
Starting with $x_n$ and working backwards, we see that the number of solutions in $k$ of $(4)$ is $p^m$, where $m\in\{0,\ldots,n\}$ is the number of non-zero entries of the diagonal of $B$.    

Here is the third (original) proof of (b). We use induction on $n\in\dbN^{\ast}$. The case $n=1$ is trivial. For $n\Ge 2$ the passage from $n-1$ to $n$ goes as follows. If $B$ is an invertible matrix, then $(4)$ is equivalent to the system of equations $x^{[p]}=B^{-1}x-B^{-1}C_0$ which defines a finite $k$-algebra of rank (order) $p^n$; as $A_1$ is an \'etale $k$-algebra (cf. (a)), we have $A_1=k^{p^n}$.  

Suppose now that the matrix $B$ is not invertible. Up to a renumbering of $x_i$'s we can assume that we have an identity $\sum_{i=1}^{n} d_iB^{(i)}=0$, where $d_i\in k$ with $d_n=1$ and where $B^{(i)}$ is the $i$-th row of $B$. The system of equations $(4)$ is equivalent to the system $\scrX$ of equations, where we keep the first $n-1$ equations of $(4)$ which involve $x_1,\ldots,x_{n-1}$ on the left hand side and where we replace the last equation of $(4)$ which involves $x_n$ on the left hand side by the equation
$$x_n+\sum_{i=1}^{n-1} d_ix_i=c_n+\sum_{i=1}^{n-1} c_id_i.\leqno (5)$$
The linear equation (5) allows us to eliminate the variable $x_n$ in $\scrX$ (and thus also in (4)). By performing this elimination we come across an Artin--Schreier system of equations in $n-1$ variables over $k$ of the form $x^{\prime}=B^{\prime}x^{\prime[p]}+C_0^{\prime}$, where $x^{\prime}:=(x_1^{\prime},\ldots,x_{n-1}^{\prime})^{\text{t}}$. The entries of $B^{\prime}\in M_{n-1\times n-1}(k)$ depend on $B$ but not on $C_0$. By induction we have $A_1=k^{p^m}$ for some $m\in\{0,\ldots,n-1\}$. This proves (b).
 
To prove (c), we note that the group scheme $\Spec\,\Lambda$ over $\Spec\,A$ is \'etale (cf. (a)). As $\Spec\,\Lambda\to\Spec\,A$ is an \'etale morphism between noetherian schemes, the existence of the stratification $U^{(1)}$, $\ldots,U^{(s)}$ of $\Spec\,A$ is a standard piece of algebraic geometry. As the fibres of $\Spec\,A_l\to\Spec\,A_{l-1}$ are non-empty (cf. (b)), the product rule $(y,x)\to y+x$ makes $\Spec\,A_l$ to be a left torsor of the affine group scheme $\Spec\,\Lambda\times_{\Spec\,A} \Spec\,A_{l-1}$ over $\Spec\,A_{l-1}$. Thus $U^{(i)}_l$ is an \'etale cover of $U^{(i)}_{l-1}$ i.e., (c) holds. 

To prove (d) and (e) we can assume that $\Spec\,A_0$ is reduced and connected. We prove (d). As $U^{(i)}_l$ is an \'etale cover of $U^{(i)}_{l-1}$ (cf. (c)), the pull-back $U_{\infty}^{(i)}$ of $U^{(i)}$ to $\Spec\,A_{\infty}$ is a pro-\'etale cover of $U^{(i)}_0$. The going up (resp. down) property holds for pro-\'etale covers (resp. pro-\'etale morphisms), cf. [Ma, Subsect. (5.E), Thm. 5] (resp. [Ma, Subsect. (5.D), Thm. 4]). We apply this to the pro-\'etale covers $U_{\infty}^{(i)}\to U_{0}^{(i)}$ and to the pro-\'etale morphism $\Spec\,A_{\infty}\to\Spec\,A_0$. We get that the image $\scrC_0$ of $\scrC_{\infty}$ in $\Spec\,A_0$ is the complement in $\Spec\,A_0$ of the union of the schematic closures in $\Spec\,A_0$ of those connected components of some $U_{0}^{(i)}$  with $i\in\{2,\ldots,s\}$ that do not intersect $\scrC_0$. As $\Spec\,A_0$ is noetherian, this union is finite and thus $\scrC_0$ is an open subscheme of $\Spec\,A_0$. If $\Spec\,A_0$ is not irreducible, then in general $\scrC_0$ is not Zariski dense in it. Here is an example communicated to us by Ofer Gabber. If $A_0=k[u,v]/(uv)$, $n=1$, and $A_1=A_0[x_1]/(ux_1^p-x_1)$, then $\scrC_0$ could be also the open subscheme $\Spec\,A_0\setminus \Spec\,k[v]$ over which $u$ is invertible. 

We prove (e). To prove (e), for $*\in\dbN\cup\{\infty\}$ we can replace $A_*$ with $A_*/(\grp\cap A_0)A_*$. Thus we can assume that $\grp\cap A_0=0$ and that $A_0$ is an integral domain. By localizing we can also assume that $U^{(1)}$ is $\Spec\,A$. The morphism $\Spec\,A_{\infty}/\grp\to\Spec\,A_0$ is an \'etale cover and thus surjective. As $\Spec\,A_0$ has points whose residue fields are algebraic over $k$ and due to the going down property, we get that $\grp$ specializes to points of $\Spec\,A_{\infty}$ whose residue fields are algebraic over $k$.\endproof 

\medskip\noindent
{\bf 2.4.2. Definition.} Either the stratification of $\Spec\,A$ in reduced, locally closed subschemes of Theorem 2.4.1 (c) or its pull-back to a stratification of $\Spec\,A_0$ in reduced, locally closed subschemes, will be called an {\it Artin--Schreier stratification}. 

\medskip\noindent
{\bf 2.4.3. Remark.} The case $C_0=0$ (i.e., the $A$-algebra $\Lambda$ of the proof of Theorem 2.4.1 (c)) is studied systematically for the first time in [Ly], using a language of $F$-modules. We came across Artin--Schreier systems of equations independently in 1998 in connection to the proof of Theorem 3.2 below.  

\medskip\smallskip\noindent
{\bf 2.5. Simple group properties.} 
Let $S\in\{W(k),B^+(W(k))\}$. Let $\iota\in S$ be $p$ if $S=W(k)$ and be $\beta$ if $S=B^+(W(k))$. Let $(t_{\alpha})_{\alpha\in\tilde\scrJ}$ be another family of tensors of $F^0(\scrT(M))[{1\over p}]$ fixed by $\phi$ and such that $G_{B(k)}$ is also the subgroup of $\pmb{GL}_{M[{1\over p}]}$ that fixes $t_{\alpha}$ for all $\alpha\in\tilde\scrJ$. For $\alpha\in\tilde\scrJ$ let $v_{\alpha}\in\scrT(H^1(D)[{1\over p}])$ be such that $\scrT(i_D[{1\over \beta}])$ maps $t_{\alpha}$ to $v_{\alpha}$, cf. Subsubsection 2.2.4. 

\medskip\noindent
{\bf 2.5.1. Fact.} {\it {\bf (a)} There exists an isomorphism $\tilde\rho_S:(M\otimes_{W(k)} S,(t_{\alpha})_{\alpha\in\tilde\scrJ})\arrowsim (H^1(D)\otimes_{\dbZ_p} S,(v_{\alpha})_{\alpha\in\tilde\scrJ})$ if and only if there exists an isomorphism $\rho_S:(M\otimes_{W(k)} S,(t_{\alpha})_{\alpha\in\scrJ})\arrowsim (H^1(D)\otimes_{\dbZ_p} S,(v_{\alpha})_{\alpha\in\scrJ})$.

\smallskip
{\bf (b)} Suppose there exists an isomorphism $\varrho_S:(M\otimes_{W(k)} S[{1\over {\iota}}],(t_{\alpha})_{\alpha\in\scrJ})\arrowsim (H^1(D)\otimes_{\dbZ_p} S[{1\over {\iota}}],(v_{\alpha})_{\alpha\in\scrJ})$. Then there exists an isomorphism $\rho_S:(M\otimes_{W(k)} S,(t_{\alpha})_{\alpha\in\scrJ})\arrowsim (H^1(D)\otimes_{\dbZ_p} S,(v_{\alpha})_{\alpha\in\scrJ})$ if and only if there exists an element $h\in G(S[{1\over {\iota}}])$ that takes the $S$-submodule $\varrho_S^{-1}(H^1(D)\otimes_{\dbZ_p} S)$ of $M\otimes_{W(k)} S[{1\over {\iota}}]$ onto $M\otimes_{W(k)} S$.}

\medskip
\proof
We prove (a). If $\tilde\rho_S$ exists, then its tensorization with $B^+(W(k))[{1\over {\beta}}]$ is of the form $i_D[{1\over {\beta}}]\circ h$, where $h\in \pmb{GL}_M(B^+(W(k))[{1\over {\beta}}])$ fixes $t_{\alpha}$ for all $\alpha\in\tilde\scrJ$. Thus $h\in G(B^+(W(k))[{1\over {\beta}}])$ and therefore $\tilde\rho_S$ takes $t_{\alpha}$ to $v_{\alpha}$ for all $\alpha\in\scrJ$. Thus we can take $\rho_S$ to be defined by $\tilde\rho_S$. Similarly, if $\rho_S$ exists, then we can take $\tilde\rho_S$ to be defined by $\rho_S$. 

Part (b) is an elementary exercise: one regains each element of the set $\{\rho_S,h\}$ from the other element via the identity $\rho_S[{1\over {\iota}}]=\varrho_S\circ h^{-1}$.\endproof

\medskip
Let $\sigma_{\phi}:=\phi\mu(p)$. It is a $\sigma$-linear automorphism of $M$ that normalizes $\Lie(G_{B(k)})\cap\End_{W(k)}(M)$. Until Subsection 2.6 we assume $k=\bar k$. Thus $M_{\dbZ_p}:=\{x\in M|\sigma_{\phi}(x)=x\}$ is a $\dbZ_p$-structure of $M$. As $\mu(p)$ and $\phi$ fix $t_{\alpha}$, we have $\sigma_{\phi}(t_{\alpha})=t_{\alpha}$. Thus we have $t_{\alpha}\in\scrT(M_{\dbZ_p}[{1\over p}])$ for all $\alpha\in\scrJ$. Thus $G$ is the pull-back to $\Spec\,W(k)$ of a flat, closed subgroup scheme $G_{\dbZ_p}$ of $\pmb{GL}_{M_{\dbZ_p}}$. If $G$ is smooth, then $\Lie(G_{\dbZ_p})=\{x\in\Lie(G)|\sigma_{\phi}(x)=x\}$. 

\medskip\noindent
{\bf 2.5.2. Lemma.} {\it {\bf (a)} Suppose that $G$ is smooth and $G_k$ is connected. Then there exists an isomorphism 
$\rho_{\dbZ_p}:(M_{\dbZ_p},(t_{\alpha})_{\alpha\in\scrJ})\arrowsim (H^1(D),(v_{\alpha})_{\alpha\in\scrJ})$ if and only if there exists an isomorphism $\rho:(M,(t_{\alpha})_{\alpha\in\scrJ})\arrowsim (H^1(D)\otimes_{\dbZ_p} W(k),(v_{\alpha})_{\alpha\in\scrJ})$.
\smallskip
{\bf (b)} Suppose that $G$ is smooth. Then there exists an isomorphism $\rho$ as in (a) if and only if there exists an isomorphism $\rho_{B^+(W(k))}:(M\otimes_{W(k)} B^+(W(k)),(t_{\alpha})_{\alpha\in\scrJ})\arrowsim (H^1(D)\otimes_{\dbZ_p} B^+(W(k)),(v_{\alpha})_{\alpha\in\scrJ})$.}

\medskip
\proof
Let $Y_{\dbQ_p}$ be the affine $\Spec\,\dbQ_p$-scheme that parametrizes isomorphisms between $(M_{\dbZ_p}[{1\over p}],(t_{\alpha})_{\alpha\in\scrJ})$ and $(H^1(D)[{1\over p}],(v_{\alpha})_{\alpha\in\scrJ})$. Let $Y$ be the schematic closure of $Y_{\dbQ_p}$ in the smooth $\dbZ_p$-scheme $Isom(M_{\dbZ_p},H^1(D))$ that parametrizes linear isomorphisms. The group scheme $G_{\dbZ_p}$ acts on $Y$ from the left. The only if part of (a) is trivial. We prove the if part of (a). As $\rho$ exists, $Y_{W(k)}$ is the trivial torsor of $G$. Thus $Y$ is a torsor of $G_{\dbZ_p}$ in the flat topology and therefore it is a smooth $\Spec\,\dbZ_p$-scheme of finite type. Thus $Y$ has points with values in finite fields. As $G_{\dbF_p}$ is connected, $Y$ has $\dbF_p$-valued points (cf. Lang's theorem). Thus $Y$ has $\dbZ_p$-valued points i.e., $\rho_{\dbZ_p}$ exists. This proves (a). The proof of (b) is similar, starting from the fact that each torsor of $G$ in the flat topology of $\Spec\,W(k)$ has $k$-valued points and therefore (as $G$ is smooth) it also has $W(k)$-valued points.\endproof 

\medskip
The next lemma is used in Subsection 2.7 which at its turn is used in the proof of Theorem 3.5.1 that reduces the proof of the Main Theorem to the case when $G$ is a torus.

\medskip\noindent
{\bf 2.5.3. Lemma.} {\it We recall that $k=\bar k$. Let $(M_{\dbZ_p},G_{\dbZ_p})$ and $\sigma_{\phi}=\phi\mu({1\over p})$ be as above. Let $\tilde G$ be a flat, closed subgroup scheme of $G$ such that the following three properties hold:

\medskip
{\bf (i)} the cocharacter $\mu:\dbG_m\to G$ factors through $\tilde G$;

\smallskip
{\bf (ii)}  we have $\phi(\Lie(\tilde G_{B(k)}))=\Lie(\tilde G_{B(k)})$;

\smallskip
{\bf (iii)} the generic fibre $\tilde G_{B(k)}$ is connected.

\medskip
Then there exists a unique closed, flat subgroup scheme $\tilde G_{\dbZ_p}$ of $G_{\dbZ_p}$ whose pull-back to $\Spec\,W(k)$ is $\tilde G$. If $\tilde G_{B(k)}$ is a reductive group, then the generic fibre $\tilde G_{\dbQ_p}$ of $\tilde G_{\dbZ_p}$ is also a reductive group and therefore there exists a set $\tilde\scrJ$ that contains $\scrJ$ and such that there exist a family of tensors $(t_{\alpha})_{\alpha\in\tilde\scrJ}$ of $\scrT(M_{\dbZ_p}[{1\over p}])$ with the property that $\tilde G_{\dbQ_p}$ is the subgroup of $\pmb{GL}_{M_{\dbZ_p}[{1\over p}]}$ which fixes $t_{\alpha}$ for all $\alpha\in\tilde\scrJ$ (thus each $t_{\alpha}$ is fixed by $\phi$).}

\medskip
\proof
We know that $\sigma_{\phi}$ normalizes $\Lie(\tilde G_{B(k)})$, cf. properties (i) and (ii). The following Lie algebra $\tilde\grg_{\dbQ_p}:=\{x\in \Lie(\tilde G_{B(k)})|\sigma_{\phi}(x)=x\}$ over $\dbQ_p$ is the unique Lie subalgebra of $\Lie(G_{\dbQ_p})$ such that we have $\Lie(\tilde G_{B(k)})=\tilde\grg_{\dbQ_p}\otimes_{\dbQ_p} B(k)$. We check that there exists a unique connected subgroup $\tilde G_{\dbQ_p}$ of $G_{\dbQ_p}$ whose Lie algebra is $\tilde\grg_{\dbQ_p}$.  The uniqueness part is implied by [Bo, Ch. II, Subsect. 7.1]. To check the existence part, we consider commutative $\dbQ_p$-algebras $A$ such that there exists a closed, flat subgroup scheme $\tilde G_{A}$ of $G_{A}$ whose Lie algebra is $\tilde\grg_{\dbQ_p}\otimes_{\dbQ_p} A$. For instance, $A$ can be $B(k)$ itself and thus we can assume $A$ is a finitely generated $\dbQ_p$-subalgebra of $B(k)$. By replacing $A$ with its quotient through a maximal ideal, we can assume $A$ is a finite field extension of $\dbQ_p$. Even more, we can assume that $A$ is a finite Galois extension of $\dbQ_p$ and that $\tilde G_{A}$ is connected. As $\Lie(\tilde G_{A})=\tilde\grg_{\dbQ_p}\otimes_{\dbQ_p} A$, from [Bo, Ch. II, Subsect. 7.1] we get that the natural action of the Galois group $\Gal(A/\dbQ_p)$ on $\tilde\grg_{\dbQ_p}\otimes_{\dbQ_p} A$ is defined naturally by a natural action of $\Gal(A/\dbQ_p)$ on $\tilde G_{A}$. This last action is free. As $\tilde G_{A}$ is an affine scheme, the quotient $\tilde G_{\dbQ_p}$ of $\tilde G_{A}$ by $\Gal(A/\dbQ_p)$ exists (cf. [BLR, Ch. 6, Sect. 6.1, Thm. 5]) and it has all the desired properties.

From [Bo, Ch. II, Subsect. 7.1] and the property (iii) we get that our notations match i.e., the pull-back of the connected group $\tilde G_{\dbQ_p}$ to $\Spec\,B(k)$ has the same Lie algebra as $\tilde G_{B(k)}$ and therefore it is $\tilde G_{B(k)}$. Let $\tilde G_{\dbZ_p}$ be the schematic closure of $\tilde G_{\dbQ_p}$ in $G_{\dbZ_p}$; it is the unique closed, flat subgroup scheme of $G_{\dbZ_p}$ whose pull-back to $\Spec\,W(k)$ is $\tilde G$. Obviously $\tilde G_{\dbQ_p}$ is reductive if and only if $\tilde G_{B(k)}$ is so. The last part of the lemma that pertains to the existence of the set $\tilde\scrJ$ follows from [De2, Prop. 3.1 c)]. Note that each $t_{\alpha}$ with $\alpha\in\tilde\scrJ$ is fixed by $\sigma_{\phi}$ (as $t_{\alpha}\in \scrT(M_{\dbZ_p}[{1\over p}])$) and by $\mu(p)$ (cf. (i) and the fact that $\tilde G_{\dbQ_p}$ fixes $t_{\alpha}$). Thus each $t_{\alpha}$ with $\alpha\in\tilde\scrJ$ is fixed by $\phi=\sigma_{\phi} \mu(p)$.\endproof

\medskip\noindent
{\bf 2.5.4. Fact.} {\it The following two statements are equivalent:

\medskip
{\bf (i)} The Main Theorem holds whenever $D$ is connected.

\smallskip
{\bf (ii)} The Main Theorem holds whenever the property (EC) holds for $D$.}

\medskip
\proof
If $D$ is connected, then the property (EC) holds for it (cf. Definition 1.1.1) and therefore (ii) implies (i). We check that (i) implies (ii). Thus we have direct sum decompositions $D=D_1\oplus D_2$ and $(M,F^1,\phi)=(M_1,F^1(M_1),\phi)\oplus (M_2,F^1(M_2),\phi)$ with the property that the two $p$-divisible groups $D_1$ and $D_2^{\text{t}}$ are connected. We consider the connected $p$-divisible group $\tilde D:=D_1\oplus D_2^{\text{t}}\oplus \mu_{p^{\infty}}$ and its filtered Dieudonn\'e module $(\tilde M,\tilde F^1,\tilde\phi)=(M_1,F^1(M_1),\phi)\oplus (M_2^*,F^0(M_2^*),p \phi)\oplus (W(k),W(k),p\sigma)$. 

As $(M_2^*,F^0(M_2^*),p \phi)=(M_2^*,F^0(M_2^*),\phi)
\otimes (W(k),W(k),p\sigma)$, it is easy to see that we can identify naturally $(\scrT(M),\phi)$ with a direct summand of $(\scrT(\tilde M),\tilde\phi)$ which is a $\pmb{GL}_{M_1}\times_{\Spec\,W(k)} \pmb{GL}_{M_2^*}\times_{\Spec\,W(k)} \pmb{GL}_{W(k)}$-submodule of $\scrT(\tilde M)$ in such a way that (a)  the filtration $F^i(\scrT(M))_{i\in\dbZ}$ of $\scrT(M)$ is the intersection of the filtration $F^i(\scrT(\tilde M))_{i\in\dbZ}$ of $\scrT(\tilde M)$ with $\scrT(M)$ and (b) the $\Gal(B(k))$-module $\scrT(H^1(D))$ gets naturally identified with a direct summand of the $\Gal(B(k))$-module $\scrT(H^1(\tilde D))$. 

Let $(t_{\alpha})_{\alpha\in\tilde\scrJ}$ be the family of tensors of $\scrT(\tilde M)$ formed by putting together $(t_{\alpha})_{\alpha\in\tilde\scrJ}$ and the projectors of $\tilde M$ on $M_1$ and $M_2^*$ along $M_2^*\oplus W(k)$ and $M_1\oplus W(k)$ (respectively). From the previous paragraph we easily get that if the Main Theorem holds in the context of $(\tilde D,(t_{\alpha})_{\alpha\in\tilde\scrJ})$, then it also holds in the context of $(D,(t_{\alpha})_{\alpha\in\scrJ})$. Thus (i) implies (ii).\endproof

\medskip\smallskip\noindent
{\bf 2.6. Two unipotent group schemes.} Let $\mu=\mu_{\text{can}}^{-1}:\dbG_m\to G$ and $M=F^1\oplus F^0$ be as in Subsection 2.1. Let 
$$\scrT(M)=\oplus_{i\in\dbZ} \tilde F^i(\scrT(M))$$ 
be the unique direct sum decomposition such that $\dbG_m$ acts via $\mu$ on $\tilde F^i(\scrT(M))$ as the $-i$-th power of the identity character of $\dbG_m$ (i.e., acts with weight $-i$). Thus each $\tilde F^i(\scrT(M))$ is a direct summand of $\scrT(M)$ that is uniquely determined by $\mu$ and $i$. The existence of $\mu:\dbG_m\to G$ implies that we have a direct sum decomposition
$$\Lie(G_{B(k)})\cap\End_{W(k)}(M)=\oplus_{i=-1}^1 \Lie(G_{B(k)})\cap \tilde F^i(\scrT(M)).\leqno (6)$$ 

Let $U_{\text{big}}$ be the smooth, unipotent, commutative, closed subgroup scheme of $\pmb{GL}_M$ whose Lie algebra is $\tilde F^{-1}(\End_{W(k)}(M)):=\End_{W(k)}(M)\cap\tilde F^{-1}(\scrT(M))$. Under the direct sum decomposition 
$$\End_{W(k)}(M)=\Hom_{W(k)}(F^1,F^0)\oplus\End_{W(k)}(F^1)\oplus\End_{W(k)}(F^0)\oplus\Hom_{W(k)}(F^0,F^1)$$ 
of $W(k)$-modules, we have $\tilde F^{-1}(\End_{W(k)}(M))=\Hom_{W(k)}(F^1,F^0)$. If $S$ is a commutative $W(k)$-algebra, then we have
$U_{\text{big}}(S)=1_{M\otimes_{W(k)} S}+\tilde F^{-1}(\End_{W(k)}(M))\otimes_{W(k)} S$. The intersection $\Lie(G_{B(k)})\cap \tilde F^{-1}(\End_{W(k)}(M))$ is a direct summand of $\tilde F^{-1}(\End_{W(k)}(M))$. Let $U$ be the smooth, closed subgroup scheme of $U_{\text{big}}$ (and thus also of $\pmb{GL}_M$) defined by the following rule on valued points:
$$U(S)=1_{M\otimes_{W(k)} S}+[\Lie(G_{B(k)})\cap \tilde F^{-1}(\End_{W(k)}(M))]\otimes_{W(k)} S.$$ 
We have $\Lie(U)=\Lie(G_{B(k)})\cap \tilde F^{-1}(\End_{W(k)}(M))$. As $U_{B(k)}$ is connected and as we have $\Lie(U_{B(k)})\subseteq\Lie(G_{B(k)})$, the group $U_{B(k)}$ is a subgroup of $G_{B(k)}$ (cf. [Bo, Ch. II, Subsect. 7.1]). As $U$ and $G$ are schematic closures in $\pmb{GL}_M$ of their generic fibres, $U$ is a smooth, closed subgroup scheme of $G$. 

\medskip\smallskip\noindent
{\bf 2.7. Review and complements to [Wi].}
Let $\scrP$ be the abelian group of functions $\zeta:\dbZ\to\dbZ$ which are periodic (i.e., for which there exists $c=c(\zeta)\in\dbN^{\ast}$ such that we have $\zeta(m)=\zeta(m+c)$ for all $m\in\dbZ$). For $\zeta\in\scrP$ let ${}^{\sigma}\zeta\in\scrP$ be defined by the rule ${}^{\sigma}\zeta(m):=\zeta(m+1)$. In [Wi] it is proved not only that we have a canonical decomposition $M=F^1\oplus F^0$ but also that there exist as well a unipotent element $u\in\pmb{GL}_M(W(k))$ and a direct sum decomposition $M=\oplus_{\zeta\in\scrP} M_{\zeta}$ of $W(k)$-modules such that the following four properties hold (cf. [Wi, p. 512]):

\medskip\noindent
{\bf (i)} if $v:=u-1_M\in\End_{W(k)}(M)$, then we have $v(F^1)\subseteq F^0$ and $v(F^0)=0$ (thus $v\in\tilde F^{-1}(\End_{W(k)}(M))$ has square $0$);

\smallskip\noindent
{\bf (ii)} for $i\in\{0,1\}$ we have $F^i=\oplus_{\zeta\in\scrP,\zeta(0)=i} M_{\zeta}$ (thus $M_{\zeta}=0$ if $\zeta(0)\notin\{0,1\}$);

\smallskip\noindent
{\bf (iii)} we have $\sigma_{\phi}:=\phi\mu(p)=uf_{\phi}$, where $f_{\phi}:M\arrowsim M$ is a $\sigma$-linear automorphism such that we have $f_{\phi}(M_{\zeta})=M_{{}^{\sigma}\zeta}$ for all $\zeta\in\scrP$ (thus $M_{\zeta}=0$ if $\text{Im}(\zeta)\nsubseteq\{0,1\}$);

\smallskip\noindent
{\bf (iv)} there exist $t\in\dbN$ and a series $0=\bar L_0\subseteq \bar L_1\subseteq\cdots\subseteq \bar L_t=M/pM$ of the $k$-vector space $M/pM$ which is normalized by $f_{\phi}$ modulo $p$ and which satisfies the relations $\bar L_j=\oplus_{\zeta\in\scrP} \bar L_j\cap M_{\zeta}/pM_{\zeta}$ and $\bar v(\bar L_j)\subseteq \bar L_{j-1}$ for all $j\in\{1,\ldots,t\}$, where $\bar v$ is $v$ modulo $p$.

\medskip
It is easy to see that properties (ii) to (iv) imply that for all $i,j\in\dbZ$ we have a direct sum decomposition $M=\oplus_{a,b\in\{0,1\}} f_{\phi}^i(F^a)\cap f_{\phi}^j(F^b)$. Thus the cocharacters $f_{\phi}^i(\mu)$ and $f_{\phi}^j(\mu)$ commute. Let $T_{B(k)}$ be the torus of $\pmb{GL}_{M[{1\over p}]}$ generated by the generic fibres of all cocharacters $f_{\phi}^i(\mu)$ with $i\in\dbZ$; it is the generic fibre of the image of a homomorphism $l:\dbG_m^m\to \pmb{GL}_{M}$ with $m>>0$ a positive integer. From this and the fact that $\Ker(l)$ is a group scheme of multiplicative type, we get that that the image of $l$ is a split torus of $\pmb{GL}_M$ which is the schematic closure $T$ of $T_{B(k)}$ in $\pmb{GL}_M$. The torus $T$ is the smallest subtorus of $\pmb{GL}_M$ through which $\mu$ factors and whose Lie algebra $\Lie(T)$ is normalized by $f_{\phi}$. 

As the $\sigma$-linear automorphism $f_{\phi}:M\arrowsim M$ and the element $u\in\pmb{GL}_M(W(k))$ are functorial (cf. loc. cit.), both of them fix each $t_{\alpha}$ with $\alpha\in\scrJ$. This implies that $T$ is a torus of $G$ and that $u\in G(W(k))$. 

Due to the property (i), if $u\neq 1_M$ then there exists a unique $\dbG_a$ subgroup scheme $U_0$ of $U$ such that we have $u\in U_0(W(k))$; it is defined by the identity $\Lie(U_0)=B(k)v\cap \Lie(U)$. If $u=1_M$, we define $U_0$ to be the identity section of $U$ (or $G$). We call $U_0$ as the {\it Wintenberger's unipotent group scheme} of either $(M,F^1,\phi)$ or $D$.

Let $g:=u^{-1}=1_M-v\in U_0(W(k))\leqslant U(W(k))\leqslant G(W(k))$. We have $g\phi=f_{\phi}\mu({1\over p})$. From this and the very definition of $T$ we get that: 

\medskip\noindent
{\bf (v)} we have $g\phi(\Lie(T))=f_{\phi}(\Lie(T))=\Lie(T)$. 

\medskip
If moreover $k=\bar k$, then there exists also a set $\scrJ_T$ which contains $\scrJ$ and a family of tensors $(t_{\alpha})_{\alpha\in\scrJ_T}$ of $F^0(\scrT(M))[{1\over p}]$ that extends $(t_{\alpha})_{\alpha\in\scrJ}$ and that has the following two extra properties (cf. Lemma 2.5.3 applied to $(M,g\phi,G,T)$ instead of $(M,\phi,G,\tilde G)$): 

\medskip\noindent
{\bf (vi)} the torus $T_{B(k)}$ is the subgroup of $\pmb{GL}_{M[{1\over p}]}$ that fixes $t_{\alpha}$ for all $\alpha\in\scrJ_T$ and 

\smallskip\noindent
{\bf (vii)} the $\sigma$-linear automorphism $g\phi$ of $\scrT(M[{1\over p}])$ fixes each $t_{\alpha}$ with $\alpha\in\scrJ_T$. 

\medskip\noindent
{\bf 2.7.1. Lemma.} {\it {\bf (a)} The multiplicity of the Newton polygon slope $0$ (resp. $1$) for $(M,h\phi)$ does not depend on $h\in 1_M+W(k)v\leqslant U_0(W(k))$ and thus it is equal to the multiplicity of the Newton polygon slope $0$ (resp. $1$) for $(M,g\phi)$ (equivalently, for $(M,\phi)$).

\smallskip
{\bf (b)} The multiplicity of the Newton polygon slope $-1$ for $(\Lie(G_{B(k)}),h\phi)$ does not depend on $h\in 1_M+W(k)v\leqslant U_0(W(k))$ and thus it is equal to the multiplicity of the Newton polygon slope $-1$ for $(\Lie(G_{B(k)}),g\phi)$ (equivalently, for $(\Lie(G_{B(k)}),\phi)$).}

\medskip
\proof
To prove (a) we can assume that $\bar v\neq 0$. Let $D_k(h)$ be the $p$-divisible group over $\Spec\,k$ whose Dieudonn\'e module is $(M,h\phi)$. As the series $0=\bar L_0\subseteq \bar L_1\subseteq\cdots\subseteq \bar L_t=M/pM$ is stable under the reductions of $h$ and $f_{\phi}$ modulo $p$ (cf. property 2.7 (iv)) and as each $\bar L_j$ with $j\in\{0,\ldots,t\}$ is normalized by $\mu_k$ (cf. properties 2.7 (ii) and (iv)), via the classical Dieudonn\'e theory we get the existence of a chain of epimorphisms 
$$D_k(h)[p]=B_t(h)\twoheadrightarrow  B_{t-1}(h)\twoheadrightarrow \cdots \twoheadrightarrow B_{1}(h)\twoheadrightarrow  B_0(h)=0$$ 
between finite flat group schemes over $\Spec\,k$ annihilated by $p$. More precisely, for each $j\in\{0,\ldots,t\}$ the Dieudonn\'e module of $B_j(h)$ is $\bar L_j$ together with the reductions modulo $p$ of $h\phi$ and of its Verschiebung map. As the reduction of $h-1_M\in\Lie(U_0)$ modulo $p$ is a scalar multiple of $\bar v$, from the property (iv) we get that $h$ modulo $p$ acts trivially on the quotients $\bar L_j/\bar L_{j-1}$ for all $j\in\{1,\ldots,t\}$. This implies that for each $j\in\{1,\ldots,t\}$ the kernel $C_j:=\Ker(B_j(h)\twoheadrightarrow B_{j-1}(h))$ does not depend indeed on $h\in 1_M+W(k)v$. Thus the $p$-rank of $D_j(h)_{\bar k}$ (resp. of $D_j(h)_{\bar k}^{\text{t}}$) equals to the $p$-rank of $\oplus_{j=1}^t C_{j,\bar k}$ (resp. of $\oplus_{j=1}^t C_{j,\bar k}^{\text{t}}$) and therefore does not depend on $h\in 1_M+W(k)v$. From this (a) follows.

To prove (b), we can assume that $k=\bar k$ and we first remark that for $z\in W(k)$ and $h=1_M+zv$, the following two properties hold.

\medskip
{\bf (i)} The multiplicity of the Newton polygon slope $-1$ for $(\End(M)[{1\over p}],h\phi)$ is the product of the multiplicities of the Newton polygon slopes $0$ and $1$ for $(M,h\phi)$ and thus (cf. (a)) it does not depend on $z$.

\smallskip
{\bf (ii)} The multiplicities of the Newton polygon slope $-1$ for $(\Lie(G_{B(k)}),h\phi)$ and $(\End(M)[{1\over p}]/\Lie(G_{B(k)}),h\phi)$ are equal to the multiplicities of the Newton polygon slope $0$ for $(\Lie(G_{B(k)}),p(h\phi))$ and $(\End(M)[{1\over p}]/\Lie(G_{B(k)}),p(h\phi))$ (respectively) and thus depend only on the reduction $\bar z$ of $z$ modulo $p$, and their sum does not depend on $\bar z\in k$ (cf. (i)).

\medskip
As $\bar z\in k=\dbG_a(k)$ and $\dbG_a$ is an integral $k$-scheme, from (i) and (ii) and their analogues over an algebraic closure of $k(\bar z)$ (with $\bar z$ as a variable), we get (via a standard specialization argument) that (b) holds. \endproof

\medskip\noindent
{\bf 2.7.2. The group scheme $G_{\text{min}}$.} In this subsubsection we assume that $k=\bar k$. Let $G_{\text{min},B(k)}$ be the smallest subgroup of $\pmb{GL}_{M[{1\over p}]}$ that contains the images of conjugates of $\mu_{B(k)}$ under powers of $\phi$; it is connected and we have $\phi(\Lie(G_{\text{min},B(k)}))=\Lie(G_{\text{min},B(k)})$. Let $G_{\text{min}}$ be the schematic closure of $G_{\text{min},B(k)}$ in $\pmb{GL}_M$. As $G$ is normalized by $\sigma_{\phi}$ and as $\mu$ factors through $G$, $G_{B(k)}$ contains the images of conjugates of $\mu_{B(k)}$ under powers of $\phi$. Thus $G_{\text{min}}$ is a closed subgroup scheme of $G$. From Lemma 2.5.3 we get that $G_{\text{min}}$ is the pull-back to $\Spec\,W(k)$ of a flat, closed subgroup scheme $G_{\text{min},\dbZ_p}$ of $G_{\dbZ_p}$. Let $G_{\text{min},\dbQ_p}$ be the generic fibre of $G_{\text{min},\dbZ_p}$. Groups like $G_{\text{min},\dbQ_p}$ were first considered in [Wi]. More precisely, $G_{\text{min},\dbQ_p}$ is the smallest subgroup of $\pmb{GL}_{M_{\dbZ_p}[{1\over p}]}$  such that its extension to $B(k)$ contains the image of $\mu_{B(k)}$ (cf. [Wi, Prop. 4.2.3]). 

\medskip\noindent
{\bf 2.7.3. Lemma.} {\it There exists a family of tensors $(t_{\alpha})_{\alpha\in\scrJ_{\text{min}}}$ of $\scrT(M_{\dbZ_p}[{1\over p}])$ such that $G_{\text{min},\dbQ_p}$ is the subgroup of $\pmb{GL}_{M_{\dbZ_p}[{1\over p}]}$ that fixes $t_{\alpha}$ for all $\alpha\in\scrJ_{\text{min}}$.}

\medskip
\proof
Let $\varepsilon\in\dbN$ be the largest number such that $\dbG_m^{\varepsilon}$ is a quotient of $G_{\text{min},\dbQ_p}$. The image of $\mu_{B(k)}$ in $\dbG_m^\varepsilon$ (extended to $B(k)$) is defined over $\dbQ_p$ and it is normalized by both $\sigma_{\phi}$ and $\phi$. From the very definition of $G_{\text{min},B(k)}$ we get that this image is $\dbG_m^\varepsilon$ (extended to $B(k)$) and therefore we have $\varepsilon\in\{0,1\}$. The relation $\varepsilon\in\{0,1\}$ also follows from the fact that each crystalline representation of $\Gal(B(k))$ of rank $1$ is isomorphic to an integral tensor power of the cyclotomic character of $\Gal(B(k))$.  

In this paragraph we assume that $\varepsilon=1$. Let $G^0_{\text{min},\dbQ_p}$  be the normal, connected subgroup of  $G_{\text{min},\dbQ_p}$ such that $G_{\text{min},\dbQ_p}/G^0_{\text{min},\dbQ_p}$ is isomorphic to $\dbG_m$. The image of $\mu$ is not contained in  $G^0_{\text{min},\dbQ_p}$ and thus $\mu$ is non-trivial. This implies that $G_{\text{min},\dbQ_p}$ is not contained in $\pmb{SL}_{M_{\dbZ_p}[{1\over p}]}:=\pmb{GL}^{\der}_{M_{\dbZ_p}[{1\over p}]}$. Therefore $G_{\text{min},\dbQ_p}/[G_{\text{min},\dbQ_p}\cap \pmb{SL}_{M_{\dbZ_p}[{1\over p}]}]$ is isomorphic to $\dbG_m$. This implies that $G^0_{\text{min},\dbQ_p}$ is the identity component of $G_{\text{min},\dbQ_p}\cap \pmb{SL}_{M_{\dbZ_p}[{1\over p}]}$. Therefore the restriction homomorphism $\Hom(\pmb{GL}_{M_{\dbZ_p}[{1\over p}]},\dbG_m)\to \Hom(G_{\text{min},\dbQ_p},\dbG_m)$ has finite cokernel. 

If $\varepsilon=0$, then the group $\Hom(G_{\text{min},\dbQ_p},\dbG_m)$ is trivial and therefore the homomorphism $\Hom(\pmb{GL}_{M_{\dbZ_p}[{1\over p}]},\dbG_m)\to \Hom(G_{\text{min},\dbQ_p},\dbG_m)$ is onto. 

Thus regardless of what $\varepsilon$ is, $\Hom(\pmb{GL}_{M_{\dbZ_p}[{1\over p}]},\dbG_m)\to \Hom(G_{\text{min},\dbQ_p},\dbG_m)$ has finite cokernel. Therefore $G_{\text{min},\dbQ_p}$ is the subgroup of $\pmb{GL}_{M_{\dbZ_p}[{1\over p}]}$ that fixes a family of tensors $(t_{\alpha})_{\alpha\in\scrJ_{\text{min}}}$ of $\scrT(M_{\dbZ_p}[{1\over p}])$, cf. [De2, Prop. 3.1 c)].\endproof

\medskip
As for each $\alpha\in\scrJ_{\text{min}}$ the tensor $t_{\alpha}$ is fixed by both $\mu_{B(k)}$ and $\sigma_{\phi}$, in fact we have $t_{\alpha}\in \{x\in F^0(\scrT(M))[{1\over p}]|\phi(x)=x\}$. 

\bigskip\smallskip\noindent
{\boldsectionfont 3. Global deformations}
\bigskip\smallskip

In this section we use the notations of Subsections 2.1, 2.6, and 2.7. Its main goal is to construct global deformations $\grG$ of $(M,F^1,\phi,G,(t_{\alpha})_{\alpha\in\scrJ})$ over $p$-adic completions $Q$ of ind-\'etale algebras over polynomial $W(k)$-algebras whose reductions modulo $p$ are geometrically connected and define spectra that have Zariski dense sets of $\bar k$-valued points. Each $\grG$ will be a filtered $F$-crystal over $\Spec\,Q/pQ$ endowed naturally with a group scheme over $\Spec\,Q$ canonically isomorphic to $G\times_{\Spec\,W(k)} \Spec\,Q$ and with a family of crystalline tensors, whose pull-back via a suitable $k$-valued point of $\Spec\,Q/pQ$ is $(M,F^1,\phi,G,(t_{\alpha})_{\alpha\in\scrJ})$. 

If $g\in U_0(W(k))$ and $T$ are as in Subsection 2.7, we will use such a global deformation $\grG$ of $(M,F^1,\phi,G,(t_{\alpha})_{\alpha\in\scrJ})$ which connects $(M,F^1,\phi,G,(t_{\alpha})_{\alpha\in\scrJ})$ with $(M,F^1,g\phi,G,(t_{\alpha})_{\alpha\in\scrJ})$ i.e., $(M,F^1,g\phi,G,(t_{\alpha})_{\alpha\in\scrJ})$ is as well the pull-back of $\grG$ via a $k$-valued point of $\Spec\,Q/pQ$. This will allow us to reduce the proof of the Main theorem to the case when $G$ is the torus $T$ and we have $g=u^{-1}=1_M$. 

In Subsection 3.1 we develop a language on connections needed to state and prove the Theorem 3.2 which pertains to the existence of moduli spaces of connections and implicitly of some global deformations of $(M,F^1,\phi,G,(t_{\alpha})_{\alpha\in\scrJ})$. Subsection 3.3 studies the number of geometric points of the fibres of our moduli spaces of connections. Subsection 3.4 translates Theorem 3.2 in terms of $p$-divisible groups as allowed by Theorem 2.3.4 and by a weaker variant of it for $p=2$. In Subsection 3.5 we apply  Subsection 3.4 to show that to prove the Main Theorem we can assume that $G$ is the torus $T$ and $g=u=1_M$. 

In this section, if $\star$ is a $W(k)$-algebra and $N$ is a $W(k)$-module, let 
$N_{\star}:=N\otimes_{W(k)} \star$. In Subsections 3.1 to 3.3 we will not use the fact that (even if $p=2$) the triple $(M,F^1,\phi)$ is the filtered Dieudonn\'e module of a $p$-divisible group $D$ over $\Spec\,W(k)$.

\medskip\smallskip\noindent
{\bf 3.1. A language on connections.}  Let $m\in\dbN$. Let $Y=\Spec\,R$ be an affine $W(k)$-scheme which is either smooth of relative dimension $m$ or $\Spec\,W(k)[[z_1,\ldots,z_m]]$. If $R$ is smooth over $W(k)$, we will assume that the $R$-module $\Omega_{R/W(k)}$ is free and we fix elements $z_1,\ldots,z_m\in R$ such that $\Omega_{R/W(k)}=\oplus_{j=1}^m Rdz_{j}$; the $W(k)[z_1,\ldots,z_m]$-algebra $R$ is \'etale. If $R=W(k)[[z_1,\ldots,z_m]]$, we will define $\Omega_{R/W(k)}:=\oplus_{j=1}^m Rdz_j$ (i.e., we will work with $(z_1,\ldots,z_m)$-adic completions of modules of relative differentials). We fix a Frobenius lift $\Phi_R$ of either $R^\wedge$ or $Y^\wedge$ and we consider an element
$$\tau\in \pmb{GL}_M(R^\wedge).$$

\noindent
Below we will identify points of $\pmb{GL}_M$ or $G^{\wedge}$  with values in a flat, affine $\Spec\,W(k)$-scheme $\Spec\,*$ with $*$-linear automorphisms of $M_{*}$. Let 
$$\tilde M:=(M+{1\over p}F^1)\otimes_{W(k)}{}_{\sigma} W(k)=(F^0\oplus {1\over p}F^1)\otimes_{W(k)}{}_{\sigma} W(k)\subset (M\otimes_{W(k)}{}_{\sigma} W(k))[{1\over p}].$$ 
We have a $W(k)$-linear isomorphism $\phi\otimes 1_{W(k)}:\tilde M\arrowsim M$.

\medskip\noindent
{\bf 3.1.1. Definition and notations.} {\bf (a)} Let $f:\Spec\,\scrQ^\wedge\to Y^\wedge$ be a formally \'etale morphism of affine $\Spec\,W(k)$-schemes whose reduction modulo $p$ is \'etale, with $Q$ a regular $W(k)$-algebra. Let $\Omega_{(\scrQ/p^n\scrQ)/W_n(k)}$ be the extension of $\Omega_{R/W(k)}/p^n\Omega_{R/W(k)}$ via $f$ modulo $p^n$. Let $\Phi_{\scrQ}$ be the only Frobenius lift of $\Spec\,\scrQ^\wedge$ (or $\scrQ^\wedge$) that satisfies the identity $f\circ\Phi_{\scrQ}=\Phi_{R}\circ f$. Let $d\Phi_{\scrQ}/p$ be the differential of $\Phi_{\scrQ}$ divided by $p$. 

Let $\phi_{\scrQ}:\tilde M_{\scrQ^\wedge}\arrowsim M_{\scrQ^\wedge}$ be the $\scrQ^\wedge$-linear isomorphism which maps $x\otimes 1$ with $x\in \tilde M$ to $\tau\circ f(\phi\otimes 1_{W(k)}(x))$. Let $\phi_{0,\scrQ}:M_{\scrQ^\wedge}\to M_{\scrQ^\wedge}$ be the $\Phi_{\scrQ}$-linear map which maps $x\otimes 1$ with $x\in M$ to $\tau\circ f(\phi(x))$. Let $\phi_{1,\scrQ}:F^1_{\scrQ^\wedge}\to M_{\scrQ^\wedge}$ be the $\Phi_{\scrQ}$-linear map which maps $x\otimes 1$ with $x\in F^1$ to ${1\over p}\phi_{0,Q}(x\otimes 1)$. Let $\scrM_{\scrQ^\wedge}:=(M_{\scrQ^\wedge},F^1_{\scrQ^\wedge},\phi_{0,\scrQ},\phi_{1,\scrQ})$. 

For $n\in\dbN^{\ast}$, we also denote by $\phi_{\scrQ}$, $\phi_{0,\scrQ}$, $\phi_{1,\scrQ}$, and $d\Phi_{\scrQ}/p$  their reductions modulo $p^n$. Let $\scrM_{\scrQ/p^n\scrQ}=(M_{\scrQ/p^n\scrQ},F^1_{\scrQ/p^n\scrQ},\phi_{0,\scrQ},\phi_{1,\scrQ})$ be the reduction modulo $p^n$ of $\scrM_{\scrQ^\wedge}$; it is an object of the category $\scrM\scrF_{[0,1]}(\scrQ)$ of Subsubsection 2.3.2. 

\smallskip
{\bf (b)} Let $\Delta$ be the flat connection on either $M_{\scrQ^\wedge}$  or $M_{\scrQ}/p^nM_{\scrQ}$ that annihilates $M\otimes 1$. 

\smallskip
{\bf (c)} Let $n\in\dbN^{\ast}$. By a {\it connection on} $\scrM_{\scrQ/p^n\scrQ}$  we will mean a connection  
$$
\nabla:M_{\scrQ/p^n\scrQ}\to M_{\scrQ/p^n\scrQ}\otimes_{\scrQ/p^n\scrQ} \Omega_{(\scrQ/p^n\scrQ)/W_n(k)}=M\otimes_{W(k)} \Omega_{(\scrQ/p^n\scrQ)/W_n(k)}
$$ 
for which the following two equations hold (to be compared with axioms 2.3 2 (c) and (d)):
$$\nabla\circ\phi_{0,\scrQ}(x)=p(\phi_{0,\scrQ}\otimes d\Phi_{\scrQ}/p)\circ\nabla(x)\;\;\;\;\forall x\in M_{\scrQ/p^n\scrQ}\;\;\;\;\text{and} \leqno (7)$$
$$\nabla\circ\phi_{1,\scrQ}(x)=(\phi_{0,\scrQ}\otimes d\Phi_{\scrQ}/p)\circ\nabla(x)\;\;\;\; \forall x\in F^1_{\scrQ/p^n\scrQ}.\leqno (8)$$
Equivalently, by a connection on $\scrM_{\scrQ/p^n\scrQ}$  we will mean a connection $\nabla$ on $M_{\scrQ/p^n\scrQ}$ for which $\phi_{\scrQ}$ is {\it horizontal} i.e., for which we have an identity of maps
$$\nabla\circ\phi_{\scrQ}=(\phi_{\scrQ}\otimes 1_{\Omega_{(\scrQ/p^n\scrQ)/W_n(k)}})\circ \tilde\Phi_{\scrQ}^*(\nabla):\tilde M_{\scrQ/p^n\scrQ}\to M\otimes_{W(k)} \Omega_{(\scrQ/p^n\scrQ)/W_n(k)}.\leqno (9)$$
Here $\tilde\Phi_{\scrQ}^*(\nabla)$ is the connection on $\tilde M_{\scrQ/p^n\scrQ}$ induced canonically by the extension $\Phi_{\scrQ}^*(\nabla)$ of $\nabla$ via $\Phi_{\scrQ}$ as follows. If  $\nabla_{\text{lift}}$ is an arbitrary connection on $M_{\scrQ/p^{n+1}\scrQ}$ that lifts $\nabla$, then $\Phi_{\scrQ}^*(\nabla_{\text{lift}})$ restricts to a connection on $\tilde M'_{n+1}:=[(F^1\oplus pF^0)/(p^{n+1}F^1\oplus p^{n+1}F^0)]\otimes_{W(k)} {}_{\sigma} \scrQ/p^{n+1}\scrQ$ and thus it induces a connection $\tilde\Phi_{\scrQ}^*(\nabla)$ on the $\scrQ/p^n\scrQ$-module $\tilde M'_{n+1}/p^n\tilde M'_{n+1}$ which is canonically identified with $\tilde M_{\scrQ/p^n\scrQ}$ (via multiplication by ${1\over p}$). The connection $\tilde\Phi_{\scrQ}^*(\nabla)$ does not depend on the choice of the lift $\nabla_{\text{lift}}$ of $\nabla$.

\smallskip
{\bf (d)} Let $\tilde G$ be a flat, closed subgroup scheme of $\pmb{GL}_M$ with the properties that $\Lie(\tilde G_{B(k)})$ is normalized by $\phi$ and $\Lie(\tilde G_{B(k)})\cap \End_{W(k)}(M)$ is normalized by $\mu$; for instance, $\tilde G$ could be either $G$ (cf. formula (6)) or $\pmb{GL}_M$. We assume that $\tau\in\tilde G(R^\wedge)$. We say the connection $\nabla$ of (c) {\it respects the $\tilde G$-action}, if for each $j\in \{1,\ldots,m\}$, the $\scrQ/p^n\scrQ$-linear endomorphism of $M_{\scrQ/p^n\scrQ}$ which for $x\in M$ maps $x\otimes 1$ to $\nabla({{\partial}\over {\partial z_j}})(x\otimes 1)\in M_{\scrQ/p^n\scrQ}$,
is an element of $(\Lie(\tilde G_{B(k)})\cap\End_{W(k)}(M))_{\scrQ/p^n\scrQ}$.

\medskip\smallskip\noindent
{\bf 3.2. Basic Theorem (the moduli of connections).} {\it Let $Y=\Spec\,R$ and $\Phi_R$ be as in Subsection 3.1. Let $\tilde G$ and $\tau\in \tilde G(R^\wedge)$ be as in Subsubsection 3.1.1 (d). Then the following three things hold:
\medskip
{\bf (a)} For each $n\in\dbN^{\ast}$ there exists an \'etale $R$-algebra $Q_n$ such that for the naturally induced formally \'etale, affine morphism 
$\ell_n:\Spec\,Q^\wedge_n\to Y^\wedge$ there exists a connection
$\nabla_n$ on $\scrM_{Q_n/p^nQ_n}$ which is integrable and nilpotent modulo $p$, which respects the $\tilde G$-action, and which satisfies the following universal property: 

\noindent
{\bf (UP)} For each formally \'etale morphism $f:\Spec\,\scrQ^\wedge\to Y^\wedge$ of affine $\Spec\,W(k)$-schemes whose reduction modulo $p$ is \'etale (with $\scrQ$ a regular $W(k)$-algebra) and for every connection $\tilde \nabla$ on $\scrM_{\scrQ/p^n\scrQ}$ which respects the $\tilde G$-action, there exists a unique morphism $f_n:\Spec\,\scrQ^\wedge\to\Spec\,Q^\wedge_n$ such that $f=\ell_n\circ f_n$ and $\tilde\nabla$ is the extension of $\nabla_n$ via $f_n$ modulo $p^n$.

\medskip
{\bf (b)} For $n\in\dbN^{\ast}$ we consider the  unique morphism $\ell^{(n)}:\Spec\,Q_{n+1}^\wedge\to \Spec\,Q_n^\wedge$ such that $\ell_{n+1}=\ell_n\circ\ell^{(n)}$ and the extension of $\nabla_n$ via $\ell^{(n)}$ modulo $p^n$ is $\nabla_{n+1}$ modulo $p^n$, cf. (a). Let $\ell^{(0)}:=\ell_1$. Then the reduction modulo $p$ of the projective system $(\ell^{(n)})_{n\in\dbN}$ is an Artin--Schreier tower in the sense of Subsection 2.4.

\smallskip
{\bf (c)} Let $\Spec\,Q_{\infty}$ be the projective limit of $\Spec\,Q_n^\wedge$ under the transition morphisms $\ell^{(n)}$, $n\in\dbN^{\ast}$. We fix a connected component $\Spec\,Q$ of $\Spec\,Q_{\infty}^\wedge$. Let $\ell:\Spec\,Q\to Y^\wedge$ 
be the resulting morphism. Then the image of $\ell$ modulo $p$ is an open subscheme of $Y_k$. If $R$ is a smooth $W(k)$-algebra, then each point of $\Spec\,Q/pQ$ specializes to a point of $\Spec\,Q/pQ$ whose residue field is an algebraic extension of $k$.} 

\medskip
\proof
To prove (a) and (b), we start by considering a morphism $f:\Spec\,\scrQ^\wedge\to Y^\wedge$ of $\Spec\,W(k)$-schemes as in (a) with the property that there exists a connection $\tilde\nabla$ on $\scrM_{\scrQ/p^n\scrQ}$ which respects the $\tilde G$-action. Writing 
$$\tilde\nabla=\Delta+\omega$$ 
with $\omega\in (\Lie(\tilde G_{B(k)})\cap \End_{W(k)}(M))\otimes_{W(k)} \Omega_{(\scrQ/p^n\scrQ)/W_n(k)}$, the equation (9) becomes
$$\omega\circ\phi_{\scrQ}+\Delta\circ\phi_{\scrQ}=(\phi_{\scrQ}\otimes 1_{\Omega_{(\scrQ/p^n\scrQ)/W_n(k)}})\circ \tilde\Phi_{\scrQ}^*(\Delta)+(\phi_{\scrQ}\otimes 1_{\Omega_{(\scrQ/p^n\scrQ)/W_n(k)}})\circ \tilde\Phi_{\scrQ}^*(\omega).$$
This equation can be rewritten as
$$\omega-\phi_{\scrQ}\circ\tilde \Phi_{\scrQ}^*(\omega)\circ\phi_{\scrQ}^{-1}=\Delta-\phi_{\scrQ}\circ\tilde\Phi_{\scrQ}^*(\Delta)\circ\phi_{\scrQ}^{-1}=:\scrR_n\leqno (10)$$ 
(here $\phi_{\scrQ}\circ\tilde \Phi_{\scrQ}^*(\omega)\circ\phi_{\scrQ}^{-1}$ and $\phi_{\scrQ}\circ\tilde\Phi_{\scrQ}^*(\Delta)\circ\phi_{\scrQ}^{-1}$ are shorter forms for extensions via the reduction modulo $p^n$ of the $\scrQ^\wedge$-linear isomorphism $\phi_{\scrQ}$). As $\phi_Q=\tau\circ f\circ (\phi\otimes 1_Q)$, $(\phi\otimes 1_Q)\tilde\Phi_{\scrQ}^*(\Delta)(\phi\otimes 1_Q)^{-1}=\Delta$, and $\tau\circ f\in\tilde G(R^\wedge)$, we get that 
$$\scrR_n=-(\tau\circ f)d(\tau\circ f)^{-1}=(\tau\circ f)^{-1}d(\tau\circ f)\in (\Lie(\tilde G_{B(k)})\cap \End_{W(k)}(M))\otimes_{W(k)} \Omega_{(\scrQ/p^n\scrQ)/W_n(k)}.\leqno (11)$$
\indent
Let $d$ be the dimension of $\tilde G_{B(k)}$. We consider a $W(k)$-basis $\{e_1,\ldots,e_d\}$ of $\Lie(\tilde G_{B(k)})\cap \End_{W(k)}(M)$. With $x_{ij}\in\scrQ/p^n\scrQ$ thought as arbitrary variables, we can write 
$$\omega=\sum_{i=1}^d\sum_{j=1}^m x_{ij}e_i\otimes dz_j\in (\Lie(\tilde G_{B(k)})\cap \End_{W(k)}(M))\otimes_{W(k)} \Omega_{(\scrQ/p^n\scrQ)/W_n(k)}.$$
The equation (10) becomes the following equation
$$\sum_{i=1}^d\sum_{j=1}^m x_{ij}e_i\otimes dz_j-\sum_{i=1}^d\sum_{j=1}^m\Phi_{\scrQ}(x_{ij})\tau\circ f((p\phi(e_i)))\otimes d\Phi_{\scrQ}/p(z_j)=\scrR_n\leqno (12)$$
between elements of $(\Lie(\tilde G_{B(k)})\cap \End_{W(k)}(M))\otimes_{W(k)} \Omega_{(\scrQ/p^n\scrQ)/W_n(k)}$.

We now construct $Q_n$ by induction on $n\in\dbN^{\ast}$ as follows.
We first assume that $n=1$. Then with respect to the $\scrQ/p\scrQ$-basis $\{e_i\otimes dz_j|i\in\{1,\ldots,d\},j\in\{1,\ldots, m\}\}$ of $(\Lie(\tilde G_{B(k)})\cap \End_{W(k)}(M))\otimes_{W(k)} \Omega_{(\scrQ/p\scrQ)/k}$, the equation (12) becomes an Artin--Schreier systems of equations in $dm$ variables over $R/pR$ of the form
$$x_{ij}=L_{ij}(x_{11}^p,\ldots,x_{dm}^p)+a_{ij}(1)\leqno (13)$$
where $a_{ij}(1)\in R/pR$ are such that $\scrR_1=\sum_{i=1}^d\sum_{j=1}^m a_{ij}(1)e_i\otimes dz_j$ and where each $L_{ij}$ is a homogeneous linear form obtained from $\sum_{i=1}^d\sum_{j=1}^m\Phi_{\scrQ}(x_{ij})\tau\circ f((p\phi(e_i)))\otimes d\Phi_{\scrQ}/p(z_j)$. Let $Q_1/pQ_1$ be the \'etale $R/pR$-algebra that parametrizes solutions of the system of equations (13).
Standard arguments that involve a lift of (13) to a system over $R$, the Jacobi criterion for \'etaleness, and localizations, show the existence of an \'etale $R$-algebra $Q_1$ such that the morphism $\ell_{1}:\Spec\,Q_1^\wedge\to Y^\wedge$ lifts the morphism $\Spec\,Q_1/pQ_1\to Y_k=\Spec\,R/pR$ defined by (13).  The universal property (UP) of (a) holds for $n=1$, cf. constructions.

The passage from $n-1$ to $n$ for $n\ge 2$ goes as follows. We assume that $Q_{n-1}$ has been constructed and that the morphism $\ell_{n-1}:\Spec\,Q_{n-1}^\wedge\to Y^\wedge$ has the universal property (UP) of (a). To construct $Q_n$, let $f_{n-1}:\Spec\,\scrQ^\wedge\to \Spec\,Q_{n-1}^\wedge$ be the unique morphism such that $\tilde\nabla$ modulo $p^{n-1}$ is the extension of $\nabla_{n-1}$ via $f_{n-1}$ modulo $p^{n-1}$. We write 
$$\omega=\omega_{n-1}^{\text{lift}}+\sum_{i=1}^d\sum_{j=1}^m p^{n-1}y_{ij}e_i\otimes dz_j,\leqno (14)$$
where $\omega_{n-1}^{\text{lift}}\in (\Lie(\tilde G_{B(k)})\cap \End_{W(k)}(M))\otimes_{W(k)} \Omega_{(Q_{n-1}/p^nQ_{n-1})/W_n(k)}$ lifts $\omega_{n-1}:=\nabla_{n-1}-\Delta\in (\Lie(\tilde G_{B(k)})\cap \End_{W(k)}(M))\otimes_{W(k)} \Omega_{(Q_{n-1}/p^{n-1}Q_n)/W_{n-1}(k)}$ and where $y_{ij}\in\scrQ/p\scrQ$ are again thought as arbitrary variables. Plugging the expression (14) of $\omega$ in the equation (12) and recalling that (10) holds with $(\omega,n)$ replaced by $(\omega_{n-1},n-1)$, one gets that the equation (12) becomes equivalent to a new Artin--Schreier system of equations in $dm$ variables over $Q_{n-1}/pQ_{n-1}$ which is of the form
$$y_{ij}=L_{ij}(y_{11}^p,\ldots,y_{dm}^p)+a_{ij}(n),\leqno (15)$$
where $a_{ij}(n)\in Q_{n-1}/pQ_{n-1}$ and where the $L_{ij}$'s are the {\it same} homogeneous forms.  We take $Q_n/pQ_n$ to be the \'etale $Q_{n-1}/pQ_{n-1}$-algebra that parametrizes solutions of the system of equations (15). Standard arguments as above
show the existence of an \'etale $Q_{n-1}$-algebra $Q_n$ such that the the morphism $\ell^{(n-1)}:\Spec\,Q_n^\wedge\to \Spec\,Q_{n-1}^\wedge$ lifts the morphism $\Spec\,Q_n/pQ_n\to \Spec\,Q_{n-1}/pQ_{n-1}$ defined naturally by (15).  From the very constructions we get that the universal property (UP) of (a) holds for $n$. 

Part (b) follows from the very definitions and from the key fact that the homogeneous linear forms showing up in (13) and (15) are the same. 

We are left to check that $\nabla_1$ is nilpotent and that each $\nabla_n$ is integrable.  We note that $\nabla_1$ annihilates the image of $\phi_{0,Q_1}$ modulo $p$ (cf. equation (7)) and it maps the image of $\phi_{1,\scrQ_1}$ to the image of $\phi_{0,Q_1}$ modulo $p$ tensored with $\Omega_{(Q_1/pQ_1)/k}$ (cf. equation (8)). From this we get that the $p$-curvature of $\nabla_1$ has square $0$. Thus $\nabla_1$ is nilpotent modulo $p$. 

Let $\nabla_{\infty}$ be the connection on $M_{Q_{\infty}^\wedge}$ whose reduction modulo $p^n$ is the natural extension of $\nabla_n$ for all $n\in\dbN^{\ast}$. To prove that each $\nabla_n$ is integrable, as the $k$-morphisms $\Spec\,Q_{\infty}/pQ_{\infty}\to \Spec\,Q_n/pQ_n$ are onto (cf. (b) and Theorem 2.4.1 (b)), it suffices to show that $\nabla_{\infty}$ is integrable i.e., its curvature $\scrH_{\infty}\in \End_{W(k)}(M)\otimes_{W(k)} \Omega^2_{Q^\wedge_{\infty}/W(k)}$ is $0$. Here  $\Omega^2_{Q^\wedge_{\infty}/W(k)}$ is the second exterior power of $\Omega_{Q^\wedge_{\infty}/W(k)}$. 

By induction on $s\in\dbN$ we prove that $\scrH_{\infty}$ is divisible by $p^s$. The case $s=0$ is obvious and the passage from $s$ to $s+1$ goes as follows. As an element of $\End_{W(k)}(F^0\oplus {1\over p}F^1)\otimes_{W(k)} \Omega^2_{Q^\wedge_{\infty}/W(k)}$, $\scrH_{\infty}$ is divisible by $p^{s-1}$ (if $s=0$, this means that it has at most $p$ in its denominators). When we pull-back via $\Phi_{Q_{\infty}}$, due to the upper index $2$ in  $\Omega^2_{Q^\wedge_{\infty}/W(k)}$, the curvature of $\tilde\Phi_{Q_{\infty}}^*(\nabla_{\infty})$ is an element of $\End_{W(k)}(F^0\oplus {1\over p}F^1)\otimes_{W(k)} \Omega^2_{Q^\wedge_{\infty}/W(k)}\otimes_{Q_{\infty}^\wedge} {}_{\Phi_{Q_{\infty}}} Q^\wedge_{\infty}=\End_{W(k)}(\tilde M)\otimes_{W(k)} \Omega^2_{Q^\wedge_{\infty}/W(k)}$ divisible by $p^{s-1}p^2=p^{s+1}$. As the $\scrQ_{\infty}^\wedge$-linear isomorphism $\phi_{Q_{\infty}}:\tilde M_{\scrQ_{\infty}^\wedge}\to M_{\scrQ_{\infty}^\wedge}$ is horizontal, one gets that $\scrH_{\infty}$ is also divisible by $p^{s+1}$. This ends the induction. Thus $\scrH_{\infty}=0$ i.e., $\nabla_{\infty}$ is integrable. Therefore (a) and (b) hold.

Part (c) follows from (b) and Theorem 2.4.1 (d) and (e). 
\endproof

\medskip\smallskip\noindent
{\bf 3.3. On special fibres of $\ell_n:\Spec\,Q^\wedge_n\to Y^\wedge$.} In this subsection we will study the fibres of the reduction modulo $p$ of the morphism $\ell_n:\Spec\,Q_n^\wedge\to Y^\wedge$ of Theorem 3.2 (a). 

Let $J_y$ be a maximal ideal of $R^\wedge$ such that $R^\wedge/J_y$ is a finite field extension $k_1$ of $k$. The corresponding closed point of either $Y_k=\Spec\,R/pR$ or $Y^\wedge=\Spec\,R^\wedge$ will be denoted simply by $y$. Let $I_y$ be a prime ideal of $R^\wedge$ contained in $J_y$ and such that $R^\wedge/I_y=W(k_1)$ and $\Phi_R(I_y)\subseteq I_y$; thus the morphism $\Spec\,W(k_1)=\Spec\,R^\wedge/I_y\hookrightarrow \Spec\,R^\wedge$ is a Teichm\"uller lift. Let $h_y\in \tilde G(W(k_1))$ be the reduction modulo $I_y$ of $\tau\in \tilde  G(R^\wedge)$. 

Let $m_y(0)$ and $m_y(1)$ be the multiplicities of the Newton polygon slopes $0$ and $1$ (respectively) for $(M_{W(k_1)},h_y(\phi\otimes\sigma_{k_1}))$. Let $m_y^{\tilde G}(-1)$ be the multiplicity of the Newton polygon slope $-1$ for $(\Lie(\tilde G_{B(k_1)}),h_y(\phi\otimes\sigma_{k_1}))$. We have the following universal relations
$$m_y(0)m_y(1)=m_y^{\pmb{GL}_M}(-1)\ge m_y^{\tilde G}(-1).\leqno (16)$$
We note that the integers $m_y(0)$, $m_y(1)$, and $m_y^{\tilde G}(-1)$ depend only on $h_y$ modulo $p$. 

Let $\Omega_y:= \Omega_{R^\wedge/W(k)}\otimes_{R^\wedge} R^\wedge/J_y\otimes_{k_1} \bar k=\oplus_{j=1}^m \bar kdz_j$. We consider the $\sigma_{\bar k}$-linear map $\omega_y:\Omega_y\to \Omega_y$ induced naturally by $d\Phi_R/p$. Let $c_y\in\{0,\ldots,m\}$ be the dimension of the $\dbF_p$-vector space of solutions of the equation $x=\omega_y(x)$ in $x\in \Omega_y$, cf. Theorem 2.4.1 (b). We recall from the proof of  Theorem 2.4.1 (b), that $c_y$ can be computed as follows. If $[\omega_y]$ is the matrix representation of $\omega_y$ with respect to the $\bar k$-basis $\{dz_1,\ldots,dz_m\}$ of $\Omega_y$ and if $[\varpi_y]\in M_m(\bar k)$ is an upper triangular matrix of the form $\ddag^{-1}[\omega_y]\ddag^{[p]}$ with $\ddag\in\pmb{GL}_m(\bar k)$, then $c_y$ is the number of non-zero entries on the main diagonal of $[\varpi_y]$.  

\medskip\noindent
{\bf 3.3.1. Theorem.} {\it For $n\in\dbN^{\ast}$ we consider the morphism $\ell_n:\Spec\,Q_n\to\Spec\,R$ of Theorem 3.2. Then with the above notations, the number of geometric points of the fibre of $\ell_n$  over the closed point $y$ of $\Spec\,R^\wedge$ defined by $J_y$ is precisely $p^{nc_ym_y^{\tilde G}(-1)}$.}

\medskip
\proof
We can assume that $k=k_1=\bar k$. Based on Theorems 2.3 (b) and 2.4.1 (b), it suffices to consider the case $n=1$ and to show that the number of solutions in $k$ of the reduction modulo $J_y$ of the system of equations
$$x_{ij}=L_{ij}(x_{11}^p,\ldots,x_{dm}^p)\leqno (17)$$
indexed by $(i,j)\in\{1,\ldots,d\}\times \{1,\ldots,m\}$, is exactly $p^{c_ym_y^{\tilde G}(-1)}$. But $(17)$ modulo $J_y$ is the same as the reduction modulo $J_y$ of the system of equations obtained from (12) by replacing its right hand side by $0$ and thus it is the same as the equation
$$ \sum_{i=1}^d\sum_{j=1}^m x_{ij}\bar e_i\otimes dz_j-\sum_{i=1}^d\sum_{j=1}^m x_{ij}^p\bar h_y(p\phi(\bar e_i))\otimes \omega_y(dz_j)=0\leqno (18)$$
in $dm$ variables over $k$. Here $\bar e_i$ and $\bar h_y$ are the reductions modulo $p$ of $e_i$ and $h_y$. Thus we are left to show that the number of solutions of the following equation obtained from $(18)$
$$x=\varphi_y\otimes \omega_y(x)\leqno (19)$$
in $x\in L_y\otimes_k \Omega_y$, is $p^{c_ym_y^{\tilde G}(-1)}$; here $L_y:=[\Lie(\tilde G_{B(k)})\cap \End_{W(k)}(M)]\otimes_{W(k)} k$ and $\varphi_y:L_y\to L_y$ is the $\sigma$-linear map induced by $h_y\phi\circ p1_{\End_{B(k)}(M[{1\over p}])}$. But the dimension of the $\dbF_p$-vector space of solutions of Equation (19) is the product of the dimensions of the $\dbF_p$-vector spaces of solutions of the equations $x=\varphi_y(x)$ in $x\in L_y$  and $x=\omega_y(x)$ in $x\in \Omega_y$. As the first dimension is the multiplicity  of the Newton polygon slope $0$ for
$(\Lie(\tilde G_{B(k)}),h_y\phi\circ p1_{\Lie(\tilde G_{B(k)})})$ and thus is the multiplicity $m_y^{\tilde G}(-1)$ of the Newton polygon slope $-1$ for $(\Lie(\tilde G_{B(k)}),h_y\phi)$ and as the second dimension is $c_y$ by very definitions, we get that Equation (19) has $p^{c_ym_y^{\tilde G}(-1)}$ solutions.
\endproof

\medskip\noindent
{\bf 3.3.2. Corollary.} {\it Suppose $R$ is a smooth $W(k)$-algebra such that $Y_k=\Spec\,R/pR$ is connected. Let $Z_k$ be the largest open subscheme of $Y_k$ with the property that for all closed points $y\in Z_k$, the product $c_ym_y^{\tilde G}(-1)$ is constant (i.e., it takes the maximal possible value). Then for all $n\in\dbN^{\ast}$, $Z_k$ is the largest open subscheme of $Y_k$ with the property that the reduction modulo $p$ of the morphism $\ell_n:\Spec\,Q_n^\wedge\to Y^\wedge$ is an \'etale cover above $Z_k$ (i.e., $Z_k$ is the open, dense stratum of the Artin--Schreier stratification of $Y_k$ defined by (13), cf. Definition 2.4.2). If for a closed point $y\in Z_k(\bar k)$ we have $c_ym_y^{\tilde G}(-1)=0$, then $Z_k=Y_k$ and in fact for all $n\in\dbN^{\ast}$ we have $Q_n^\wedge=R^\wedge$.}

\medskip\smallskip\noindent
{\bf 3.4. Geometric translation of a concrete case of Theorem 3.2.} 
Now we will translate Theorem 3.2 (a) and (b) in terms of $p$-divisible groups. Let $\nabla$ be the unique connection on $M_Q$ whose reduction modulo $p^n$ is the natural extension of $\nabla_n$ for all $n\in\dbN^{\ast}$; this makes sense due to the universal property (UP) of Theorem 3.2 (a). For each $n\in\dbN^{\ast}$ the reduction modulo $p^n$ of $\nabla$ respects the $\tilde G$-action, cf. Theorem 3.2 (a). Thus if $\tilde G=G$, then  $\nabla$ annihilates $t_{\alpha}\in\scrT(M\otimes_{W(k)} Q[{1\over p}])$ for all $\alpha\in\scrJ$.

As $Q/pQ$ is either an ind-\'etale $k[z_1,\ldots,z_m]$-algebra or an ind-\'etale $k[[z_1,\ldots,z_m]]$-algebra, the set $\{z_1,\ldots,z_m\}$ is a finite $p$-basis for $Q/pQ$ in the sense of [BM, Def. 1.1.1]. 

We consider the following condition:

\medskip
{\bf (*)} for each closed point $y\in Y_k$ whose residue field is an algebraic extension of $k$, we have $m_y(0)m_y(1)=0$.

\medskip\noindent
{\bf 3.4.1. Theorem.} {\it If $R=W(k)[[z_1,\ldots,z_m]]$, then we choose the connected component $\Spec\, Q$ of $\Spec\, Q_{\infty}^\wedge$ such that the morphism $\Spec\,Q/pQ\to\Spec\,R/pR$ is onto (equivalently, it is a pro-\'etale cover). If $p=2$ we assume that the condition 3.4 (*) holds. Then the following two properties hold:

\medskip
{\bf (a)} There exists a unique $p$-divisible group $\scrD_{Q/pQ}$ over $\Spec\,Q/pQ$ such that the evaluation of $\dbD(\scrD_{Q/pQ})$ at the thickening $\grQ:=(\Spec\,Q/pQ\hookrightarrow\Spec\,Q,\delta(p))$ is the triple $(M_Q,\phi_{0,Q},\nabla)$. 

\smallskip
{\bf (b)} There exists a unique $p$-divisible group $\scrD$ over $\Spec\,Q$ such that the evaluation of the filtered Dieudonn\'e crystal of $\scrD$  at the thickening $\grQ$ is the quadruple $(M_Q,F^1_Q,\phi_{0,Q},\nabla)$.}

\medskip
\proof
We first prove (a) and (b) for $p>2$ and $R$ a smooth $W(k)$-algebra. For $p>2$ there exists a unique finite, flat, commutative group scheme $\scrD_n$
over $\Spec\,Q^\wedge_n$ of $p$-power rank such that the object $\dbD(\scrD_n)$ of $\scrM\scrF_{[0,1]}^\nabla(Q_n)$ is $(\scrM_{Q_n/p^nQ_n},\nabla_n)=(M_{Q_n/p^nQ_n},F^1_{Q_n/p^nQ_n},\break\phi_{0,Q_n},\phi_{1,Q_n},\nabla_n)$, cf. Theorem 2.3.4. Due to the universal property enjoyed by $\nabla_n$ (cf. Theorem 3.2 (a)), we can identify $\dbD((\ell^{(n)})^*(\scrD_n))=\dbD(\scrD_{n+1}[p^n])$. Thus we can also identify $(\ell^{(n)})^*(\scrD_n)=\scrD_{n+1}[p^n]$, cf. Theorem 2.3.4. Therefore there exists a unique $p$-divisible group $\scrD$
over $\Spec\,Q$ such that  for all $n\in\dbN^{\ast}$ we have $\scrD[p^n]=\scrD_{n,Q}=\scrD_{n+1,Q}[p^n]$. 

The evaluation at the thickening $\grQ$ of the filtered Dieudonn\'e crystal of $\scrD$ is the quadruple $(M_Q,F^1_Q,\phi_{0,Q},\nabla)$ (we recall from Section 2 that we suppress the Verschiebung maps of such evaluations). As $Q/pQ$ has a finite $p$-basis, from [BM, Prop. 1.3.3] we get that each (filtered) $F$-crystal on $CRIS(\Spec\,Q/pQ/\Spec\,W(k))$ is uniquely determined by its evaluation at the thickening $\grQ$. Therefore $\scrD_{Q/pQ}$ is uniquely determined by $(M_Q,\phi_{0,Q},\nabla)$, cf. [BM, Thm. 4.1.1]. Thus, if $R$ is a smooth $W(k)$-algebra, then (a) holds for $p>2$ and from Grothendieck--Messing deformation theory we get that (b) also holds for $p>2$.

Next we include two extra ways of proving (a) that work for all primes $p\Ge 2$. The first way does not assume that for $p=2$ the condition 3.4 (*) holds and it goes as follows. We work with an arbitrary prime $p\Ge 2$. The existence and the uniqueness of $\scrD_{Q/pQ}$ can be deduced from [dJ, Main Thm. 1]. Strictly speaking, loc. cit. is stated in a way that applies only to smooth $k$-algebras. But as the field $K_{Q/pQ}$ has a finite $p$-basis, loc. cit. applies to show that $\scrD_{K_{Q/pQ}}$ exists and is unique. Descent and extension arguments as in [dJ, Subsect. 4.4] show that $\scrD_{Q/pQ}$ itself exists and it is unique.

We describe with full details the second way (as it is the simplest). Again we work with an arbitrary prime $p\Ge 2$ but if $p=2$ we do assume that the condition 3.4 (*) holds. 

Let $\grA$ be the set of points of $\Spec\,Q/pQ$ whose residue fields are algebraic extensions of $k$. Let $Q_0$ be the localization of $Q$ with respect to a point $\grp_0\in \grA$. The residue field $k_0$ of $\grp_0$ is an algebraic extension of $k$ and thus it is a perfect field. The ring $S_0:=Q_0/pQ_0$ has a finite $p$-basis as $Q/pQ$ does. Let $Q_0^{\text{h}}$ and $\hat Q_0$ be the henselization and the completion of $Q_0$ (respectively). Thus $S_0^{\text{h}}:=Q_0^{\text{h}}/pQ_0^{\text{h}}$ is the henselization of $S_0$ (this follows easily from [BLR, Ch. 2, Sect. 2.3, Prop. 4]). Let $\hat S_0:=\hat Q_0/p\hat Q_0$. Let $x_1,\ldots,x_m\in\hat Q_0$ be such that we can identify $\hat Q_0=W(k_0)[[x_1,\ldots,x_m]]$. Let $\Phi^{\prime}_{\hat Q_0}$ be the Frobenius lift of $\hat Q_0$ that is compatible with $\sigma_{k_0}$ and that takes $x_j$ to $x_j^p$ for all $j\in\{1,\ldots,m\}$. 

Let $(M_{\hat Q_0},F^1_{\hat Q_0},\phi^{\prime}_{0,\hat Q_0},\nabla_0)$ be the extension of $(M_Q,F^1_Q,\phi_{0,Q},\nabla)$ via $Q\hookrightarrow\hat Q_0$ with $\phi^{\prime}_{0,\hat Q_0}$ as a $\Phi^{\prime}_{\hat Q_0}$-linear map. The connection $\nabla_0$ on $M\otimes_{W(k)} \hat Q_0$ is uniquely determined by the equality $\nabla_0\circ\phi^{\prime}_{0,\hat Q_0}=(\phi^{\prime}_{0,\hat Q_0}\otimes d\Phi^{\prime}_{\hat Q_0})\circ\nabla_0$. We recall what $\phi^{\prime}_{0,\hat Q_0}$ is. 

The $\hat Q_0$-linear isomorphism $\phi^{\prime}_{\hat Q_0}:\tilde M_{\hat Q_0}\arrowsim M_{\hat Q_0}$ defined naturally by $\phi^{\prime}_{0,\hat Q_0}$ is the composite of a correction $\hat Q_0$-linear automorphism 
$$\grK:\tilde M_{\hat Q_0}\arrowsim \tilde M_{\hat Q_0}$$ 
with the $\hat Q_0$-linear isomorphism $\phi_{\hat Q_0}:\tilde M_{\hat Q_0}\arrowsim M_{\hat Q_0}$. For $j\in\{1,\ldots,m\}$ let $s_j:=\Phi^{\prime}_{\hat Q_0}(z_j)-z_j^p\in p\hat Q_0$. For $x\in M+{1\over p}F^1$ we have (cf. [De1, Formula (1.1.3.4)])
$$\grK(x\otimes 1):=\sum_{i_1,i_2,\ldots,i_m\in\dbN} (\prod_{j=1}^m \nabla({d\over dz_j})^{i_j})(x\otimes 1)\prod_{j=1}^m {{s_j^{i_j}}\over {i_j!}}.\leqno (20)$$ 
\indent
The reduction of $(M_{\hat Q_0},F^1_{\hat Q_0},\phi^{\prime}_{0,\hat Q_0},(t_{\alpha})_{\alpha\in\scrJ})$ modulo the ideal $(x_1,\ldots,x_m)$ of $\hat Q_0$, is a quadruple of the form 
$$(M_{W(k_0)},F^1_{W(k_0)},h^{\prime}(\phi\otimes\sigma_{k_0}),(t_{\alpha})_{\alpha\in\scrJ})$$ 
for some element $h^{\prime}\in \pmb{GL}_M(W(k_0))$. The triple $(M_{W(k_0)},F^1_{W(k_0)},h^{\prime}(\phi\otimes\sigma_{k_0}))$ is the filtered Dieudonn\'e module of a uniquely determined $p$-divisible group $D(h^{\prime})$ over $\Spec\,W(k_0)$ (for $p=2$, cf. Proposition 2.2.6 and the fact that the condition 3.4 (*) holds). If $\tilde G=G$, then $\nabla$ annihilates each tensor $t_{\alpha}$ and thus from Formula $(20)$ we get that $\grK$ fixes each $t_{\alpha}$ viewed as a tensor of $\scrT(\tilde M_{\hat Q^0})[{1\over p}]$; therefore $h^\prime$ fixes each $t_{\alpha}$ with $\alpha\in\scrJ$ and thus $h^{\prime}\in G(W(k_0))$.

There exists a unique $p$-divisible group $\scrD_{\hat Q_0}$ over $\Spec\,\hat Q_0$ that lifts $D(h^{\prime})$ and such that the evaluation of its filtered Dieudonn\'e crystal at the thickening $(\Spec\,\hat Q_0/p\hat Q_0\hookrightarrow \Spec\,\hat Q_0,\delta(p))$ is $(M_{\hat Q_0},F^1_{\hat Q_0},\phi^{\prime}_{0,\hat Q_0},\nabla_0)$ with respect to the Frobenius endomorphism $\Phi_{\hat Q_0}^\prime$ of $\hat Q_0$ or is $(M_{\hat Q_0},F^1_{\hat Q_0},\phi_{0,\hat Q_0},\nabla)$ with respect to the Frobenius endomorphism $\Phi_{\hat Q_0}$ of $\hat Q_0$ (here and below we denote also by $\nabla$ its natural extensions). The argument for this (working with $\Phi_{\hat Q_0}^\prime$) goes as follows. The existence of $\scrD_{\hat Q_0}$ is implied by [Fa2, Thm. 10]. The uniqueness of $\scrD_{\hat Q_0}\times_{\Spec\,\hat Q_0} \Spec\,\hat Q_0/p\hat Q_0$ is implied by [BM, Prop. 1.3.3 and Thm. 4.1.1]. Thus the uniqueness of $\scrD_{\hat Q_0}$ itself follows from Grothendieck--Messing deformation theory and the fact that the natural divided power structure of the ideal $p(x_1,\ldots,x_m)$ of $\hat Q_0$ is nilpotent modulo $p(x_1,\ldots,x_m)^q$ for all $q\in\dbN^{\ast}$. We recall from Section 1 that $\scrD_{\hat S_0}:=\scrD_{\hat Q_0}\times_{\Spec\,\hat Q_0} \Spec\,\hat S_0$. In the next two paragraphs we will use descent in order to show that $\scrD_{\hat S_0}$ is defined over $\Spec\,S_0^{\text{h}}$.

Let $\hat S_0^{(p)}$ be $\hat S_0$ but viewed as a $\hat S_0$-algebra via the Frobenius endomorphism of $\hat S_0$. The $\hat S_0\otimes_{S_0^{\text{h}}} \hat S_0$-module $\hat S_0^{(p)}\otimes_{{S_0^{\text{h}}}^{(p)}} \hat S_0^{(p)}$ is free and has $\{\prod_{j=1}^m z_j^{n_j}\otimes\prod_{j=1}^m z_{j}^{m_j}|j\in\{1,\ldots,m\},n_j,m_j\in\{0,\ldots,p-1\}\}$ as a $\hat S_0\otimes_{S_0^{\text{h}}} \hat S_0$-basis. Thus $\hat S_0\otimes_{S_0^{\text{h}}} \hat S_0$ has a finite $p$-basis with $m^2$ elements. The normalization of $S_0^{\text{h}}$ in $\hat S_0$ is $S_0^{\text{h}}$ itself and thus $\hat S_0\otimes_{S_0^{\text{h}}} \hat S_0$ is an integral domain. We check that $\hat S_0\otimes_{S_0^{\text{h}}} {}_{}\hat S_0$ is normal. We write $\hat S_0$ as an inductive limit $\hat S_0=\text{lim.}\text{ind.}_{\diamondsuit\in\dbI} {}_{}S_{0\diamondsuit}$ of normal $\hat S_0^{\text{h}}$-algebras of finite type indexed by the set of objects $\dbI$ of a filtered category. Thus $\hat S_0\otimes_{S_0^{\text{h}}} \hat S_0=\text{lim.}\text{ind.}_{\diamondsuit\in\dbI} S_{0\diamondsuit}\otimes_{S_0^{\text{h}}} \hat S_0$. As $R/pR$ is an excellent ring (see [Ma, Ch. 13, Sect. 34]), the homomorphism $S_{0\diamondsuit}\otimes_{S_0^{\text{h}}} \hat S_0\to S_{0\diamondsuit}$ is regular. Thus $S_{0\diamondsuit}\otimes_{S_0^{\text{h}}} \hat S_0$ is normal, cf. [Ma, Ch. 13, Sect. 33, Lemmas 2 and 4]. Therefore $\hat S_0\otimes_{S_0^{\text{h}}} \hat S_0$ is normal. Let $s_1$, $s_2:\Spec\,\hat S_0\otimes_{S_0^{\text{h}}} \hat S_0\to\Spec\,\hat S_0$ be the two natural projection morphisms. 

Both $\dbD(s_1^*(\scrD_{\hat S_0}))$ and $\dbD(s_2^*(\scrD_{\hat S_0}))$ are defined naturally by the triple $(M_Q,\phi_{0,Q},\nabla)$ and thus we have a canonical identification $\dbD(s_1^*(\scrD_{\hat S_0}))=\dbD(s_2^*(\scrD_{\hat S_0}))$. Let 
$$\theta:s_1^*(\scrD_{\hat S_0})\arrowsim s_2^*(\scrD_{\hat S_0})$$ 
be the unique isomorphism such that $\dbD(\theta)$ is this identification $\dbD(s_1^*(\scrD_{\hat S_0}))=\dbD(s_2^*(\scrD_{\hat S_0}))$, cf. [BM, Thm. 4.1.1]. The local rings of $\hat S_0\otimes_{S_0^{\text{h}}}\hat S_0\otimes_{S_0^{\text{h}}} \hat S_0$ are normal and have finite $p$-bases (this is argued as for $\hat S_0\otimes_{S_0^{\text{h}}}\hat S_0$). Thus based on loc. cit. we get naturally a descent datum on $\scrD_{\hat S_0}$ with respect to the faithfully flat morphism $\Spec\,\hat S_0\to\Spec\,S_0^{\text{h}}$ (i.e., $\theta$ satisfies the cocycle condition $s_{23}^*(\theta)\circ s_{12}^*(\theta)=s_{13}^*(\theta)$). Thus standard descent of coherent sheaves of modules (see [BLR, Ch. 6, Sect. 6.1, Thm. 4]) applied to the finite $\Spec\,\hat S_0$-scheme $\scrD_{\hat S_0}[p^n]$ and to the evaluation of $\dbD(\scrD_{\hat S_0}[p^n])$ at the thickening $(\Spec\,\hat S_0\hookrightarrow\Spec\,\hat Q_0/p^n\hat Q_0,\delta(p))$, shows that $\scrD_{\hat S_0}$ is the pull-back of a $p$-divisible group $\scrD_{S_0^{\text{h}}}$ over $\Spec\,S_0^{\text{h}}$ whose Dieudonn\'e crystal is uniquely determined by its evaluation $(M_{Q_0^{\text{h}\wedge}},\phi_{0,Q_0^{\text{h}\wedge}},\nabla)$ at the thickening $(\Spec\,S_0^{\text{h}}\hookrightarrow\Spec\,Q_0^{\text{h}\wedge},\delta(p))$. 

As in the previous paragraph we argue that there exists a natural descent datum on $\scrD_{S_0^{\text{h}}}$ with respect to the faithfully flat (pro-\'etale) morphism $\Spec\,S_0^{\text{h}}\to\Spec\,S_0$ and that $\scrD_{S_0^{\text{h}}}$ is the pull-back of a $p$-divisible group $\scrD_{S_0}$ over $\Spec\,S_0$ whose Dieudonn\'e crystal is uniquely determined by its evaluation $(M_{Q_0^\wedge},\phi_{0,Q_0},\nabla)$ at the thickening $(\Spec\,S_0\hookrightarrow\Spec\,Q_0^\wedge,\delta(p))$. 

Each point of $\Spec\,Q/pQ$ specializes to a point in $\grA$, cf. Theorem 3.2 and our hypotheses. Based on [BM, Thm. 4.1.1] we get that the $\scrD_{Q_0/pQ_0}$'s glue together to define a $p$-divisible group $\scrD_{Q/pQ}$ over $\Spec\,Q/pQ$ whose Dieudonn\'e crystal is uniquely determined by its evaluation $(M_Q,\phi_{0,Q},\nabla)$ at the thickening $\grQ$. Thus (a) holds. If $p>2$, then as above we argue that (b) follows from (a) and the Grothendieck--Messing deformation theory. 

We are left to prove (b) in the case when $p=2$ and the condition 3.4 (*) holds. It suffices to prove the existence and the uniqueness of the lift $\scrD_{Q/p^2Q}$ of $\scrD_{Q/pQ}$ to $\Spec\,Q/p^2Q$. We fix a lift $\scrD^\prime_{Q/p^2Q}$ of $\scrD_{Q/pQ}$ to $\Spec\,Q/p^2Q$, cf. [Il, Thm. 4.4 a) and f)]. Let $\delta(p)^{\text{tr}}$ be the trivial divided power structure of the ideal $(p)$ of $Q/p^2Q$ defined by the identities $(p)^{[s]}=0$, $s\in\dbN^{\ast}\setminus\{1\}$. The $W(k)$-module $\tilde F^{-1}(\End_{W(k)}(M))$ is the maximal direct summand of $\End_{W(k)}(M)$ on which $\dbG_m$ acts via $\mu$ as the identity character of $\dbG_m$ (i.e., acts with weight $1$), cf. Subsection 2.6. Let $\scrL_{\text{crys-lift}}$ (resp. $\scrL_{\text{lift}}$) be the free $Q/pQ$-module of lifts of $F^1_{Q/pQ}$ to direct summands of $M_{Q/p^2Q}$, the zero element corresponding to the Hodge filtration defined by $\scrD^\prime_{Q/p^2Q}$ and by the divided power structure $\delta(p)$ (resp. $\delta(p)^{\text{tr}}$) of the ideal $(p)$ of $Q/p^2Q$. The $Q/pQ$-module structure of $\scrL_{\text{crys-lift}}$ (resp. of $\scrL_{\text{lift}}$) is defined naturally by identifying $\scrL_{\text{crys-lift}}$ (resp. $\scrL_{\text{lift}}$) with the set of images of the lift of $M_{Q/p^2Q}$ that defines the zero element of $\scrL_{\text{crys-lift}}$ (resp. of $\scrL_{\text{lift}}$) through elements of the form $1_{M_{Q/p^2Q}}+p\upsilon\in \pmb{GL}_M(Q/p^2Q)$, where $\upsilon\in \tilde F^{-1}(\End_{W(k)}(M))\otimes_{W(k)} Q/p^2Q$. Let $L\in\scrL_{\text{crys-lift}}$ be such that it corresponds to $F^1_{Q/p^2Q}$. 

We define a natural map of sets 
$$\scrN_{Q/pQ}:\scrL_{\text{lift}}\to\scrL_{\text{crys-lift}}$$ 
as follows. Let $x\in\scrL_{\text{lift}}$ and let $\scrD^x_{Q/p^2Q}$ be the lift of $\scrD_{Q/pQ}$ defined by $x$ and Grothendieck--Messing deformation theory (the divided power structure of the ideal $(p)$ of $Q/p^2Q$ being $\delta(p)^{\text{tr}}$). We define $\scrN_{Q/pQ}(x)$ to be the Hodge filtration of $M_{Q/p^2Q}$ defined by $\scrD^x_{Q/p^2Q}$ using the divided power structure $\delta(p)$ of the ideal $(p)$ of $Q/p^2Q$. 

The map $\scrN_{Q/pQ}$ has a functorial aspect with respect to pull-backs of $\scrD_{Q/pQ}$. Let 
$$\grT:\Spec\,W(\tilde k)\to \Spec\,\hat Q_0$$ 
be a Teichm\"uller lift with respect to $\phi_{\hat Q_0}$ whose special fibre is dominant. Here $\tilde k$ is a big enough perfect field that contains the field of fractions $K_{\hat S_0}$ of $\hat S_0$. The map $\scrN_{\tilde k}$ which is the analogue of $\scrN_{Q/pQ}$ but obtained using lifts of $\scrD_{\tilde k}$ to $W_2(\tilde k)$, is injective (see proof of Proposition 2.2.6). Thus $\scrN_{Q/pQ}$ is injective. Thus, if $\scrD_{Q/p^2Q}$ exists, then it is unique. To $\scrD_{\hat Q_0}\times_{\Spec\,\hat Q_0} \Spec\,\hat Q_0/p^2\hat Q_0$ corresponds an element $L_0\in\scrL_{\text{lift}}\otimes_{Q/pQ}\hat S_0$. The images of $L_0$ in $\scrL_{\text{lift}}\otimes_{Q/pQ} \hat S_0\otimes_{S_0^{\text{h}}} \hat S_0$ via $s_1$ and $s_2$ are equal and thus we have $L_0\in\scrL_{\text{lift}}\otimes_{Q/pQ} S_0^{\text{h}}$, cf. [BLR, Ch. 6, Sect. 6.1, Lem. 2]. Working with the pro-\'etale morphism $\Spec\,S_0^{\text{h}}\to\Spec\,S_0$, a similar argument shows that $L_0\in\scrL_{\text{lift}}\otimes_{Q/pQ} S_0$. From the injectivity of $\scrN_{Q/pQ}$ and of its analogue ``localization" $\scrN_{S_0}$ (thought to be $\scrN_{Q/pQ}\otimes 1_{S_0}$) and from the fact that $\scrN_{S_0}(L_0)=L$, we get that $L_0\in\scrL_{\text{lift}}\otimes_{Q/pQ} \cap_{\grp_0\in \grA} S_0$; here we use quotations for ``localization" as we will not stop to check that the map $\scrN_{Q/pQ}$ is indeed $Q/pQ$-linear. As each point of $\Spec\,Q/pQ$ specializes to a point in $\grA$, we have $Q/pQ=\cap_{\grp_0\in \grA} S_0$. Thus $L_0\in\scrL_{\text{lift}}$. The element $\scrN_{Q/pQ}(L_0)-L\in\scrL_{\text{crys-lift}}$ is $0$ as this is so after tensorization with $\hat S_0$ over $Q/pQ$. Thus $\scrN_{Q/pQ}(L_0)=L$ i.e., $\scrD_{Q/p^2Q}$ exists. This ends the proof of (b).\endproof

\medskip\noindent
{\bf 3.4.2. \'Etale Tate-cycles.} We continue to assume that the hypotheses of Theorem 3.4.1 hold.  Let the $p$-divisible group $\scrD$ be as in Theorem 3.4.1 (b). 

We will also assume that we have $\phi_R(z_j)=z_j^p$ for all $j\in\{1,\ldots,m\}$ and that $\tilde G=G$; thus $\tau\in G(R^\wedge)$ and $\nabla$ annihilates $t_{\alpha}\in\scrT(M_Q)[{1\over p}]$ for all $\alpha\in\scrJ$. 

We consider the $B^+(Q)$-linear monomorphism
$$i_{\scrD}:M_{B^+(Q)}=M_Q\otimes_Q {}_{i_Q} B^+(Q)\to H^1(\scrD_{K_Q})\otimes_{\dbZ_p} B^+(Q)$$ constructed as in (3), the Frobenius endomorphism of $M_{B^+(Q)}=M_Q\otimes_Q {}_{i_Q} B^+(Q)$ being $\phi_{0,Q}\otimes\Phi_{Q/pQ}$ (here $i_Q:Q\hookrightarrow B^+(Q)$ and $\Phi_{Q/pQ}$ are as in Subsection 2.3 and Subsubsection 2.3.5 but for $Q$; thus the $W(k)$-monomorphism $i_Q$ is compatible with the Frobenius lifts $\Phi_Q$ and $\Phi_{Q/pQ}$). Let $\scrV_{\alpha}\in\scrT(H^1(\scrD_{K_Q})[{1\over p}])\otimes_{\dbQ_p} B^+(Q)[{1\over {\beta}}]$ correspond to $t_{\alpha}$ via the $B^+(Q)[{1\over {\beta}}]$-linear isomorphism $i_{\scrD}[{1\over {\beta}}]$. We check that $\scrV_{\alpha}\in\scrT(H^1(\scrD_{K_Q})[{1\over p}])$. 

Let $\grT:\Spec\,W(\tilde k)\to \Spec\,\hat Q_0$ be as in the last paragraph of the proof of Theorem 3.4.1. We denote by $\grT_1$ the $W(k)$-homomorphism $B^+(Q)\to B^+(W(\tilde k))$ defined naturally by $\grT$ and the choice of a $W(k)$-homomorphism $\grT_0:\bar Q\otimes_Q W(\tilde k)\to\overline{W(\tilde k)}$ (see Subsection 2.3 and the end  of Subsection 2.1 for $\bar Q$ and $\overline{W(\tilde k)}$). As $\grT$ is dominant, the following restriction $\grT_{00}:\bar Q\to \overline{W(\tilde k)}$ of $\grT_0$ is injective. Based on Lemma 2.3.1 (a) and (b) we easily get that the $W(k)$-homomorphism $W(A_{Q/pQ})\to W(A_{\tilde k})$ defined naturally by $\grT_{00}$ is also injective. This implies that the $W(k)$-homomorphism $\grT_1[{1\over {\beta}}]:B^+(Q)[{1\over {\beta}}]\to B^+(W(\tilde k))[{1\over {\beta}}]$ is also injective. The image of $\scrV_{\alpha}$ via $\grT_1[{1\over {\beta}}]$ is a tensor of $\scrT(H^1(\scrD_{B(\tilde k)})[{1\over p}])$, cf. end of Subsubsection 2.2.4. As we have a canonical identification $H^1(\scrD_{K_Q})[{1\over p}]=H^1(\scrD_{B(\tilde k)})[{1\over p}]$ of $\dbZ_p$-modules, the relation $\scrV_{\alpha}\in\scrT(H^1(\scrD_{K_Q})[{1\over p}])$ is implied by the injectivity of $\grT_1[{1\over {\beta}}]$. The tensor $\scrV_{\alpha}\in\scrT(H^1(\scrD_{K_Q})[{1\over p}])$ is fixed by $\Gal(K_{\bar Q}/K_Q)$ (as $t_{\alpha}$ is so) and thus it is an \'etale Tate-cycle of $\scrD$ (more precisely, of $\scrD_{Q[{1\over p}]}$).

We assume now that there exists a section $a:\Spec\,W(k)\to Y^\wedge$ such that $\tau\circ a\in G(W(k))$ is the identity section and the ideal $I_a$ that defines $a$ contains the ideal $(z_1,\ldots,z_m)$ of $R^\wedge$. Let $y:\Spec\,k\hookrightarrow Y^\wedge$ be the closed point defined by $a$. As $d\Phi_R/p(dz_j)=z_j^{p-1}dz_j$ is $0$ modulo the ideal $J_y$, the $\sigma$-linear map $\omega_y$ of the Subsection 3.3 is $0$ and thus $c_y=0$. Let $a_{\infty}:\Spec\,W(k)\to\Spec\,Q_{\infty}^\wedge$ be the unique morphism such that its composite with the morphism $\Spec\,Q_{\infty}^\wedge\to Y^\wedge$ is $a$, cf. Theorem 3.3.1 applied with $\tilde G=G$ and $c_y=0$. We will choose $\Spec\,Q$ to be the unique connected component of $\Spec\,Q_{\infty}^\wedge$ through which $a_{\infty}$ factors; we denote also by $a_{\infty}:\Spec\,W(k)\to\Spec\,Q$  the resulting factorization. The $\Spec\,k$-scheme $\Spec\,Q/pQ$ is geometrically connected. 

As $a_{\infty}:\Spec\,W(k)\to\Spec\,Q$ is a Teichm\"uller lift with respect to $\Phi_Q$, the pull-back of $(\scrD,(t_{\alpha})_{\alpha\in\scrJ})$ (resp. of $(\scrD,(\scrV_{\alpha})_{\alpha\in\scrJ})$) via $a_{\infty}$ is $(D,(t_{\alpha})_{\alpha\in\scrJ})$ (resp. is $(D,(v_{\alpha})_{\alpha\in\scrJ})$). 

\medskip\noindent
{\bf 3.4.3. Lemma.} {\it We work under all assumptions of Subsubsection 3.4.2. Then to prove the Main Theorem for $(D,(t_{\alpha})_{\alpha\in\scrJ})$ it is enough to prove the Main Theorem for an arbitrary pull-back of $(\scrD,(t_{\alpha})_{\alpha\in\scrJ})$ via a $W(k)$-valued point of $\Spec\,Q$.}

\medskip
\proof
As $Q[{1\over p}]$ is an integral domain and $\scrD_{\bar Q[{1\over p}]}$ is isomorphic to $(\dbQ_p/\dbZ_p)^{d_M}$, all pull-backs of $(\scrD,(\scrV_{\alpha})_{\alpha\in\scrJ})$ via $\overline{B(k)}$-valued points of $\Spec\,Q$ are isomorphic to $(H^1(D),(v_{\alpha})_{\alpha\in\scrJ})$. This is so due to the last sentence before the  lemma and the fact that all the mentioned pull-backs can be canonically identified with the pull-back of $(\scrD,(\scrV_{\alpha})_{\alpha\in\scrJ})$  via the natural dominant morphism $\Spec\,K_{\bar Q}\to \Spec\,Q$. 

Moreover each pull-back of $(\scrD,(t_{\alpha})_{\alpha\in\scrJ})$ via a $W(k)$-valued point of $\Spec\,Q$ is a pair $(D(h^{\prime}),(t_{\alpha})_{\alpha\in\scrJ})$, where the $p$-divisible group $D(h^{\prime})$ over $\Spec\,W(k)$ has a filtered Dieudonn\'e module of the form $(M,F^1,h^{\prime}\phi)$ for some element $h^{\prime}\in G(W(k))$ (see the paragraph after Formula $(20)$ of the proof of Theorem 3.4.1, applied with $k_0=k$). 

The lemma follows from the last two paragraphs.\endproof

\medskip\noindent
{\bf 3.4.4. Remarks.} {\bf (a)} Suppose that $p=2$ and that the condition 3.4 (*) does not hold. One can use Artin approximation theory to show that $\scrD$ always exists, provided modulo $2$ we work in the \'etale topology of $\Spec\,Q/2Q$. In general, $\scrD$ is not unique and there exists nothing to guarantee that we can work with a single connected \'etale cover of $\Spec\,Q/2Q$. 

\smallskip
{\bf (b)} If the condition 3.4 (*) holds, then we have $m_y^{\tilde G}(-1)=0$ (cf. inequality (16)) and thus $Q=R^\wedge$ (cf. Corollary 3.3.2).

\medskip\smallskip\noindent
{\bf 3.5. Application.} Until Section 4 we assume that $k=\bar k$. Let $u=1_M+v,g=u^{-1}=1_M-v\in U_0(W(k))$ be as in Subsection 2.7. We take $m=1$, $z=z_1$, $R=W(k)[z]$, $\Phi_R(z)=z^p$, $\tilde G=G$, and 
$$\tau=1_{M_{R^\wedge}} +(g-1_{M_{R^\wedge}})\otimes z=1_{M_{R^\wedge}}-v\otimes z\in U_0(R^\wedge)\leqslant G(R^\wedge);$$ 
thus $\tau$ modulo $z$ is $1_M$ and $\tau$ modulo $(z-1)$ is $g$. From Lemma 2.7.1 (a) we get that if $p=2$, then condition 3.4 (*) holds if and only if either $D$ or $D^{\text{t}}$ is connected. Let $a:\Spec\,W(k)=\Spec\,R^\wedge/(z)\hookrightarrow \Spec\,R^\wedge=Y^\wedge$. Let $y:\Spec\,k\hookrightarrow \Spec\,R^\wedge=Y^\wedge$ be defined by $a$. 

As $d\Phi_R/p(dz)=z^{p-1}dz$, we have $\omega_y=0$ and $c_y=0$. Let the connected component $\Spec\,Q$ of $\Spec\,Q_{\infty}^\wedge$ be chosen as in Subsubsection 3.4.2 and let $a_{\infty}:\Spec\,W(k)\to\Spec\,Q$ be as in Subsubsection 3.4.2.
Thus all the assumptions of Subsubsection 3.4.2 hold if and only if, in the case when $p=2$, either $D$ or $D^{\text{t}}$  is connected.

For $\tilde y\in Y^\wedge(k)=k$ with $\tilde y\neq y$ (i.e., with $\tilde y\neq 0$), we have $c_{\tilde y}=1$. Thus from Corollary 3.3.2 and Lemma 2.7.2 we get that the morphism $\Spec\,Q/pQ[{1\over z}]\to \Spec\,R/pR[{1\over z}]$ induced by $\ell[{1\over z}]$ is a pro-\'etale cover. Thus the reduction modulo $p$ of $\ell:\Spec\,Q\to\Spec\,R^\wedge$ is onto. Thus there exists a section $b_{\infty}:\Spec\,W(k)\to \Spec\,Q$ such that the section $b:=\ell\circ b_{\infty}:\Spec\,W(k)\to Y^\wedge$ is defined by the ideal $(z-1)$ of $R^\wedge$. As $\Phi_R(z-1)=z^p-1\in (z-1)$, $b_{\infty}$ is a Teichm\"uller lift with respect to $\Phi_Q$. Thus the pull-back of $(M_{Q},F^1_Q,\phi_Q,(t_{\alpha})_{\alpha\in\scrJ})$ via $b_{\infty}$ is $(M,F^1,g\phi,(t_{\alpha})_{\alpha\in\scrJ})$. This achieves the connection of $(M,F^1,\phi,G,(t_{\alpha})_{\alpha\in\scrJ})$ with $(M,F^1,g\phi,G,(t_{\alpha})_{\alpha\in\scrJ})$ via a global deformation over $\Spec\,Q$ with $\Spec\,Q/pQ[{1\over z}]$ a suitable connected, pro-\'etale cover of $\Spec\,k[z][{1\over z}]$.

\medskip\noindent
{\bf  3.5.1. Theorem.} {\it It suffices to prove the Main Theorem under the extra hypotheses that $G$ is the torus $T$, that $u=g=1_M$, and, in the case when $p=2$, that $D$ is connected.}

\medskip
\proof 
If $p=2$, then we can assume that $D$ is connected (cf. Fact 2.5.4). Thus all the assumptions of Subsubsection 3.4.2 hold. Let $D(g)$ be the $p$-divisible group over $\Spec\,W(k)$ whose filtered Dieudonn\'e module is $(M,F^1,g\phi)$. If $p=2$, then $D(g)_k$ is connected (cf. Lemma 2.7.1 (a)) and thus the existence and the uniqueness of $D(g)$ follow from Proposition 2.2.6. 

The pull-back of $(\scrD,(t_{\alpha})_{\alpha\in\scrJ})$ through $b_{\infty}$ is $(D(g),(t_{\alpha})_{\alpha\in\scrJ})$. Thus to prove the Main Theorem for $(D,(t_{\alpha})_{\alpha\in\scrJ})$ is equivalent to proving the Main Theorem for $(D(g),(t_{\alpha})_{\alpha\in\scrJ})$, cf. Lemma 3.4.3. If $T$ and $(t_{\alpha})_{\alpha\in\scrJ_T}$ are as in Subsection 2.7 and if the Main Theorem holds for $(D(g),(t_{\alpha})_{\alpha\in\scrJ_T})$, then the Main Theorem also holds for $(D(g),(t_{\alpha})_{\alpha\in\scrJ})$ and thus for $(D,(t_{\alpha})_{\alpha\in\scrJ})$. As $T_{B(k)}$ is the subgroup of $\pmb{GL}_{M[{1\over p}]}$ that fixes $t_{\alpha}$ for all $\alpha\in\scrJ_T$ (cf. property 2.7 (vi)), we get that to prove the Main Theorem for $(D,(t_{\alpha})_{\alpha\in\scrJ})$ we can also assume that $G=T$ and $g=u=1_M$.\endproof

\bigskip\smallskip
\noindent
{\boldsectionfont 4. Proof of the Main Theorem}
\bigskip\smallskip 

In Subsections 4.1 and 4.2 we prove the Main Theorem. Example 4.3 illustrates the computations of Subsections 4.1 and 4.2. Subsection 4.4 proves Corollary 1.4. Corollary 4.5 pertains to a $\dbZ_p$ form of the Main Theorem. We will use the notations of Subsection 2.1 and of Subsubsections 2.2.1 to 2.2.4. Until Section 5 we assume that $k=\bar k$. 

\medskip\smallskip\noindent
{\bf 4.1. Notations and simple properties.}
We start the proof of the Main Theorem. To prove the Main Theorem we can assume that $G$ is a torus $T$ and, in the case when $p=2$, that $D$ is connected (cf. Theorem  3.5.1). This case follows by combining Lemma 4.1.1 with the recent result [Ki3, Cor. (1.4.3)] that relies heavily on the new crystalline machinery developed in [Ki1] and [Ki2]. 

To be self-contained, we present a proof of the Main Theorem in the case $G=T$ which provides more information, which is shorter and simpler, which relies only on the crystalline results of Subsection 2.2, and which precedes the proof of [Ki3, Cor. (1.4.3)].

Let $\mu=\mu_{\text{can}}^{-1}:\dbG_m\to T=G$ be as in Subsection 2.1. Let $\sigma_{\phi}$, $M_{\dbZ_p}$, and $G_{\dbZ_p}$ be as before Lemma 2.5.2. Let $T^{\acute{et}}:=G^{\acute{et}}$ and let $T_{\dbZ_p}:=G_{\dbZ_p}$. Let $T_{0,\dbZ_p}$ be the smallest subtorus of $T_{\dbZ_p}$ such that $\mu:\dbG_m\to T$ factors through $T_0:=T_{0,W(k)}$. Thus $T_0$ is the subtorus of $T$ generated by the images of the conjugates of $\mu$ under integral powers of $\phi$. Based on Subsection 2.7 and the proof of Theorem 3.5.1 we can assume that $T=T_0$ but for the sake of flexibility, in what follows it is more convenient not to assume this. 

Let $\scrE$ (resp. $\scrE^{\acute{et}}$) be the $\dbZ_p$-subalgebra of endomorphisms of $M_{\dbZ_p}$ (resp. of $H^1(D)$ and thus also of $D^{\text{t}}$) fixed by $T_{\dbZ_p}$ (resp. by $T^{\acute{et}}$). As $T_{\dbZ_p}$ is a torus, $\scrE$ is a semisimple $\dbZ_p$-algebra. Let $T_{1,\dbZ_p}$ be the double centralizer of $T_{\dbZ_p}$ in $\pmb{GL}_{M_{\dbZ_p}}$ (i.e., the centralizer of $\scrE$ in $\pmb{GL}_{M_{\dbZ_p}}$). Thus $T_1:=T_{1,W(k)}$ is the double centralizer of $T$ in $\pmb{GL}_M$. Both $T_{1,\dbZ_p}$ and $T_1$ are tori. Let $T_{2,\dbZ_p}$ be a maximal torus of $\pmb{GL}_{M_{\dbZ_p}}$ that contains $T_{1,\dbZ_p}$. Thus $T_2:=T_{2,W(k)}$ is a maximal torus of $\pmb{GL}_M$ that contains $T_1$ and such that we have $\phi(\Lie(T_2))=\Lie(T_2)$. Let $n$ and $m$ be the ranks of $F^0$ and $F^1$ (respectively) i.e., the codimension and the dimension (respectively) of $D$. 

Let $\scrB:=\{a_1,\ldots,a_{n+m}\}$ be a $W(k)$-basis for $M=F^1\oplus F^0$ such that we have inclusions $\{a_1,\ldots,a_n\}\subseteq F^0$ and $\{a_{n+1},\ldots,a_{n+m}\}\subseteq F^1$ and moreover $T_2$ normalizes $W(k)x$ for all $x\in\scrB$. For $i\in\{1,\ldots,n\}$ let $\eps_i:=0$. For $i\in\{n+1,\ldots,n+m\}$ let $\eps_i:=1$. As $\phi(\Lie(T_2))=\Lie(T_2)$, there exists a permutation $\pi$
of the set $\{1,\ldots,n+m\}$ such that we have $\phi(a_i)\in \dbG_m(W(k))p^{\eps_i}a_{\pi(i)}$ for all $i\in\{1,\ldots,n+m\}$. By replacing each $a_i$ with a suitable $\dbG_m(W(k))$-multiple of it, we can assume that we have 
$$\phi(a_i)=p^{\eps_i}a_{\pi(i)}\;\;\;\forall i\in\{1,\ldots,n+m\}.$$ 
\noindent
Let $\phi_1:F^1\to M$ be the $\sigma$-linear map such that for all $i\in\{n+1,\ldots,n+m\}$ we have $\phi_1(a_i):=a_{\pi(i)}$. Thus $(M,F^1,\phi,\phi_1)$ modulo $p$ is $\dbD(D[p])$ of Subsubsection 2.3.3. 

\medskip\noindent
{\bf 4.1.1. Lemma.} {\it We have $\scrE^{\acute{et}}=\scrE$ and $T^{\acute{et}}$ is a torus.}

\medskip
\proof
Each element of $\scrE$ (resp. of $\scrE^{\acute{et}}$) when viewed as a tensor of $\scrT(M[{1\over p}])$ (resp. of $\scrT(H^1(D))$) is fixed by both $\phi$ and $\mu$ (resp. by $\Gal(B(k))$). The functorial aspect of Fontaine comparison theory allows us to identify $\scrE^{\acute{et}}[{1\over p}]=\scrE[{1\over p}]$ in such a way that $\scrE^{\acute{et}}$ is a $\dbZ_p$-subalgebra of $\scrE$. We have $\scrE^{\acute{et}}=\scrE^{\acute{et}}[{1\over p}]\cap\End(H^1(D))=\scrE[{1\over p}]\cap\End((M,F^1,\phi))=\scrE$, cf. Lemma 2.2.5. Therefore $\scrE^{\acute{et}}$ is a semisimple $\dbZ_p$-algebra.

The centralizer of $\scrE^{\acute{et}}$ in $\pmb{GL}_{H^1(D)}$ is a torus $T^{\acute{et}}_{1,\dbZ_p}$ that contains $T^{\acute{et}}$. As $i_D[{1\over {\beta}}]$ is an isomorphism, the group schemes $T$ and $T^{\acute{et}}$ are isomorphic over $\Spec\,B^+(W(k))[{1\over {\beta}}]$; thus $T^{\acute{et}}_{\dbQ_p}$ is a torus. As $T^{\acute{et}}_{\dbQ_p}$ is a subtorus of $T^{\acute{et}}_{1\dbQ_p}$, the schematic closure $T^{\acute{et}}$ of $T^{\acute{et}}_{\dbQ_p}$ in the torus $T^{\acute{et}}_{1,\dbZ_p}$ is a torus over $\Spec\,\dbZ_p$; to check this, we can work in the \'etale topology of $\Spec\,\dbZ_p$ and thus we can assume that $T^{\acute{et}}_{1,\dbZ_p}$ is split and this case is trivial. Thus $T^{\acute{et}}$ is a torus.\endproof

\medskip\noindent
{\bf 4.1.2. Decomposing $D$.} Let $\pi=\prod_{l\in\scrC(\pi)} \pi_l$ be the decomposition of $\pi$ into cycles. If $l\in\scrC(\pi)$ and $\pi_l=(i_1,\ldots,i_q)$ with $q\in\dbN^{\ast}$, let $M_l:=\oplus_{s=1}^q W(k)a_{i_s}$. We get a direct sum decomposition $(M,F^1,\phi)=\oplus_{l\in\scrC(\pi)} (M_l,M_l\cap F^1,\phi)$. The projection of $M$ on $M_l$ along $\oplus_{l^\prime\in\scrC(\pi)\setminus\{l\}} M_{l^\prime}$ is an element of $\scrE=\scrE^{\acute{et}}$. Thus the direct sum decomposition $(M,F^1,\phi)=\oplus_{l\in\scrC(\pi)} (M_l,M_l\cap F^1,\phi)$ defines a product decomposition  $D=\prod_{l\in\scrC(\pi)} D_l$ into $p$-divisible groups over $\Spec\,W(k)$ whose special fibres have a unique Newton polygon slope. Below we will often use this fact in order to reduce our computations to the simpler case when $\pi$ is a cycle.

\medskip\noindent
{\bf 4.1.3. Lubin--Tate quintuples.} We say that $(M,F^1,\phi,T,(t_{\alpha})_{\alpha\in\scrJ})$ is a Lubin--Tate quintuple if $m=1$ and $\pi$ is a cycle (equivalently, if the $W(k)$-module $F^1$ has rank $1$ and the $F$-isocrystal $(M[{1\over p}],\phi)$ over $k$ is simple of Newton polygon slope ${1\over {n+1}}$).

We check that if $(M,F^1,\phi,T,(t_{\alpha})_{\alpha\in\scrJ})$ is a Lubin--Tate quintuple, then we have $T=G=T_0=T_1=T_2$. The cocharacter $\mu:\dbG_m\to \pmb{GL}_M$ of $(M,F^1,\phi)$ acts trivially on $a_i$ for $i\in\{2,\ldots,n+m\}$ and non-trivially on $a_1$. Thus $T_0$ contains the rank $1$ subtorus of $T_2$ that fixes $a_2, a_3,\ldots$, and $a_{n+m}$. But $\phi$ normalizes $\Lie(T_{0,B(k)})$ (cf. the very definition of $T_0$) and thus $T_0$ contains the subgroup scheme of $T_2$ generated by the images of the conjugates of $\mu$ under powers of $\phi$. Therefore by induction on $i\in\{1,\ldots,n+m\}$ we get that $T_0$ contains the rank $1$ subtorus of $T_2$ that fixes $a_{1+i},\ldots, a_{n+m}$, $a_1,\ldots$, and $a_{i-1}$. Thus $T_2\leqslant  T_0$. Therefore we have $T_0=T=G= T_1=T_2$. 

\medskip\noindent
{\bf 4.1.4. Theorem.} {\it We recall that $G$ is a torus $T$ and, in the case when $p=2$, that $D$ is connected. There exists $t\in T(B^+(W(k))[{1\over {\beta}}])$ that takes the $B^+(W(k))$-submodule $(i_D[{1\over {\beta}}])^{-1}(H^1(D)\otimes_{\dbZ_p} B^+(W(k)))$ of $M\otimes_{W(k)} B^+(W(k))[{1\over {\beta}}]$ onto $M\otimes_{W(k)} B^+(W(k))$.} 

\medskip\smallskip\noindent
{\bf 4.2. Proofs of 1.2 and 4.1.4.} To prove the Main Theorem is equivalent to proving Theorem 4.1.4, cf. Lemma 2.5.2 (b) and Fact 2.5.1 (b). We will prove the Main Theorem and Theorem 4.1.4 in the next eleven subsubsections.  

\medskip\noindent
{\bf 4.2.1. Lemma.} {\it Suppose that $T=T_1$ (for instance, this holds if $(M,F^1,\phi,T,(t_{\alpha})_{\alpha\in\scrJ})$ is a Lubin--Tate quintuple). Then the Main Theorem and Theorem 4.1.4 hold.}

\medskip
\proof
As $T=T_1$, $T$ is the centralizer of $\scrE$ in $\pmb{GL}_M$. Thus to prove the Main Theorem and Theorem 4.1.4 we can assume that $\scrE=\{t_{\alpha}|\alpha\in\scrJ\}$ and $\scrE^{\acute{et}}=\{v_{\alpha}|\alpha\in\scrJ\}$, cf. Fact 2.5.1 (a) and (b). Let $\grE:=\scrE^{\acute{et}}\otimes_{\dbZ_p} W(k)=\scrE\otimes_{\dbZ_p} W(k)$, cf. Lemma 4.1.1. The $W(k)$-algebra $\grE$ is a finite product $\prod_{i\in\scrI} M_{r_i}(W(k))$ of matrix $W(k)$-algebras (here each $r_i\in\dbN^{\ast}$). Each representation of $\prod_{i\in\scrI} M_{r_i}(W(k))$ on a free $W(k)$-module of finite rank is a direct sum indexed by $i\in\scrI$ of a finite number $c_i$ of copies of the standard representation of $M_{r_i}(W(k))$ of rank $r_i$. The representations of $\grE$ on $M$ and $H^1(D)\otimes_{\dbZ_p} W(k)$ involve the same numbers $c_i$, as the tensorizations of these representations with $B^+(W(k))[{1\over {\beta}}]$ are isomorphic (cf. Subsubsection 2.2.4). Thus the representations of $\grE$ on $M$ and $H^1(D)\otimes_{\dbZ_p} W(k)$ are isomorphic i.e., there exists an isomorphism $\rho:(M,(t_{\alpha})_{\alpha\in\scrJ})\arrowsim (H^1(D)\otimes_{\dbZ_p} W(k),(v_{\alpha})_{\alpha\in\scrJ})$. Therefore the Main Theorem and Theorem 4.1.4 hold.\endproof

\medskip\noindent
{\bf 4.2.2. The $w_i$ elements.} Let $t\in T_1(B^+(W(k))\fracwithdelims[]1{\beta})$ be an element that takes the $B^+(W(k))$-module $(i_D[{1\over {\beta}}])^{-1}(H^1(D)\otimes_{\dbZ_p} B^+(W(k)))$ onto $M\otimes_{W(k)} B^+(W(k))$, cf. Lemma 4.2.1 applied to $(M,F^1,\phi,T_1)$ instead of $(M,F^1,\phi,T)$. As $i_D$ is a $B^+(W(k))$-monomorphism and $\beta$ annihilates $\Coker(i_D)$, there exist elements $w_i\in B^+(W(k))\cap \dbG_m(B^+(W(k))\fracwithdelims[]1{\beta})$ such that we have
$$t(a_i)=w_ia_i,\;\;\forall i\in\{1,\ldots,n+m\}.$$
\indent
Let $q_k:B^+(W(k))\twoheadrightarrow V(k)/pV(k)$ be as in Subsubsection 2.2.1. To compute the $q_k(w_i)$'s (see Proposition 4.2.5 below) we need few extra preliminaries. In all that follows we assume that the different roots of $p$ which will show up are powers of a fixed high order root (like the $(p^{m+n})!$-th root) of $p$ and that $q_k(\xi)\in V(k)/pV(k)$ is $p^{1\over p}$ modulo $p$. We consider the system $\scrS$ of equations
$$X_i^p=(-p)^{\eps_{\pi(i)}\over p}X_{\pi(i)},\;\; i\in\{1,\ldots,n+m\}$$
in $n+m$-variables $X_1,\ldots,X_{n+m}$ over $V(k)$. We fix a non-zero solution $(Z_1,\ldots,Z_{n+m})$ of $\scrS$. Let $v:K(k)\setminus\{0\}\to\dbQ$ be the valuation of $K(k)$ normalized by $v(p)=1$. 

\medskip\noindent
{\bf 4.2.3. Lemma.} {\it {\bf (a)} If $\pi$ is a cycle, then any other non-zero solution of $\scrS$ is of the form $$(\eta Z_1,\eta^{p^{\gamma_2}} Z_2,\ldots,\eta^{p^{\gamma_{n+m}}} Z_{n+m}),\leqno (21)$$
where $\eta\in\mu_{p^{n+m}-1}(W(k))\subseteq \dbG_m(V(k))$ and where for $i\in\{2,\ldots,n+m\}$ the number $\gamma_i\in\{1,\ldots,n+m-1\}$ is uniquely determined by the equality $\pi^{\gamma_i}(1)=i$.

\smallskip
{\bf (b)} For $i\in\{1,\ldots,n+m\}$ there exists a rational function $Q_i(x)\in\dbQ(x)$ that depends on $\pi$ and $n$ but not on $p$ and such that we have $v(Z_i)=Q_i(p)\in [0,{1\over {p(p-1)}}]\cap\dbQ$. 

\smallskip
{\bf (c)} If $(i_1,\ldots,i_q)$ is a cycle of $\pi$, then $v(Z_{i_1})=Q_{i_1}(p)$ is equal to ${1\over {p(p-1)}}$ if and only if $\eps_{i_1}=\cdots=\eps_{i_q}=1$.} 

\medskip
\proof
Part (a) is trivial. To prove (b) and (c), let $(i_1,\ldots,i_q)$ be a cycle of $\pi$ of length $q$ and let $i_{q+1}:=i_1$. Let
$$Q_{i_1}(x):={{\sum_{j\in\{1,\ldots,q\},\,\eps_{i_{j+1}}=1} x^{q-j-1}}\over {x^q-1}}.\leqno (22)$$
From the shape of $\scrS$ we get $Z_{i_1}^{p^q}=(-p)^{\sum_{j\in\{1,\ldots,q\},\,\eps_{i_{j+1}}=1} p^{q-j-1}}Z_{i_1}$. As $Z_{i_1}\in V(k)\setminus\{0\}$, we get that $v(Z_{i_1})=Q_{i_1}(p)\le \sum_{j=1}^q {p^{q-j-1}\over {p^q-1}}={1\over {p(p-1)}}$. Thus $Q_{i_1}(p)\in [0,{1\over {p(p-1)}}]\cap\dbQ$ and therefore (b) holds. We have $Q_{i_1}(p)={1\over {p(p-1)}}$ if and only if $\sum_{j\in\{1,\ldots,q\},\,\eps_{i_{j+1}}=1} x^{q-j}=\sum_{j=1}^q x^{q-j}$ i.e., if and only if $\eps_{i_1}=\cdots=\eps_{i_q}=1$. Thus (c) also holds.\endproof 

\medskip\noindent
{\bf 4.2.4. Proposition.} {\it We recall that $i_D^*\colon T_p(D_{B(k)})\otimes_{\dbZ_p} B^+(W(k))\hookrightarrow M^*\otimes_{W(k)} B^+(W(k))$ is the transpose (dual) of the $B^+(W(k))$-linear monomorphism $i_D:M\otimes_{W(k)} B^+(W(k))\hookrightarrow H^1(D)\otimes_{\dbZ_p} B^+(W(k))$ of (1). Let $r_D:M\hookrightarrow H^1(D)\otimes_{\dbZ_p} B^+(W(k))$ and $r_D^*:T_p(D_{B(k)})\hookrightarrow M^*\otimes_{W(k)} B^+(W(k))$ be the natural restrictions of $i_D$ and $i_D^*$ (respectively). Then the $\dbF_p$-linear maps $r_D$ modulo $p$ and $r_D^*$ modulo $p$ are injective.}

\medskip
\proof 
We recall that $\beta$ annihilates $\Coker(i_D)$ and $\Coker(i_D^*)$ and that $q_k(\beta)$ is $p^{1\over {p-1}}$ modulo $p$ times a unit of $V(k)/pV(k)$ (see Subsubsections 2.2.1 and 2.2.3). Thus if $p>2$, then the $\dbF_p$-linear maps $r_D$ modulo $p$ and $r_D^*$ modulo $p$ are injective. 

We check that $r_D$ modulo $p$ is injective even for $p=2$; thus $D$ is connected. To check this, we can assume that $\pi$ is a cycle (cf. Subsubsection 4.1.2). We show that the assumption that $r_D$ modulo $p$ is not injective, leads to a contradiction. This assumption implies that there exists $i\in\{1,\ldots,n+m\}$ such that $w_i\in pB^+(W(k))$. 

By induction on $j\in\{0,\ldots,n+m-1\}$ we check that $w_{i+j}\in pB^+(W(k))$ (here $w_{n+m+s}:=w_s$ for $s\in\{1,\ldots,i-1\}$). The case $j=0$ is obvious and the passage from $j$ to $j+1$ goes as follows. If $i+j\Le n$, then $i_D(a_{i+j+1}\otimes 1)=(1_{H^1(D)}\otimes\Phi_k)(i_D(a_{i+j}\otimes 1))\in (1_{H^1(D)}\otimes\Phi_k)(H^1(D)\otimes_{\dbZ_p} pB^+(W(k)))\subseteq H^1(D)\otimes_{\dbZ_p} pB^+(W(k))$. If $i+j>n$, then due to the fact that $i_D$ respects filtrations we have $i_D(a_{i+j+1}\otimes 1)=(1_{H^1(D)}\otimes\Phi_{1k})(i_D(a_{i+j}\otimes 1))\in (1_{H^1(D)}\otimes\Phi_{1k})(H^1(D)\otimes_{\dbZ_p} pF^1(B^+(W(k)))\subseteq H^1(D)\otimes_{\dbZ_p} pB^+(W(k))$. Thus regardless of what $i+j$ is, we have $w_{i+j+1}\in pB^+(W(k))$. This ends the induction. 

Thus $r_D$ modulo $p$ is the $0$ map. Referring to (2), we get that $\im(j_D)\subseteq pH^1(D)\otimes_{\dbZ_p} \text{gr}^1$. As $\Coker(j_D)$ is annihilated by $\psi$ and as $\psi\in p\dbG_m(V(k))$ for $p=2$, we also have $pH^1(D)\otimes_{\dbZ_p} \text{gr}^1\subseteq \im(j_D)$. Thus $\im(j_D)=pH^1(D)\otimes_{\dbZ_p} \text{gr}^1$ and for each $j\in\{1,\ldots,n\}$, the element $b_j\in H^1(D)\otimes_{\dbZ_p} B^+(W(k))$ defined by the equality $i_D(a_j\otimes 1)=pb_j$ is such that its image in $H^1(D)\otimes_{\dbF_p} k$ is non-zero. Thus $(1_{H^1(D)}\otimes \Phi_k)^s(b_j)$ modulo $p$ is non-zero for all $s\in\dbN^{\ast}$. Taking $s>>0$ we get that the Newton polygon slope ${m\over {m+n}}$ of $(M,\phi)$ is $0$. Thus $D$ is \'etale and connected and therefore trivial. Thus $r_D$ modulo $p$ is injective as its domain is $0$. Contradiction. Thus $r_D$ modulo $p$ is injective even if $p=2$. This implies that $r_D^*$ modulo $p$ is injective even if $p=2$.\endproof

\medskip\noindent
{\bf 4.2.5. Proposition.} {\it Let $(x_1,\ldots,x_{n+m})$ be the reduction modulo $p$ of the solution $(Z_1,\ldots,Z_{n+m})$ we fixed in Subsubsection 4.2.2. Then for each $i\in\{1,\ldots,n+m\}$ there exists $v_i\in\dbG_m(V(k)/pV(k))$ such that we have 
$$q_k(w_i)=q_k(\xi)^{\eps_i}x_iv_i.\leqno (23)$$}
\noindent
{\it Proof:}  To check Formula $(23)$ we can assume that $\pi$ is a cycle $(i_1,\ldots,i_{m+n})$, cf. Subsubsection 4.1.2. Thus $(M,\phi)$ has only one Newton polygon slope ${m\over {m+n}}$. We first check $(23)$ in the case when $n=0$. As $n=0$, we have $D=\mu_{p^{\infty}}^m$ and thus each $Q_{i}(p)$ is ${1\over {p(p-1)}}$ (cf. Lemma 4.2.3 (c)). We also have $\im(i_D)=H^1(D)\otimes_{\dbZ_p} \beta B^+(W(k))$ (see Subsubsection 2.2.3) and thus each $w_1,\ldots,w_m$ is a $\dbG_m(B^+(W(k)))$-multiple of $\beta$. But up to $\dbG_m(V(k)/pV(k))$-multiples, $x_i$ is $p^{1\over {p(p-1)}}$ modulo $p$ (cf. Lemma 4.2.3 (c)) and $q_k(\beta)$ is $p^{1\over {p-1}}$ modulo $p$ (see Subsubsection 2.2.1). As $p^{1\over {p-1}}=p^{1\over p}p^{1\over {p(p-1)}}$, we easily get that Formula $(23)$ holds if $n=0$. 

We prove Formula $(23)$ in the case when $n>0$. As $n>0$ and as $\pi$ is a cycle, we have $Q_i(p)\in [0,{1\over {p(p-1)}})$ for all $i\in\{1,\ldots,n+m\}$ (cf. Lemma 4.2.3 (b) and (c)).

The $r_D$ modulo $p$ and the extension $<,>$ to $B^+(W(k))/pB^+(W(k))$ of the perfect pairing $<,>:T_p(D_{B(k)})/pT_p(D_{B(k)})\times H^1(D)/pH^1(D)\to\dbF_p$ define an $\dbF_p$-linear map $$j_1:T_p(D_{B(k)})/pT_p(D_{B(k)})\to\Hom(\dbD(D[p]),\dbD(B^+(W(k))/pB^+(W(k))))$$
via the formula $(j_1(x))(y)=<x,z>\in B^+(W(k))/pB^+(W(k))$, where $y\in M/pM$ and $x\in T_p(D_{B(k)})/pT_p(D_{B(k)})$ and where $z\in H^1(D)/pH^1(D)\otimes_{\dbF_p} B^+(W(k))/pB^+(W(k))$ is the image of $y$ through $r_D$ modulo $p$. The epimorphism $\dbD(B^+(W(k))/pB^+(W(k)))\twoheadrightarrow\dbD(V(k)/pV(k))$ defined by $q_k$, defines an $\dbF_p$-linear map 
$$j_2:\Hom(\dbD(D[p]),\dbD(B^+(W(k))/pB^+(W(k))))\to\Hom(\dbD(D[p]),\dbD(V(k)/pV(k))).$$ 
As $r_D$ modulo $p$ is injective, $j_1$ is an $\dbF_p$-linear monomorphism. We check that $j_2$ is also an $\dbF_p$-linear monomorphism.

Let $\tilde x:M/pM\to B^+(W(k))/pB^+(W(k))$ be a $k$-linear map that defines an element of $\Ker(j_2)$. The kernel of the $k$-epimorphism $B^+(W(k))/pB^+(W(k))\twoheadrightarrow V(k)/pV(k)$ defined by $q_k$ is annihilated by the Frobenius endomorphism of $B^+(W(k))/pB^+(W(k))$. Thus $\tilde x$ annihilates $\im(\phi)$ modulo $p$. As $n>0$, we get that there exists $i\in\{1,\ldots,n+m\}$ such that $\tilde x$ annihilates $a_i$ modulo $p$. As $\im(\phi)$ modulo $p$ is contained in $\Ker(\tilde x)$ and as $\phi_1$ modulo $p$ maps $\Ker(\tilde x)\cap (F^1/pF^1)$ to $\Ker(\tilde x)$, by induction on $s\in\{0,\ldots,n+m-1\}$ we get that $\tilde x$ annihilates $a_{\pi^s(i)}$ modulo $p$. Thus, as $\pi$ is a cycle and as $\Ker(\tilde x)$ is a $k$-vector space, we have $\tilde x=0$. Therefore $j_2$ is an $\dbF_p$-linear monomorphism. 

Let $x:M/pM\to V(k)/pV(k)$ be a $\sigma^{-1}$-linear map that defines an element of $\break\Hom(\dbD(D[p]),\dbD(V(k)/pV(k)))$ (we recall that $q_k$ is $\sigma^{-1}$-linear). As $x(F^1/pF^1)$ belongs to $q_k(\xi)V(k)/pV(k)$, the image through $x$ of $a_i$ modulo $p$ is of the form $q_k(\xi)^{\eps_i}y_i$, where $y_i\in V(k)/pV(k)$. But $x$ takes $\phi_{\eps_i}(a_i)$ modulo $p$ to $q_k(\xi)^{\eps_{\pi(i)}}y_{\pi(i)}$ as well as to $(-y_i)^p$ (cf. the definition of $\bar\Phi_k$ and $\bar\Phi_{1k}$ in Subsubsection 2.2.1). Thus $(y_1,\ldots,y_{n+m})$ is a  solution of the reduction of $\scrS$ modulo $p$. Conversely, if $(y_1,\ldots,y_{n+m})$ is a solution of the reduction of $\scrS$ modulo $p$, then the $\sigma^{-1}$-linear map $x:M/pM\to V(k)/pV(k)$ that takes $a_i$ modulo $p$ to $q_k(\xi)^{\eps_i}y_i$ defines an element of $\Hom(\dbD(D[p]),\dbD(V(k)/pV(k)))$. Thus, as a set, $\Hom(\dbD(D[p]),\dbD(V(k)/pV(k)))$ is in bijection with the set of solutions of the reduction of $\scrS$ modulo $p$. Thus, as $V(k)$ is strictly henselian and $\scrS$ defines a finite, flat $V(k)$-algebra of degree $p^{n+m}$, the dimension  $\kappa:=\dim_{\dbF_p}(\Hom(\dbD(D[p]),\dbD(V(k)/pV(k))))$ is at most $n+m=\dim_{\dbF_p}(T_p(D_{B(k)})/pT_p(D_{B(k)}))$. As $j_1$ and $j_2$ are $\dbF_p$-linear monomorphisms, by reasons of dimensions we get that $\kappa=n+m$ and that both $j_1$ and $j_2$ are isomorphisms. 

Let $\{a_1^*,\ldots,a_{n+m}^*\}$ be the $W(k)$-basis for $M^*$ which is the dual of the $W(k)$-basis  $\scrB=\{a_1,\ldots,a_m\}$ for $M$. Let 
$$j_3:T_p(D_{B(k)})/pT_p(D_{B(k)})\to M^*\otimes_{W(k)} {}_{\sigma^{-1}}V(k)/pV(k)$$
be the $\dbF_p$-linear map defined naturally by $(1_{M^*}\otimes q_k)\circ i_D^*$ (equivalently, by $j_2\circ j_1$). As $j_2\circ j_1$ is an isomorphism, from the description of $\Hom(\dbD(D[p]),\dbD(V(k)/pV(k)))$ we get that $[V(k)/pV(k)]\im(j_3)$ is the $V(k)/pV(k)$-submodule of $M^*\otimes_{W(k)} V(k)/pV(k)$ generated by elements of the form $\sum_{i=1}^{n+m}  a_i^*\otimes q_k(\xi)^{\eps_i}y_i$,
where $(y_1,\ldots,y_{n+m})$ runs through all reductions modulo $p$ of solutions of $\scrS$ (for $p\Ge 3$ this result is a particular case of [Fa2, Sect. 4, p. 128]). The {\it Moore determinant} of the square matrix of $M_{n+m}(W(k))$ whose rows are $(\eta,\eta^{p^{\gamma_2}},\ldots,\eta^{p^{\gamma_{n+m}}})$, with $\eta$ running through $\mu_{p^{n+m}-1}(W(k))$, is invertible (cf. [Go, Def. 1.3.2 and Lem. 1.3.3]). From this and (21) we get that $[V(k)/pV(k)]\im(j_3)$ is generated by all $a_i^*\otimes q_k(\xi)^{\eps_i}x_i$'s with $i\in\{1,\ldots,m+n\}$. As ${1\over p}+Q_i(p)<{1\over p}+{1\over {p(p-1)}}={1\over {p-1}}\Le 1$, each element $q_k(\xi)^{\eps_i}x_i\in V(k)/pV(k)$ is non-zero; thus $j_3$ is injective (even if $p=2$). But as $i_D^*$ is the transpose of $i_D$ and as we have $t(a_i)=w_ia_i$ (see Subsubsection 4.2.2), we easily get that $[V(k)/pV(k)]\im(j_3)$ is generated by all $a_i^*\otimes q_k(w_i)$'s with $i\in\{1,\ldots,m+n\}$. 

By comparing the two expressions of $[V(k)/pV(k)]\im(j_3)$ we got, we conclude that for each $i\in\{1,\ldots,n+m\}$ the element $q_k(w_i)$ is a $\dbG_m(V(k)/pV(k))$-multiple of $q_k(\xi)^{\eps_i}x_i$. Thus $(23)$ holds even if $n>0$.\endproof

\medskip\noindent
{\bf 4.2.6. A Lubin--Tate quintuple.} To the quintuple $(M,F^1,\phi,T,(t_{\alpha})_{\alpha\in\scrJ})$ we will associate a Lubin--Tate quintuple $(\tilde M,\tilde F^1,\tilde\phi,\tilde T,(\tilde t_{\alpha})_{\alpha\in\tilde\scrJ})$ to which we will apply Subsubsections 4.2.1 to 4.2.5. Let $o(\pi)$ be the order of $\pi$. Let $\eta(\pi)\in\dbN^{\ast}$ be the smallest number with the property that for each cycle $(i_1,\ldots,i_q)$ of $\pi$, there exist at most $\eta(\pi)$ distinct cycles of $\pi$ of the form $(i^\prime_1,\ldots,i^\prime_q)$ and such that up to a cyclic rearrangement we have $\eps_{i^\prime_j}=\eps_{i_j}$ for all $j\in \{1,\ldots,q\}$. Let $r\in o(\pi)\eta(\pi)\dbN^{\ast}$. Let $\{\tilde a_1,\ldots,\tilde a_r\}$ be a $W(k)$-basis for $\tilde M:=W(k)^r$. 

For $s\in\dbN^{\ast}$ and $i\in\{1,\ldots,r\}$ let $\tilde a_{sr+i}:=\tilde a_i$. Let $\tilde D$ be the $p$-divisible group over $\Spec\,W(k)$ whose filtered Dieudonn\'e module is $(\tilde M,\tilde F^1,\tilde\phi)$, where $\tilde F^1:=W(k)\tilde a_1$, $\tilde\phi(\tilde a_1)=p\tilde a_2$, and $\tilde\phi(\tilde a_i)=\tilde a_{i+1}$ for $i\in\{2,\ldots,r\}$ (cf. Lemma 2.2.5 and Proposition 2.2.6). 

Let $(F^i(\scrT(\tilde M)))_{i\in\dbZ}$ be the filtration of $\scrT(\tilde M)$ defined by $\tilde F^1$.  Let 
$$(\tilde M,\tilde F^1,\tilde\phi)^{\otimes s}:=(\tilde M^{\otimes s},(F^i(\scrT(\tilde M))\cap \tilde M^{\otimes s})_{i\in\dbZ},\tilde\phi).$$ Let $(\tilde t_{\alpha})_{\alpha\in\tilde\scrJ}$ be the set of all tensors of $F^0(\scrT(\tilde M))[{1\over p}]$ fixed by $\tilde\phi$. Let $\tilde v_{\alpha}\in\scrT(H^1(\tilde D)[{1\over p}])$ be the element that corresponds to $\tilde t_{\alpha}$ via the following $B^+(W(k))[{1\over {\beta}}]$-linear isomorphism $i_{\tilde D}[{1\over {\beta}}]:\tilde M\otimes_{W(k)} B^+(W(k))[{1\over {\beta}}]\arrowsim H^1(\tilde D)\otimes_{\dbZ_p} B^+(W(k))[{1\over {\beta}}]$, where $i_{\tilde D}$ is the analogue of $i_D$ of Formula (1) but for $\tilde D$ instead of $D$.

Let $\tilde G$ be the schematic closure in $\pmb{GL}_{\tilde M}$ of the subgroup $\tilde G_{B(k)}$ of $\pmb{GL}_{\tilde M[{1\over p}]}$ that fixes $\tilde t_{\alpha}$ for all $\alpha\in\tilde\scrJ$. Let $\tilde T$ be the maximal torus of $\pmb{GL}_{\tilde M}$ that normalizes $W(k)\tilde a_i$ for all $i\in\{1,\ldots,r\}$. Let $\tilde\mu:\dbG_m\to \tilde G$ be the inverse of the canonical split cocharacter  of $(\tilde M,\tilde F^1,\tilde\phi)$, cf. end of Subsection 2.1. As $\tilde G$ contains the conjugates of $\tilde\mu$ under $\tilde\phi$, it also contains $\tilde T$. As $\Lie(\tilde T)$ is generated by elements of $F^0(\scrT(\tilde M))$ fixed by $\tilde\phi$, $\tilde G$ fixes $\Lie(\tilde T)$ and thus $\tilde G$ is contained in $\tilde T$. Thus $\tilde G=\tilde T$ (this also follows from Subsubsection 4.1.3). 

\medskip\noindent
{\bf 4.2.7. The $\tilde w_i$ elements.} The double centralizer of $\tilde T$ in $\pmb{GL}_{\tilde M}$ is $\tilde T$ itself. Thus there exists $\tilde t\in\tilde T(B^+(W(k))[{1\over {\beta}}])$ such that it takes $(i_{\tilde D}[{1\over {\beta}}])^{-1}(H^1(\tilde D)\otimes_{\dbZ_p} B^+(W(k)))$ onto $\tilde M\otimes_{W(k)} B^+(W(k))$, cf. Lemma 4.2.1  applied to $(\tilde D,(\tilde t_{\alpha})_{\alpha\in\tilde\scrJ})$. As in Subsubsection 4.2.2, for $i\in\{1,\ldots,r\}$ let $\tilde w_i\in B^+(W(k))\cap\dbG_m(B^+(W(k))[{1\over {\beta}}])$ be such that $\tilde t(\tilde a_i)=\tilde w_i\tilde a_i$. Let $\tilde w_{i+r}:=\tilde w_i$. Thus $q_k(\tilde w_i)\in V(k)/pV(k)$ is well defined for all $i\in\{1,\ldots,2r\}$. We apply Lemma 4.2.3 and Proposition 4.2.5 to $(\tilde M,\tilde F^1,\tilde\phi)$ (thus $(n,m)$ and $\pi$ get replaced by the pair $(r-1,1)$ and by the cycle $(1\, 2\,\ldots\,r)$). As $q_k(\xi)$ is $p^{1\over p}$ modulo $p$, from $(22)$ and $(23)$ we get that $q_k(\tilde w_1)$ is a $\dbG_m(V(k)/pV(k))$-multiple of $p^{{1\over p}+{p^{-1}\over {p^r-1}}}$ modulo $p$ and that for $i\in\{2,\ldots,r\}$ the element $q_k(\tilde w_i)$ is a $\dbG_m(V(k)/pV(k))$-multiple of $p^{{p^{-2+i}\over {p^r-1}}}$ modulo $p$.

\medskip\noindent
{\bf 4.2.8. Proposition.} {\it The triple $(M,F^1,\phi)$ is a direct summand of $\oplus_{s=1}^r (\tilde M,\tilde F^1,\tilde\phi)^{\otimes s}$.}

\medskip
\proof 
We fix a cycle $\pi_1=(i_1,\ldots,i_q)=(i_{1,1},\ldots,i_{q,1})$ of $\pi$. Let $\eps(\pi_1):=\sum_{s=1}^q \eps_{i_s}$. Let $\pi_1$, $\pi_2,\ldots,\pi_\chi$ be all distinct cycles of $\pi$ of length $q$ that are of the form $\pi_j=(i_{1,j},\ldots,i_{q,j})$ with $(\eps_{i_{1,j}},\ldots,\eps_{i_{q,j}})=(\eps_{i_1},\ldots,\eps_{i_q})$ (here $\chi\in\dbN^{\ast}$ and $j\in\{1,\ldots,\chi\}$). We have ${r\over q}\in\dbN^{\ast}$ and $\chi\Le {r\over q}$, cf. the definitions of $r$, $o(\pi)$, and $\eta(\pi)$ in Subsubsection 4.2.6. For $j\in\{1,\ldots,{r\over q}\}$ let $\varpi_j\in\dbG_m(W(k))$ be such that $\sigma^r(\varpi_j)=\varpi_j$ and the reductions modulo $p$ of $\varpi_1,\ldots,\varpi_{r\over q}$ are linearly independent over $\dbF_{p^q}$. 

For $d\in\{1,\ldots,\eps(\pi_1)\}$ let $s_d\in\{1,\ldots,q\}$ be such that we have $s_1<s_2<\cdots <s_{\eps(\pi_1)}$ and $\{s_d|d\in\{1,\ldots,\eps(\pi_1)\}\}=\{s\in\{1,\ldots,q\}|\eps_{i_s}=1\}$. For $(d,i)\in\{1,\ldots,\eps(\pi_1)\}\times\{0,\ldots,{{r-q}\over q}\}$ let 
$$l_{d+i\eps(\pi_1)}:=r+2-s_d+(r-q)i.$$ 
We note that the numbers $l_1$ modulo $r,\ldots,l_{{r\over q}\eps(\pi_1)}$ modulo $r$ are all distinct. 

For $j\in\{1,\ldots,\chi\}$ let $c_j:(\oplus_{s=1}^q W(k)a_{i_{s,j}},F^1\cap (\oplus_{s=1}^q W(k)a_{i_{s,j}}),\phi)\to (\tilde M,\tilde F^1,\tilde\phi)^{\otimes {r\over q}\eps(\pi_1)}$
be the unique morphism that maps $a_{i_{1,j}}$ to the sum
$$\sum_{v=0}^{{{r-q}\over q}} \sigma^{vq}(\varpi_j)\tilde a_{vq+l_1}\otimes\tilde a_{vq+l_2}\otimes \cdots\otimes\tilde a_{vq+l_{{r\over q}\eps(\pi_1)}}.$$
For $s\in\{1,\ldots,q\}$, we have $\tilde\phi(\tilde a_{s-1+vq+l_1}\otimes \cdots\otimes \tilde a_{s-1+vq+l_{{r\over q}\eps(\pi_1)}})=p^{\chi(s)}\tilde a_{s+vq+l_1}\otimes \cdots\otimes \tilde a_{s+vq+l_{{r\over q}\eps(\pi_1)}}$, where $\chi(s)\in\{0,1\}$. We have $\chi(s)=1$ if and only if there exists $(d,i)\in\{1,\ldots,\eps(\pi_1)\}\times\{0,\ldots,{{r-q}\over q}\}$ such that $s-1+vq+l_{d+i\eps(\pi_1)}=s-1+vq+r+2-s_d+(r-q)i$ is congruent modulo $r$ to $1$ i.e., if and only if $s\in\{s_1,\ldots,s_{\eps(\pi_1)}\}$. This implies that $c_j$ exists; the uniqueness of $c_j$ follows from the fact that $\pi_j$ is a cycle.

For $j\in\{\chi+1,\ldots,{r\over q}\}$, let $\im(c_j)$ be the natural analogue of $\im(c_1),\ldots,\im(c_\chi)$. As the reductions modulo $p$ of $\varpi_j$'s are linearly independent over $\dbF_{p^q}$, the Moore determinant of the square matrix of $M_{{r\over q}}(k)$ whose $j$-th row is the reduction modulo $p$ of $(\varpi_j,\sigma^q(\varpi_j),\ldots,\sigma^{r-q}(\varpi_j))$ is invertible (cf. [Go, Def. 1.3.2 and Lem. 1.3.3]). Thus the morphism $\sum_{j=1}^\chi c_j$ is injective and its image has $\oplus_{s=\chi+1}^{r\over q} \im(c_j)\oplus (S(\pi_1),(F^i(\scrT(\tilde M))\cap S(\pi_1))_{i\in\dbZ},\tilde\phi)$ as a direct supplement in $(\tilde M,\tilde F^1,\tilde\phi)^{\otimes {r\over q}\eps(\pi_1)}$. Here $S(\pi_1)$ is the $W(k)$-span of those elements $\tilde a_{j_1}\otimes \cdots  \otimes\tilde a_{j_{{r\over q}\eps(\pi_1)}}$ of $\tilde M^{\otimes {r\over q}\eps(\pi_1)}$ that are not of the form $\tilde a_{w+l_1}\otimes\cdots\otimes \tilde a_{w+l_{{r\over q}\eps(\pi_1)}}$ for some $w\in\{1,\ldots,r\}$. 

We note that ${r\over q}\eps(\pi_1)\Le r$. We consider a different sequence $\pi_1^\prime,\ldots,\pi_{\chi^\prime}^\prime$ of cycles of $\pi$ that is analogue to the sequence $\pi_1,\ldots,\pi_\chi$. The numbers $q$, $\eps(\pi_1)$, $l_1,\ldots,l_{{r\over q}\eps(\pi_1)}$ are canonically associated to $\pi_1$ and in fact they determine uniquely $\pi$. Thus the morphism $\sum_{j=1}^{\chi^\prime} c_j^\prime$ analogue to $\sum_{j=1}^\chi c_j$ but obtained using the cycles $\pi_1^\prime,\ldots,\pi_{\chi^\prime}^\prime$, is such that its image is contained either in $(S(\pi_1),(F^i(\scrT(\tilde M))\cap S(\pi_1))_{i\in\dbZ},\tilde\phi)$ or in some $(\tilde M,\tilde F^1,\tilde\phi)^{\otimes {r\over q}\eps(\pi_1^\prime)}$ with $\eps(\pi_1^\prime)\neq\eps(\pi_1)$. Thus by adding such morphisms $\sum_{j=1}^\chi c_j$ together we get that $(M,F^1,\phi)$ is a direct summand of $\oplus_{s=1}^r (\tilde M,\tilde F^1,\tilde\phi)^{\otimes s}$.\endproof

\medskip\noindent
{\bf 4.2.9. The projector $\tilde t_{\alpha_0}$.} We fix a direct sum decomposition 
$$\oplus_{s\in\dbN^{\ast}} (\tilde M,\tilde F^1,\tilde\phi)^{\otimes s}=(M,F^1,\phi)\oplus (M^\perp,(F^i(\scrT(\tilde M))\cap M^\perp)_{i\in\dbZ},\tilde\phi)$$ 
such that $M$ is a direct summand of $\oplus_{s=1}^r \tilde M^{\otimes s}$ (cf. Proposition 4.2.8) and thus also of $\scrT(\tilde M)$. Therefore there exists $\alpha_0\in\tilde\scrJ$ such that $\tilde t_{\alpha_0}$ is a projector of $\scrT(\tilde M)$ on $M$. As $M$ is a direct summand of $\oplus_{s=1}^r \tilde M^{\otimes s}$, we get that $M^*$ is a direct summand of $\oplus_{s=1}^r \tilde M^{*\otimes s}$. Thus we can identify $\scrT(M)$ with a direct summand of $\scrT(\tilde M)$ and under this identification each $t_{\alpha}$ is identified with a $\tilde t_{\tilde\alpha}$ for some $\tilde\alpha\in\tilde\scrJ$. Let $L^1(D)$ be the direct summand of $\oplus_{s=1}^r H^1(\tilde D)^{\otimes s}$ that corresponds to $(M,F^1,\phi)$ via $i_{\tilde D}[{1\over {\beta}}]$. The element $\tilde t$ (of Subsubsection 4.2.7) takes $i_{\tilde D}^{-1}(L^1(D)\otimes_{\dbZ_p} B^+(W(k)))$ onto $M\otimes_{W(k)} B^+(W(k))$. We have $L^1(D)=\tilde v_{\alpha_0}(\scrT(H^1(\tilde D)))$ and thus $\Gal(B(k))$ normalizes $L^1(D)$. As we have $\tilde w_i\in B^+(W(k))$ for all $i\in\{1,\ldots,r\}$ and as $M\subseteq\oplus_{s=1}^r\tilde M^{\otimes s}$, the $B^+(W(k))[{1\over {\beta}}]$-linear isomorphism 
$$\tilde i_D[{1\over {\beta}}]:M\otimes_{W(k)} B^+(W(k))[{1\over {\beta}}]\arrowsim L^1(D)\otimes_{\dbZ_p} B^+(W(k))[{1\over {\beta}}]$$ 
induced naturally by $i_{\tilde D}[{1\over {\beta}}]$, induces via restriction a $B^+(W(k))$-linear monomorphism 
$$\tilde i_D:M\otimes_{W(k)} B^+(W(k))\hookrightarrow L^1(D)\otimes_{\dbZ_p} B^+(W(k)).$$
\noindent 
{\bf 4.2.10. Proposition.} {\it The $B^+(W(k))$-linear monomorphism $\tilde i_D:M\otimes_{W(k)} B^+(W(k))\hookrightarrow L^1(D)\otimes_{\dbZ_p} B^+(W(k))$ defined by $i_{\tilde D}$ is $i_D$. In particular, we have $L^1(D)=H^1(D)$.}

\medskip
\proof
 As the crystalline representations of $\Gal(B(k))$ over $\dbQ_p$ are stable under subobjects that are direct summands (see [Fo4, Subsect. 5.5]), we can identify naturally the triples $(i_D[{1\over {\beta}}],H^1(D)[{1\over p}],(v_{\alpha})_{\alpha\in\scrJ})=(\tilde i_D[{1\over {\beta}}],L^1(D)[{1\over p}],(\tilde v_{\tilde\alpha})_{\alpha\in\scrJ})$ and thus we only need to check that $H^1(D)=L^1(D)$. To check this it suffices to consider the case when $\pi$ is a cycle $\pi_1=(i_1,\ldots,i_q)$ (cf. Subsubsection 4.1.2); thus $q=n+m$ and $\eps(\pi_1)=m$. If $w\in\{1,\ldots,r\}$, then the element $\tilde t$ of Subsubsection 4.2.7 takes $\tilde a_{w+l_1}\otimes\tilde a_{w+l_2}\otimes\cdots\otimes \tilde a_{w+l_{{r\over q}\eps(\pi_1)}}\in B^+(W(k))$ to itself times 
$$\tilde w(w):=\prod_{j=1}^{{r\over q}\eps(\pi_1)} \tilde w_{w+l_j}\in B^+(W(k)).$$ 
We have $q_k(\tilde w(w))=\prod_{j=1}^{{r\over q}\eps(\pi_1)} q_k(\tilde w_{w+l_j})$. But $q_k(\tilde w_{w+l_j})$ up to a $\dbG_m(V(k)/pV(k))$-multiple is of the form $p^{{{\eps_{l_{w,j}}}\over p} +{{p^{\gamma_{l_{w,j}}}\over {p^r-1}}}}$ modulo $p$, where $\eps_{l_{w,j}}\in\{0,1\}$ and $\gamma_{l_{w,j}}\in\{-1,\ldots,r-2\}$ (see the end  of Subsubsection 4.2.7). At most one of $\eps_{l_{w,j}}$'s is $1$. The $\gamma_{l_{w,j}}$'s are distinct (as $l_j$ modulo $r$ is uniquely determined by $j\in\{1,\ldots,{r\over q}\eps(\pi_1)\}$) and their number is exactly ${r\over q}\eps(\pi_1)$ and thus it is at most $r$. Thus $q_k(\tilde w(w))$ is the image in $V(k)/pV(k)$ of an element of $V(k)$ whose $p$-adic valuation is at most 
$${1\over p}+\sum_{\gamma=-1}^{r-2} {{p^{\gamma}}\over {p^r-1}}={1\over p}+{1\over {p(p-1)}}={1\over {p-1}}.$$ 
This last $p$-adic valuation is ${1\over {p-1}}$ if and only if ${r\over q}\eps(\pi_1)=r$ i.e., if and only if we have $n+m=q=\eps(\pi_1)=m$. The last equalities are equivalent to $n=0$. In particular we get that if either $n>0$ or $p>2$, then $q_k(\tilde w(w))\neq 0$.

We check that $H^1(D)=L^1(D)$ in the case when $(n,p)=(0,2)$. In this case we can assume $(n,m)=(0,1)$. As $H^1(D)$ has rank $1$, there exists $s\in\dbZ$ such that $H^1(D)=2^sL^1(D)$. As $i_{\tilde D}$ is strict with filtrations, we have $\tilde w_1\in F^1(B^+(W(k)))$. Therefore $\tilde w(w)=\prod_{i=1}^r \tilde w_i\in F^1(B^+(W(k))$ is congruent modulo $F^2(B^+(W(k)))$ to $\xi$ times an element of $B^+(W(k))$ whose image through $q_k$ is $2^{1\over {2(2-1)}}$ modulo $2$, cf. the last part of the previous paragraph and Subsubsection 2.2.1. Thus $\tilde w(w)$ is congruent modulo $F^2(B^+(W(k)))$ to a $\dbG_m(B^+(W(k))/F^2(B^+(W(k))))$-multiple of $\beta$ modulo $F^2(B^+(W(k)))$ (cf. end of Subsubsection 2.2.1). On the other hand we have $\im(i_D)=2^sL^1(D)\otimes_{\dbZ_2} \beta B^+(B(k))$ and therefore $\tilde w(w)$ is a $\dbG_m(B^+(W(k)))$-multiple of $2^s\beta$. As the $V(k)$-module $\text{gr}^1=F^1(B^+(W(k)))/F^2(B^+(W(k)))$ is torsion free, we have $s=0$. Therefore we have $H^1(D)=L^1(D)$ if $n=0$ and $p=2$. 

We check that $H^1(D)=L^1(D)$ in the case when either $n>0$ or $p>2$. We first show that we have an inclusion $H^1(D)\subseteq L^1(D)$ i.e., $L^1_1(D):=H^1(D)\cap L^1(D)$ is $H^1(D)$. If $p=2$ and $n>0$, then $L^1_1(D)$ is a $\dbZ_2$-submodule of $H^1(D)$ such that $i_D$ factors through the inclusion $L^1_1(D)\otimes_{\dbZ_2} B^+(W(k))\hookrightarrow H^1(D)\otimes_{\dbZ_2} B^+(W(k))$ (cf. the fact that $\tilde i_D$ is a $B^+(W(k))$-linear monomorphism); as in the proof of Lemma 2.2.5 applied with $L^1_1(D)$ instead of $H^1(D_1)$, we argue first that $2$ annihilates $H^1(D)/L^1_1(D)$ and second that $H^1(D)=L^1_1(D)\subseteq L^1(D)$ (note that $n>0$ implies that $D^{\text{t}}$ is connected). If $p>2$, then the equality $L_1^1(D)=H^1(D)$ is implied by the fact that $\text{Coker}(j_D)$ is annihilated by $\psi\in V(k)\setminus pV(k)$ (cf. end of Subsubsection 2.2.3). 

Thus we always have $H^1(D)\subseteq L^1(D)$ and therefore in the transpose context we have an inclusion $j_{D,\tilde D}:L^1(D)^*\hookrightarrow T_p(D)=H^1(D)^*$. As either $n>0$ or $p>2$, we have $q_k(\tilde w(w))\neq 0$ for all $w\in\{1,\ldots,r\}$. This implies that the $\dbF_p$-linear map $j_4:L^1(D)^*/pL^1(D)^*\to M^*\otimes_{W(k)} {}_{\sigma^{-1}}V(k)/pV(k)$ that is the natural analogue of the map $j_3$ (of the proof of Proposition 4.2.5) is injective. But $j_4$ is the composite of $j_{D,\tilde D}$ modulo $p$ with $j_3$. Thus $j_{D,\tilde D}$ modulo $p$ is injective and therefore $j_{D,\tilde D}$ is an isomorphism. Thus we have $H^1(D)=L^1(D)$ if either $n>0$ or $p>2$. If we have ${r\over q}\eps(\pi_l)\Le p-2$ for all cycles $\pi_l$ of $\pi$ (with $l\in\scrC(\pi)$), then the identity $H^1(D)=L^1(D)$ is also implied by [Fa2, Thm. 5].\endproof 

\medskip\noindent
{\bf 4.2.11. End of the proof of 1.2 and 4.1.4.} As $\tilde T=\tilde G$ is generated by conjugates of $\tilde\mu$ via powers of $\tilde\phi$ (see Subsubsection 4.2.6), the image of $\tilde T$ in $\pmb{GL}_M$ is the subtorus $T_0$ of $T$. Thus the automorphism of $M\otimes_{W(k)} B^+(W(k))[{1\over {\beta}}]$ defined by the element $\tilde t$ of Subsubsection 4.2.7 is an element $t\in T_0(B^+(W(k))[{1\over {\beta}}])$ that takes the $B^+(W(k))$-submodule $(i_D[{1\over {\beta}}])^{-1}(H^1(D)\otimes_{\dbZ_p} B^+(W(k)))$ of $M\otimes_{W(k)} B^+(W(k))[{1\over {\beta}}]$ onto $M\otimes_{W(k)} B^+(W(k))$; here we used the fact (see Proposition 4.2.10) that we have $H^1(D)=L^1(D)$ regardless of what $\pi$ is. As $T_0$ is a subtorus of $T$, we have $t\in T(B^+(W(k))[{1\over {\beta}}])$ and thus Theorem 4.1.4 holds. This ends the proofs of Theorem 4.1.4 and the Main Theorem.\endproof 

\medskip\smallskip\noindent
{\bf 4.3. Example.} 
Suppose $T=T_0$ has rank $2$ and $\phi^2(\mu)=\mu$. There exist two subcases:

\medskip
{\bf (a)} all Newton polygon slopes of $(M,\phi)$ are $1\over 2$, or

\smallskip
{\bf (b)} the Newton polygon slopes of $(M,\phi)$ are elements of the set $\{0,{1\over 2},1\}$, the multiplicity of the Newton polygon slope ${1\over 2}$ is positive, and either the multiplicity of the Newton polygon slope 0 or of the Newton polygon slope $1$ is as well positive.

\medskip
In the subcase (a) we have $n=m$ and $T_1$ has rank $2$. This is so as we can choose $\scrB$ such that we have $\pi(i)=i+n$ for all $i\in\{1,\ldots,n\}$. As $T_0=T_1$, Lemma 4.2.1  applies.

We consider the subcase (b). Let $M=M_0\oplus M_{1\over 2}\oplus M_1$ be the Newton polygon slope decomposition of $(M,\phi)$; thus for $i\in\{0,{1\over 2},1\}$ all Newton polygon slopes of $(M_i,\phi)$ are $i$. For $s\in\{0,1\}$ let $M_{1\over 2}^s:=M_{1\over 2}\cap F^s$. We have a direct sum decomposition $M_{1\over 2}=M_{1\over 2}^1\oplus M_{1\over 2}^0$. We have no a priori relation between $n$ and $m$. The torus $T_1$ is the center of $\pmb{GL}_{M_0}\times_{\Spec\,W(k)} \pmb{GL}_{M_{1\over 2}^1}\times_{\Spec\,W(k)} \pmb{GL}_{M_{1\over 2}^0}\times_{\Spec\,W(k)} \pmb{GL}_{M_1}$. Thus $T_1$ has rank 4 if both $M_1$ and $M_0$ are non-zero and has rank $3$ otherwise. To fix the ideas we will assume that there exist a pair $(q_0,q_1)\in\{0,\ldots,n\}\times\{0,\ldots,m\}$ such that $\{1,2,\ldots,q_0,n+m+1-q_1,n+m+2-q_1,\ldots,n+m\}$ is the set of elements fixed by $\pi$ and we have $\pi(q_0+s)=n+m+1-q_1-s$ for all $s\in\{1,\ldots,n-q_0\}$. We have $q_0+q_1>0$ and $n-q_0=m-q_1>0$. We can choose $q_k(w_i)\in V(k)/pV(k)$ to be a $\dbG_m(V(k)/pV(k))$-multiple of the reduction modulo $p$ of the following elements of $V(k)$ (cf. Proposition 4.2.5):

\medskip
(i) $1$ for $i\in\{1,\ldots,q_0\}$;

(ii) $(-p)^{1\over {p-1}}$ for $i\in\{n+m+1-q_1,\ldots,n+m\}$;

(iii) $(-p)^{1\over {p(p^2-1)}}$ for $i\in\{q_0+1,\ldots,n\}$;

(iv) $(-p)^{{1\over p}+{1\over {p^2-1}}}$ for $i\in\{n+1,\ldots,n+m-q_1\}$.

\medskip
The value of (ii) is the product of the values of (iii) and (iv). We have $o(\pi)=2$ and $\eta(\pi)=\max\{q_0,n-q_0,q_1\}$. If $q_1>0$, then $T_0$ is the closed subgroup scheme of $\pmb{GL}_M$ that fixes the following three types of tensors fixed by $\phi$: 

\medskip
-- all endomorphisms of $(M,\phi)$,

\smallskip
-- the elements $a_i\in M$ for $i\in\{1,\ldots,q_0\}$, and

\smallskip
-- the elements $a_{q_0+s}\otimes a_{n+m+1-q_1-s}\otimes a_i^* + a_{n+m+1-q_1-s}\otimes a_{q_0+s}\otimes a_i^*\in M^{\otimes 2}\otimes_{W(k)} M^*$, with $s\in\{1,\ldots,n-q_0\}$ and with $i\in\{n+m+1-q_1,\ldots,n+m\}$. 

\medskip
Until Section 5 we will continue to assume that $k=\bar k$ but we come back to the general situation of Subsection 1.1; thus $G$ is not any more assumed to be a torus.  

\medskip\smallskip\noindent
{\bf 4.4. Proof of Corollary 1.4.}
We recall that $\mu:=\mu_{\text{can}}^{-1}$, cf. Subsection 2.1. Let $\sigma_{\phi}$, $M_{\dbZ_p}$, and $G_{\dbZ_p}$ be as in Subsubsection 2.5.1. If $G_{\text{min}}$ is the closed subgroup of $G$ introduced in Subsubsection 2.7.2 and if the family of tensors $(t_{\alpha})_{\alpha\in\scrJ_{\text{min}}}$ of $\scrT(M_{\dbZ_p}[{1\over p}])\subseteq \scrT(M[{1\over p}])$ is as in Lemma 2.7.3, then to prove the Corollary 1.4 we can assume that $\scrJ=\scrJ_{\text{min}}$ is a subset of $\scrJ_{\text{twist}}$, that $G=G_{\text{min}}$, and that $(t_{\alpha})_{\alpha\in\scrJ}$ is  in fact the family of all tensors of $\{x\in F^0(\scrT(M))[{1\over p}]|\phi(x)=x\}$.

Let $W\in\scrW$, cf. the notations of Corollary 1.4. From the functorial aspects of the canonical split cocharacters we get that $W$ is normalized by the image of $\mu_{\text{can}}=\mu^{-1}$. From this and the fact that $\phi(W[{1\over p}])=W[{1\over p}]$, we get that $W[{1\over p}]$ is normalized by all conjugates of $\mu_{B(k)}$ under powers of $\phi$. This implies that each $W$ is normalized by $G=G_{\text{min}}$. 

Let $\rho:(M,(t_{\alpha})_{\alpha\in\scrJ})\arrowsim (H^1(D)\otimes_{\dbZ_p} W(k),(v_{\alpha})_{\alpha\in\scrJ})$ be an isomorphism as in the Main Theorem. We denote also by $\rho$ the isomorphism $\rho:\scrT(M)\arrowsim \scrT(H^1(D)\otimes_{\dbZ_p} W(k))$ induced by it. 

We check that for each $W\in\scrW$ we have $\rho(W)=W^{\acute et}\otimes_{\dbZ_p} W(k)$. As $i_D[{1\over {\beta}}]:M\otimes_{W(k)} B^+(W(k))[{1\over {\beta}}]\arrowsim H^1(D)\otimes_{\dbZ_p} B^+(W(k))[{1\over {\beta}}]$ maps $W\otimes_{W(k)} B^+(W(k))[{1\over {\beta}}]$ to $W^{\acute et}\otimes_{\dbZ_p} B^+(W(k))[{1\over {\beta}}]$, to check that we have $\rho(W)=W^{\acute et}\otimes_{\dbZ_p} W(k)$ it suffices to show that the $B^+(W(k))[{1\over {\beta}}]$-linear automorphism $h:=i_D[{1\over {\beta}}]^{-1}\circ\rho_{B^+(W(k))[{1\over {\beta}}]}$ of $M\otimes_{W(k)} B^+(W(k))[{1\over {\beta}}]$ normalizes $W\otimes_{W(k)} B^+(W(k))[{1\over {\beta}}]$. As $h$ fixes each element $t_{\alpha}$ with $\alpha\in\scrJ$, we have $h\in G(B^+(W(k))[{1\over {\beta}}])$. As $W$ is normalized by $G$, we get that $h$ normalizes $W\otimes_{W(k)} B^+(W(k))[{1\over {\beta}}]$. Thus we have $\rho(W)=W^{\acute et}\otimes_{\dbZ_p} W(k)$ for all $W\in\scrW$.

To end the proof of Corollary 1.4 we are left to show that we can choose $\rho$ such that it maps $t_{\alpha}$ to $v_{\alpha}$ for all $\alpha\in\scrJ_{\text{twist}}$. To check this let $\tilde D:=D\oplus \mu_{p^{\infty}}$; if $p=2$, then the property (EC) holds for $\tilde D$ as well. The Dieudonn\'e module of $\tilde D$ is $(\tilde M,\tilde\phi):=(M,\phi)\oplus (W(k),p\sigma)$.  The Hodge filtration of $\tilde M$ defined by $\tilde D$ is $\tilde F^1:=F^1\oplus W(k)$. Let $(F^i(\scrT(\tilde M)))_{i\in\dbZ}$ be the filtration of $\scrT(\tilde M)$ defined by $\tilde F^1$. Let $\tilde\mu:\dbG_m\to\pmb{GL}_{\tilde M}$ be the inverse of the canonical split cocharacter of $(\tilde M,\tilde F^1,\tilde\phi)$; it acts on $M$ as $\mu$ does (cf. the functorial properties of canonical split cocharacters). Let $(\tilde t_{\alpha})_{\alpha\in\tilde\scrJ}$ be the family of all tensors of $F^0(\scrT(\tilde M))[{1\over p}]$ fixed by $\phi$. Let $(\tilde v_{\alpha})_{\alpha\in\tilde\scrJ}$ be the family of tensors of $\scrT(H^1(\tilde D))[{1\over p}]$ that corresponds to $(\tilde t_{\alpha})_{\alpha\in\tilde\scrJ}$ via Fontaine comparison theory for $\tilde D$. Let $\tilde\rho:(\tilde M,(\tilde t_{\alpha})_{\alpha\in\tilde\scrJ})\arrowsim (H^1(\tilde D),(\tilde v_{\alpha})_{\alpha\in\tilde\scrJ})$, cf. Main Theorem applied to $\tilde D$. By composing $\tilde \rho$ with an automorphism 
of $(\tilde M,(\tilde t_{\alpha})_{\alpha\in\tilde\scrJ})$ defined by an element of the image of $\tilde\mu(W(k)):\dbG_m(W(k))\to \pmb{GL}_{\tilde M}(W(k))$ we can assume that $\tilde\rho$ takes each element $t$ of the set $\cup_{i\in\dbZ} \{x\in F^i(\scrT(W(k))|(p\sigma)(x)=p^ix\}$ to the element $v$ of $\scrT(H^1(\mu_{p^{\infty}}))$ that corresponds to $t$ via Fontaine comparison theory. Let $\rho:(M,(t_{\alpha})_{\alpha\in\scrJ})\arrowsim (H^1(D)\otimes_{\dbZ_p} W(k),(v_{\alpha})_{\alpha\in\scrJ})$ be defined by the restriction of $\tilde \rho$ to $M$. For $\alpha\in\scrJ_{\text{twist}}$, let $i\in\dbZ$ be such that $t_{\alpha}\in \{x\in F^i(\scrT(M))[{1\over p}]|\phi(x)=p^ix\}$. Let $t_{-i}$ be a generator of the $\dbZ_p$-module $\{x\in F^{-i}(\scrT(W(k))|(p\sigma)(x)=p^{-i}x\}$.

We have $t_{\alpha}\otimes t_{-i}\in \{x\in F^0(\scrT(\tilde M))[{1\over p}]|\tilde\phi(x)=x\}$. Thus $\tilde \rho$ maps $t_{\alpha}\otimes t_{-i}$ to $v_{\alpha}\otimes v_{-i}$, where $v_{-i}$ corresponds to $t_{-i}$ via Fontaine comparison theory for $\mu_{p^{\infty}}$. As $\tilde\rho$ also maps $t_{-i}$ to $v_{-i}$, we conclude that $\tilde \rho$ and thus also $\rho$ maps $t_{\alpha}$ to $v_{\alpha}$. Therefore we can take $\rho_{\text{twist}}:(M,(t_{\alpha})_{\alpha\in\scrJ_{\text{twist}}})\arrowsim (H^1(D)\otimes_{\dbZ_p} W(k),(v_{\alpha})_{\alpha\in\scrJ_{\text{twist}}})$ to be defined by $\rho$.\endproof

\medskip\noindent
{\bf 4.4.1. Example.} Suppose that the hypotheses of the Main Theorem hold and we have an isogeny $\lambda_D:D\to D^{\text{t}}$. Let $\lambda_M:M\times M\to W(k)$ and $\lambda_{H^1(D)}:H^1(D)\times H^1(D)\to \dbZ_p$  be bilinear forms on $M$ and $H^1(D)$ defined naturally by $\lambda_D$. We can naturally view $\lambda_M$ as a tensor of $\{x\in F^{-1}(M^{* \otimes 2})|\phi(x)=p^{-1}x\}$. From Corollary 1.4 we get that there exists an isomorphism $\rho:(M,(t_{\alpha})_{\alpha\in\scrJ})\arrowsim (H^1(D),(v_{\alpha})_{\alpha\in\scrJ})$ such that for all $x$, $y\in M$ we have $\lambda_M(x,y)=\lambda_{H^1(D)}(\rho(x),\rho(y))$.

\medskip
We have the following corollary of the Main Theorem and its proof. 

\medskip\smallskip\noindent
{\bf 4.5. Corollary.}  {\it Suppose that the hypotheses of the Main Theorem hold. Then there exists an isomorphism $\rho_{\dbZ_p}:(M_{\dbZ_p},(t_{\alpha})_{\alpha\in\scrJ})\arrowsim (H^1(D),(v_{\alpha})_{\alpha\in\scrJ})$ and therefore the group schemes $G_{\dbZ_p}$ and $G_{\dbZ_p}^{\acute{et}}$ over $\Spec\,\dbZ_p$ are isomorphic.}

\medskip
\proof
We recall from [BLR, Ch. 7, Sect. 7.1, Thm. 5] that there exists a unique smooth group scheme $G^\prime$ over $\Spec\,W(k)$ which is equipped with a homomorphism $i_G:G^\prime\to G$ for which the following universal property holds: if $Y$ is a smooth $\Spec\,W(k)$-scheme, then each morphism $Y\to G$ of $\Spec\,W(k)$-schemes factors uniquely through $i_G$. By considering smooth $B(k)$-schemes $Y$, one gets that the generic fibre of $i_G$ is an isomorphism. See [BLR, Ch. 3, Sect. 3.2] for dilatations. As $i_G$ is obtained using a sequence of dilatations centered on special fibres (see [BLR, Ch. 7, Sect. 7.1,  pp. 174--175]), $G^\prime$ is affine. Let $G^{\prime 0}_k$ be the identity component of $G^\prime_k$.  Let $G^{\prime 0}$ be the open subgroup scheme of $G^\prime$ whose generic and special fibres are $G_{B(k)}$ and $G^{\prime 0}_k$ (respectively). As $G^{\prime 0}$ is the complement in $G^{\prime}$ of a divisor of $G^{\prime}$, it is an affine, smooth  group scheme over $\Spec\,W(k)$. As the monomorphism $\dbZ_p\hookrightarrow  W(k)$ is of index of ramification $1$  in the sense of [BLR, Ch. 3, Sect. 3.6, Def. 1], the analogue $G^{\prime 0}_{\dbZ_p}\to G_{\dbZ_p}$ of the homomorphism $G^{\prime 0}\to G$ is a $\dbZ_p$-structure for $G^{\prime 0}\to G$ (cf. [BLR, Ch. 3, Sect. 3.6, Cor. 6 and Ch. 7, Sect. 7.1, Thm. 5]).

 Let $U$, $g\in U(W(k))$, and $T$ be as in Subsections 2.6 and 2.7. Due to the universal property of $i_G$ and the fact that $U_k$ is connected (being a product of $\dbG_a$'s), one gets that the unipotent group scheme $U$ is naturally a closed subgroup scheme of $G^{\prime 0}$.

The $\dbZ_p$-structure $(M_{\dbZ_p},(t_{\alpha})_{\alpha\in\scrJ})$ of $(M,(t_{\alpha})_{\alpha\in\scrJ})$ defined by $\sigma_{\phi}=\phi\mu(p):M\arrowsim M$ is isomorphic to the $\dbZ_p$-structure of $(M,(t_{\alpha})_{\alpha\in\scrJ})$ defined by $g\sigma_{\phi}:M\arrowsim M$. This is so as a standard application of Lang theorem to $G^{\prime 0}_{\dbZ_p}$ and $g\in G^{\prime 0}(W(k))\leqslant G(W(k))$ shows that there exists an element $\tilde g\in G^{\prime 0}(W(k))\leqslant G(W(k))$ such that we have $g\sigma_{\phi}=\tilde g\sigma_{\phi} \tilde g^{-1}$ (to be compared with [NV, Prop. 2.1]). Thus, as in the proof of Theorem 3.5.1 we argue that to prove the existence of $\rho_{\dbZ_p}$ we can assume that $D$ is connected if $p=2$ and we can replace $(M,\phi,G,(t_{\alpha})_{\alpha\in\scrJ})$ by $(M,g\phi,T,(t_{\alpha})_{\alpha\in\scrJ_T})$, with $(t_{\alpha})_{\alpha\in\scrJ_T}$ as in Subsection 2.7. In other words, we can assume that $G=T$ is a torus and in this case the existence of the isomorphism $\rho_{\dbZ_p}$ follows from Lemma 2.5.2 (a) and the Main Theorem.\endproof

\bigskip\smallskip
\noindent
{\boldsectionfont 5. The ramified context}
\bigskip\smallskip

Let $k$ be again an arbitrary perfect field of characteristic $p>0$. Let $K$ be a finite, totally ramified field extension of $B(k)$. Let $V$ be the ring of integers of $K$ i.e., the normalization of $W(k)$ in $K$. Let $e:=[V:W(k)]$. Let $\pi_V$ be a uniformizer of $V$. Let $X$ be a free variable. The minimal polynomial $f_e\in W(k)[X]$ of $\pi_V$ over $W(k)$ is an Eisenstein polynomial. Let $\grS:=W(k)[[X]]$. Let $S_e$ be the $\grS$-subalgebra of $B(k)[[X]]$ generated by all ${{X^{en}}\over {n!}}$ with $n\in\dbN$; it is the divided power hull of any one of the ideals $(X^e)$, $(f_e)$, or $(p,X^e)=(p,f_e)$ of $\grS$. Let $J_e$ (resp. $K_e$) be the ideal of $S_e$ generated by all ${{f_e^{n}}\over {n!}}$ (resp. by all ${{X^{en}}\over {n!}}$) with $n\in\dbN^{\ast}$. Let $R_e:=S_e^\wedge$. By mapping $X$ to $\pi_V$, we get identifications $V=W(k)[X]/(f_e)=S_e/J_e$ and a $W(k)$-epimorphism $e_V:R_e\twoheadrightarrow V$. Let $F^1(R_e):=\Ker(e_V)$. Let $\Rtil_e$ be the completion of $S_e$ with respect to the decreasing filtration given by its ideals $K_e^{[n]}$, $n\in\dbN$. Thus $\Rtil_e=\text{proj}.\text{lim}._{n\in\dbN^{\ast}} S_e/K_e^{[n]}$. We recall that $K_e^{[0]}:=S_e$ and that for $n\Ge 1$ the ideal $K_e^{[n]}$ of $S_e$ is generated by all products $\frac{\diamondsuit_1^{a_1}}{a_1!}\cdots\frac{\diamondsuit_m^{a_m}}{a_m!}$, with $\diamondsuit_1,\ldots,\diamondsuit_m\in K_e$ and $m$, $a_1,\ldots,a_m\in\dbN$ such that $a_1+\cdots+a_m\ge n$.

In this section we study ramified analogues of the Main Theorem over $R_e$. Subsection 5.1 lists properties of $R_e$ and $\tilde R_e$. Subsection 5.2 presents the basic setting on $p$-divisible groups. Theorem 5.3 refines the deformation theory of [Fa2, Ch. 7]. See Subsections 5.3 to 5.5 for the main results of this section. In particular, Subsection 5.5 presents a crystalline variant (converse) to [Ki3, Cor. (1.4.3)] for $p>2$. The counterexamples of Subsection 5.6 emphasize that for $p>2$ the hypotheses of Theorem 5.4 are needed in general.  

\medskip\smallskip\noindent
{\bf 5.1. On $R_e$ and $\tilde R_e$.} If $m\in\dbN^{\ast}$, $m\Ge n$, then we have ${{X^{em}}\over {m!}}\in K_e^{[n]}$. If $p>2$, then $\Rtil_e$ is also the completion of $S_e$ with respect to its decreasing filtration $(J_e^{[n]})_{n\in\dbN}$; thus for $p>2$ we have as well a $W(k)$-epimorphism $\tilde e_V:\Rtil_e\twoheadrightarrow V$ that maps $X$ to $\pi_V$. For $p=2$, we only have a $W(k)$-epimorphism $\tilde e_{V/2V}:\Rtil_e\twoheadrightarrow V/2V$ that maps $X$ to $\pi_V$ modulo $2$.

The $W(k)$-algebra $\Rtil_e$ (resp. $R_e$) consists
of formal power series $\sum_{n\Ge 0} a_nX^n$ such that the sequence $([{n\over e}]!a_n)_{n\in\dbN}$ is formed by elements of $W(k)$ (resp. is formed by elements of $W(k)$ and converges to $0$). We have $W(k)$-epimorphisms from $S_e$, $R_e$, and $\Rtil_e$ onto $W(k)$ defined by the rule $\sum_{n\ge 0} a_nX^n\to a_0$. 

With the notations of Subsection 2.2, we have $s_k(f_e(\underline{\pi_V},0,0,\ldots))=f_e(\pi_V)=0$ and one can easily check that $f_e(\underline{\pi_V},0,0,\ldots)$ generates $\Ker(s_k)$ as well. There exists a $W(k)$-monomorphism 
$i_e:R_e\hookrightarrow B^+(W(k))$ 
defined by the rule: $X\mapsto (\underline{\pi_V},0,0,\ldots)\in W(A_k)$. Let $\Phi_k$ be the Frobenius lift of $\grS$, $S_e$, $R_e$, or $\Rtil_e$ that is compatible with $\sigma$ and such that
$\Phi_k(X)=X^p$ (it makes sense to denote it also by $\Phi_k$, as each $i_e:R_e\hookrightarrow B^+(W(k))$ respects Frobenius lifts). The sequence $({{[{{np}\over e}]!}\over {[{n\over e}]!}})_{n\in\dbN^{\ast}}$ of integers in $W(k)$ converges to $0$ in the $p$-adic topology and thus we have $\Phi_k(\Rtil_e)\subseteq R_e$. 

As $R_edx$ is the $p$-adic completion of the quotient of $\Omega_{R_e/W(k)}$ by its $p$-torsion submodule (equivalently, by its $R_e$-submodule generated by all relations of the form $d\delta_q(y)=\delta_{q-1}(y)dy$ with $y\in F^1(R_e)$, $q\in\dbN^{\ast}$, and $\delta_q(y)={{y^q}\over {q!}}$ defining the natural divided power structure on $F^1(R_e)$), below we will use $R_edx$ in relation to connections. 

For $q\in\dbN^{\ast}$ we consider the ideals $\Itil_e(q):=\{\sum_{n=0}^{\infty} a_nX^n\in \Rtil_e|a_0=\cdots=a_{q-1}=0\}$ and $I_e(q):=\{\sum_{n=0}^{\infty} a_nX^n\in R_e|a_0=\cdots=a_{q-1}=0\}$ of $\Rtil_e$ and $R_e$ (respectively). We have $I_e(q)=R_e\cap \tilde I_e(q)$.

\medskip\smallskip\noindent
{\bf 5.2. The setting.} Let $H$ be a $p$-divisible group over $\Spec\,V$. Let $(N,\phi_N,\nabla_N)$ be the evaluation of $\dbD(H_{V/pV})$ at the thickening associated naturally to the closed embedding $\Spec\,V/pV\hookrightarrow \Spec\,R_e$. Thus $N$ is a free $R_e$-module of rank equal to the height of $H$, $\phi_N:N\to N$ is a $\Phi_k$-linear endomorphism, and $\nabla_N:N\to N\otimes_{R_e} R_edX$ is a connection on $N$ with respect to which $\phi_N$ is horizontal i.e., we have $\nabla_N\circ\phi_N=(\phi_N\otimes d\Phi_k)\circ\nabla_N$. The connection $\nabla_N$ is integrable and nilpotent modulo $p$. 

\medskip\noindent
{\bf 5.2.1. Uniqueness of connections.} If $\tilde\nabla_N$ is a connection on $N$ such that we have $\tilde\nabla_N\circ\phi_N=(\phi_N\otimes d\Phi_k)\circ\tilde\nabla_N$, then as $\Phi_k(X)=X^p$ by induction on $q\in\dbN^{\ast}$ we get $\nabla_N-\tilde\nabla_N\in X^{q}\End_{R_e}(N)\otimes_{R_e} dXR_e[{1\over p}]$. Thus $\nabla_N=\tilde\nabla_N$. A similar argument shows that $\nabla_N$ modulo $I_e(q)$ is uniquely determined by $\phi_N$ modulo $I_e(q)$ and that the extension of $\nabla_N$ to a connection on $N\otimes_{R_e} \Rtil_e$ (to be denoted also by $\nabla_N$) is the unique connection such that the Frobenius endomorphism $\phi_N\otimes\Phi_k$ of $N\otimes_{R_e} \Rtil_e$ is horizontal.

\medskip\noindent
{\bf 5.2.2. Tensors.} Let $(F^i(\scrT(N/F^1(R_e)N)))_{i\in\dbZ}$ be the filtration of $\scrT(N/F^1(R_e)N)$ defined by the direct summand $F^1_H$ of $N/F^1(R_e)N$ that is the Hodge filtration of $H$. We consider a family $(t_{H\alpha})_{\alpha\in\scrJ}$ of tensors of $\scrT(N[{1\over p}])$ fixed by $\phi_N$ and whose images in $\scrT(N/F^1(R_e)N)[{1\over p}]$ belong to $F^0(\scrT(N/F^1(R_e)N))[{1\over p}]$. Let $H^1(H):=T_p(H^{\text{t}}_{K})(-1)$. Let $v_{H\alpha}\in\scrT(H^1(H)[{1\over p}])$ correspond to $t_{H\alpha}$ via the $B^+(W(k))$-linear monomorphism $i_H:N\otimes_{R_e} {}_{i_e}B^+(W(k))\hookrightarrow H^1(H)\otimes_{\dbZ_p} B^+(W(k))$ obtained as $i_D$ of (1) was (this time $i_H$ is $\Gal(K)$-invariant; see [Fa2, Sect. 4, p. 127] for the canonical action of $\Gal(K)$ on $N\otimes_{R_e} {}_{i_e}B^+(W(k))$). 

Let $(M,\phi,(t_{\alpha})_{\alpha\in\scrJ}):=(N,\phi_N,(t_{H\alpha})_{\alpha\in\scrJ})\otimes_{R_e} R_e/I_e(1)$. Let $G_{B(k)}$ be the subgroup of $\pmb{GL}_{M[{1\over p}]}$ that fixes $t_{\alpha}$ for all $\alpha\in\scrJ$. Let $G$ be the schematic closure of $G_{B(k)}$ in $\pmb{GL}_M$. 

It is well known that there exist isomorphisms 
$$\scrK:(M\otimes_{W(k)} R_e[{1\over p}],\phi\otimes\Phi_k)\arrowsim (N\otimes_{R_e} R_e[{1\over p}],\phi_N\otimes\Phi_k)$$ 
(to be compared with [Fa2, Sect. 6, p. 132] which works with $\Rtil_e$ instead of with $R_e$). We choose such an isomorphism $\scrK$ that lifts $1_{M[{1\over p}]}$. There exists no element of $\scrT(M)[{1\over p}]\otimes I_e(1)[{1\over p}]$ fixed by $\phi\otimes\Phi_k$. This implies that (i) $\scrK$ is the unique such isomorphism that lifts $1_{M[{1\over p}]}$ and that (ii) each $t_{H\alpha}$ is the extension of $t_{\alpha}$ via the $R_e[{1\over p}]$-linear isomorphism $\scrT(M)[{1\over p}]\otimes_{B(k)} R_e[{1\over p}]\arrowsim\scrT(N)[{1\over p}]$ induced by $\scrK$. 

Thus the subgroup scheme $G_{R_e[{1\over p}]}$ of $\pmb{GL}_{N[{1\over p}]}$ that fixes $t_{H\alpha}$ for all $\alpha\in\scrJ$ is isomorphic to $G_{B(k)}\times_{\Spec\,B(k)} \Spec\,R_e[{1\over p}]$ and therefore it is smooth. Let $\scrG$ be the schematic closure of $G_{R_e[{1\over p}]}$ in $\pmb{GL}_N$; $\scrG$ is not a priori a closed subgroup scheme of $\pmb{GL}_{N}$.
Let $\scrG^{\acute et}$ be the schematic closure in $\pmb{GL}_{H^1(H)}$ of the subgroup of $\pmb{GL}_{H^1(H)[{1\over p}]}$ that fixes $v_{H\alpha}$ for all $\alpha\in\scrJ$.

\medskip\noindent
{\bf 5.2.3. Three possible ramified analogues of the Main Theorem.} If $k=\bar k$, then one would like to know if any one of the following three conditions holds:

\medskip
{\bf (i)} there exists an isomorphism $(N,(t_{H\alpha})_{\alpha\in\scrJ})\arrowsim (H^1(H)\otimes_{\dbZ_p} R_e,(v_{H \alpha})_{\alpha\in\scrJ})$;

\smallskip
{\bf (ii)} there exists an isomorphism $(M,(t_{\alpha})_{\alpha\in\scrJ})\arrowsim (H^1(H)\otimes_{\dbZ_p} W(k),(v_{H \alpha})_{\alpha\in\scrJ})$;

\smallskip
{\bf (iii)} there exists an isomorphism $(N/F^1(R_e)N,(t_{H\alpha})_{\alpha\in\scrJ})\arrowsim (H^1(H)\otimes_{\dbZ_p} V,(v_{H \alpha})_{\alpha\in\scrJ})$.

\medskip\noindent
Obviously (i) implies both (ii) and (iii). The goals of this section are to prove that:

\medskip
$\bullet$ for $p>2$, the condition (i) holds provided the pair $(H,(t_{H\alpha})_{\alpha\in\scrJ})$ has a good deformation to $W(k)$ (i.e., it has a lift to $R_e$ in the below sense) (cf. Theorem 5.4 below);

\smallskip
$\bullet$ for $p>2$, $\scrG$ is a reductive group scheme if and only if $\scrG^{\acute et}$ is so (cf. Corollary 5.5.1);

\smallskip
$\bullet$ weaker versions of the above two results hold for $p=2$;

\smallskip
$\bullet$ in general, the condition (i) does not hold (cf. Subsection 5.6).

\medskip\noindent
{\bf 5.2.4.  Definitions.} {\bf (a)} By a {\it filtration lift} of the pair $(H,(t_{H\alpha})_{\alpha\in\scrJ})$ to $R_e$ we mean a direct summand $F^1_N$ of $N$ that lifts $F^1_H$ and such that we have $t_{H\alpha}\in F^0(\scrT(N))[{1\over p}]$ for all $\alpha\in\scrJ$, where $(F^i(\scrT(N)))_{i\in\dbZ}$ is the filtration of $\scrT(N)$ defined by $F^1_N$. 

\smallskip
{\bf (b)} By a {\it lift} of the pair $(H,(t_{H\alpha})_{\alpha\in\scrJ})$ to $R_e$ we mean a pair $(H_{R_e},(t_{H\alpha})_{\alpha\in\scrJ})$, where $H_{R_e}$ is a $p$-divisible group over $\Spec\,R_e$ whose reduction modulo $F^1(R_e)$ is $H$ and whose Hodge filtration is a direct summand of $N$ that is a filtration lift of the pair $(H,(t_{H\alpha})_{\alpha\in\scrJ})$.

\medskip\noindent
{\bf 5.2.5. Remarks.} {\bf (a)} If $e\le p-1$, then there exists an isomorphism $(M\otimes_{W(k)} V,(t_{\alpha})_{\alpha\in\scrJ})\arrowsim (N/F^1(R_e)N,(t_{H\alpha})_{\alpha\in\scrJ})$ induced by $\scrK$ modulo $F^1(R_e)[{1\over p}]$ (cf. [La, Thm. 2.1]). Thus condition 5.2.3 (ii) implies condition 5.2.3 (iii) if $e\le p-1$. If $e\le p-1$ and $G$ is smooth, then as in the proof of Lemma 2.5.2 (b) one argues that conditions 5.2.3 (ii) and (iii) are in fact equivalent. 

\smallskip
{\bf (b)} If $p>2$, then the Grothendieck--Messing deformation theory (or the Theorem 5.3 below) implies that we have a canonical identification between the filtration lifts of the pair $(H,(t_{H\alpha})_{\alpha\in\scrJ})$ to $R_e$ and lifts of the pair $(H,(t_{H\alpha})_{\alpha\in\scrJ})$ to $R_e$.

\smallskip
{\bf (c)} Suppose $k=\bar k$. At a first glance one would be inclined to think that at least for $p>2$, there exists a lift of the pair $(H,(t_{H\alpha})_{\alpha\in\scrJ})$ to $R_e$ if and only if the condition 5.2.3 (i) holds. In what follows we will check the ``only if'' part (cf. Theorem 5.4 below). We hope that future work will clarify if or when the ``if'' part holds.  

\medskip\noindent
{\bf 5.2.6. Lemma.} {\it We assume that $\scrG$ is a split reductive, closed subgroup scheme of $\pmb{GL}_N$. Then there exists a filtration lift of the pair $(H,(t_{H\alpha})_{\alpha\in\scrJ})$ to $R_e$.}

\medskip
\proof
This is only a variant of either [Va1, Subsubsects. 5.3.1 and 5.3.2] or [Ki3, Prop. (1.1.5) and Lem. (1.4.5)].\endproof

\medskip\noindent
{\bf 5.2.7. Extra notations.} Until Subsection 5.5 we will assume that there exists a filtration lift $F^1_N$ of the pair $(H,(t_{H\alpha})_{\alpha\in\scrJ})$ to $R_e$. Let $F^1:=F^1_N\otimes_{R_e} R_e/I_e(1)$. Let  $M=F^1\oplus F^0$ and $\mu:\dbG_m\to G$  be as in Subsection 2.1; these notations make sense even if $p=2$ and $(M,F^1,\phi)$ is not the filtered Dieudonn\'e module of a $p$-divisible group over $\Spec\,W(k)$. 

Let $m$ be the relative dimension of the unipotent group scheme $U$ of Subsection 2.6. Let $R:=W(k)[[z_1,\ldots, z_m]]$ be equipped with the Frobenius lift $\Phi_R$ which is compatible with $\sigma$ and which maps each $z_j$ to $z_j^p$. We identify $Y=Y^\wedge=\Spec\,R$ with the spectrum of the completion of the local ring of $U$ at the identity element of $U_k$ and thus we can speak about the universal element $\tau:=\tau_{\text{univ}}\in U(R)$. Let $\Spec\,Q$ be the unique connected component of $\Spec\,Q_{\infty}^\wedge$ of Theorem 3.2 (c) with $\tilde G=G$ which has the property that the morphism $\Spec\,Q\to\Spec\,R$ is a pro-\'etale cover; we can identify canonically $R=Q$ (cf. Theorem 3.3.1). Let $M_R:=M\otimes_{W(k)} R$ and $F^1_R:=F^1\otimes_{W(k)} R$ and let $\phi_{0,R}$ be as in Subsubsection 3.1.1 (a). Let the connection $\nabla$ on $M_R=M_Q$ be as in the beginning of  Subsection 3.4. Let $a_{\infty}:\Spec\,W(k)\hookrightarrow \Spec\,Q=\Spec\,R$ be the section defined by the ideal $I:=(z_1,\ldots,z_m)$ of $R$. Let $\scrM_R^0:=(M_R,F^1_R,\phi_{0,R})$. 

\medskip\noindent
{\bf 5.2.8. Extensions.} Let $\nu:R\to\Rtil_e$ be a $W(k)$-homomorphism that maps $I$ to $\Itil_e(1)$. We denote also by $\nu$ the resulting morphism $\Spec\,\Rtil_e\to\Spec\,R$.  By the extension of $(\scrM_R^0,\nabla,(t_{\alpha})_{\alpha\in\scrJ})$ through $\nu$ we mean the quintuple 
$$(M\otimes_{W(k)} \Rtil_e,F^1\otimes_{W(k)} \Rtil_e,\phi_{\nu},\nabla_{\nu},(t_{\alpha})_{\alpha\in\scrJ}),$$ 
where $(M\otimes_{W(k)} \Rtil_e,\phi_{\nu},\nabla_{\nu})$ is the evaluation at the thickening $(\Spec\,\tilde R_e/p\tilde R_e\hookrightarrow\Spec\,\tilde R_e,\delta(p))$ of the pull-back to $CRIS(\Spec\,\Rtil_e/p\Rtil_e/\Spec\,W(k))$ of the (uniquely determined) $F$-crystal on $CRIS(\Spec\,R/pR/\Spec\,W(k))$ whose evaluation at the thickening $(\Spec\,R/pR\hookrightarrow \Spec\,R),\delta(p))$ is $(M_R,\phi_{0,R},\nabla)$. Thus the $\Phi_k$-linear endomorphism $\phi_{\nu}$ of $M\otimes_{W(k)} \Rtil_e$ is defined using a correction automorphism $\grK$ as in the proof of Theorem 3.4.1 and $\nabla_{\nu}$ is the unique connection on $M\otimes_{W(k)} \Rtil_e$ with respect to which $\phi_{\nu}$ is horizontal. 

Similarly for $q\in\dbN^{\ast}$ we define the extension 
$$(M\otimes_{W(k)} \Rtil_e/\Itil_e(q),F^1\otimes_{W(k)} \Rtil_e/\Itil_e(q),\Phi_{\nu(q)},\nabla_{\nu(q)},(t_{\alpha})_{\alpha\in\scrJ})$$ 
of $(\scrM_R^0,\nabla,(t_{\alpha})_{\alpha\in\scrJ})$ through a $W(k)$-homomorphism $\nu(q):R\to\Rtil_e/\Itil_e(q)=R_e/I_e(q)$ that maps $I$ to $\Itil_e(1)/\Itil_e(q)=I_e(1)/I_e(q)$. 

\medskip\smallskip\noindent
{\bf 5.3. Theorem.} {\it We assume that there exists a filtration lift $F^1_N$ of the pair $(H,(t_{H\alpha})_{\alpha\in\scrJ})$ to $R_e$. If $p=2$ we also assume that $(M,F^1,\phi)$ is the filtered Dieudonn\'e module of a $p$-divisible group $D$ over $\Spec\,W(k)$. Then there exists a $W(k)$-homomorphism $\nu:R\to\Rtil_e$ such that $\nu(I)\subseteq \Itil_e(1)$ and the quintuple $(N\otimes_{R_e} \Rtil_e,F^1_N\otimes_{R_e} \Rtil_e,\phi_N\otimes\Phi_k,\nabla_N,(t_{H\alpha})_{\alpha\in\scrJ})$ is isomorphic to the extension of $(\scrM_R^0,\nabla,(t_{\alpha})_{\alpha\in\scrJ})$ through $\nu$, under an isomorphism which modulo $\Itil_e(1)$ is $1_M$.}

\medskip
\proof
This proof is a group theoretical refinement of [Fa2, Sect. 7, pp. 135--136]. If $p\Ge 3$ let $D$ be the $p$-divisible group over $\Spec\,W(k)$ whose filtered Dieudonn\'e module is $(M,F^1,\phi)$. Thus $D$ is well defined for all primes $p$.

Let $\scrD$ be the unique $p$-divisible group over $\Spec\,R$ that lifts $D$ and such that the evaluation of its filtered Dieudonn\'e crystal at the thickening $(\Spec\,R/pR\hookrightarrow\Spec\,R,\delta(p))$ is (defined by) $(\scrM_{R}^0,\nabla)$, cf. Theorem 3.4.1 (b). The natural divided power structure of the ideal $(pX)$ of $\grS/(X^e)$ is nilpotent. Thus there exists a unique $p$-divisible group $H_{\grS/(X^e)}$ over $\Spec\,\grS/(X^e)$ that lifts both $H_{V/pV}$ and $D$ and such that 
the evaluation of the filtered Dieudonn\'e crystal of $H_{\grS/(X^e)}$ at the thickening $(\Spec\,V/pV\hookrightarrow\Spec\,\grS/(X^e),\delta(p))$ is $(N,F_N^1,\phi_N,\nabla_N)$ modulo $I_e(e)$, cf. Grothendieck--Messing deformation theory. 

By induction on $q\in\dbN^{\ast}$ we construct a $W(k)$-homomorphism $\nu(q):R\to \Rtil_e/\Itil_e(q)$ that maps $I$ to $\Itil_e(1)/\Itil_e(q)$ and such that the extension of $(\scrM_R^0,\nabla,(t_{\alpha})_{\alpha\in\scrJ})$ through $\nu(q)$ is isomorphic to $(N\otimes_{R_e} \Rtil_e,F^1_N\otimes_{R_e} \Rtil_e,\phi_N\otimes\Phi_k,\nabla_N,(t_{H\alpha})_{\alpha\in\scrJ})$ modulo $\Itil_e(q)$, under an isomorphism $\scrI_q$ which modulo $\Itil_e(1)/\Itil_e(q)$ is defined by $1_M$. Such an isomorphism $\scrI_q$ is unique as we have $\Phi_k^q(\Itil_e(1)/\Itil_e(q))=0$. The construction of $\nu(1)$ is obvious. 

For $q\ge 2$ the passage from $q-1$ to $q$ goes as follows. We lift $\nu(q-1):R\to \Rtil_e/\Itil_e(q-1)$ to an arbitrary $W(k)$-homomorphism $\tilde \nu(q):R\to \Rtil_e/\Itil_e(q)$. We endow the ideal $\Jtil_e(q):=\Itil_e(q-1)/\Itil_e(q)$ of $\Rtil_e/\Itil_e(q)$ with the trivial divided power structure: $\Jtil_e(q)^{[l]}=0$ if $l\Ge 2$. We consider an identification $\tilde \scrI_q$ between  a quadruple of the form $(M\otimes_{W(k)} \Rtil_e/\Itil_e(q),F^1_q,\phi^\prime_{\tilde \nu(q)},\nabla^\prime_{\tilde \nu(q)})$ and $(N\otimes_{R_e} \Rtil_e,F^1_N\otimes_{R_e} \Rtil_e,\phi_N\otimes\Phi_k,\nabla_N)$ modulo $\Itil_e(q)$ which modulo $\Jtil_e(q)$ is $\scrI_{q-1}$. Here $F^1_q$ is a direct summand of $M\otimes_{W(k)} \Rtil_e/\Itil_e(q)$ that lifts $F^1\otimes_{W(k)} \Rtil_e/\Itil_e(q-1)$, $\phi^\prime_{\tilde \nu(q)}$ lifts $\phi_{\nu(q-1)}$, and $\nabla^\prime_{\tilde \nu(q)}$ lifts $\nabla_{\nu(q-1)}$. In the next two paragraphs we show that we can choose $\tilde \scrI_q$ such that we have 
$$(\Phi^\prime_{\tilde \nu(q)},\nabla^\prime_{\tilde \nu(q)})=(\Phi_{\tilde \nu(q)},\nabla_{\tilde \nu(q)}).\leqno (24)$$
\indent
We first consider the case when $q-1<e$. By a second induction on $s\in\{1,\ldots,q-1\}$ we show that the pull-backs of $\scrD$ and $H_{\grS/(X^e)}$ to $\Spec\,\grS/(X^s)$ are isomorphic, under a unique isomorphism whose evaluation at the thickening $(\Spec\,k[[X]]/(X^s)\hookrightarrow \Spec\,\grS/(X^s),\delta(p))$ is $\scrI_s$; here the $W(k)$-homomorphism $R\to \grS/(X^s)$ used for the pull-back of $\scrD$ is $\nu(s)$. The case $s=1$ holds by constructions. Due to the existence of $\scrI_{s-1}$ for $2\le s<q$, the passage from $s-1$ to $s$ follows from the fact that $\Jtil_e(s)$ has a nilpotent divided power structure and from the Grothendieck--Messing deformation theory. This ends the second induction. Thus the Dieudonn\'e crystals of the pull-backs of $\scrD$ and $H_{\grS/(X^e)}$ to $\Spec\,\grS/(X^{q-1})$ are canonically identified. Thus as $\Jtil_e(q)$ has a nilpotent divided power structure, we can choose $\tilde \scrI_q$ such that (24) holds.

We next consider the case $q-1\Ge e$. We have $\Phi_k(\Itil_e(q-1))\subseteq p\Itil_e(q)$ and thus for $g_q\in\Ker(\pmb{GL}_N(\Rtil_e)\to \pmb{GL}_N(\Rtil_e/\Itil_e(q-1)))$ the $\Phi_k$-linear map $g_q(\phi_N\otimes\Phi_k) g_q^{-1}$ is of the form $h_q(\phi_N\otimes\Phi_k)$, where $h_q\in \pmb{GL}_N(\Rtil_e)$ is congruent modulo $\Itil_e(q)$ to $g_q$. When $g_q$ varies, $\Phi^\prime_{\tilde \nu(q)}$ varies by a left multiple of it by an arbitrary element $g_q^{\text{left}}\in\Ker(\pmb{GL}_M(\Rtil_e/\Itil_e(q))\to \pmb{GL}_M(\Rtil_e/\Itil_e(q-1)))$. We show that $\Phi^\prime_{\tilde \nu(q)}$ is of the form 
$$\tilde g_q\Phi_{\tilde \nu(q)},\;\;\; \text{where}\;\;\; \tilde g_q\in\Ker(\pmb{GL}_M(\Rtil_e/\Itil_e(q))\to \pmb{GL}_M(\Rtil_e/\Itil_e(q-1))).$$ 
As $\Phi_k(\Itil_e(q-1))\subseteq \Itil_e(q)$ and as $F^1_q$ and $F^1\otimes_{W(k)} \Rtil_e/\Itil_e(q)$ coincide modulo $\Jtil_e(q)$, we have $\Phi^\prime_{\tilde \nu(q)}(F^1\otimes_{W(k)} \Rtil_e/\Itil_e(q))=\Phi^\prime_{\tilde \nu(q)}(F^1_q)\subseteq pM\otimes_{W(k)} \Rtil_e/\Itil_e(q)$. Thus both $\Phi^\prime_{\tilde \nu(q)}$ and $\Phi_{\tilde \nu(q)}$ induce $\Rtil_e/\Itil_e(q)$-linear isomorphisms $(M+{1\over p}F^1)\otimes_{W(k)} {}_{\sigma} \Rtil_e/\Itil_e(q)\arrowsim M\otimes_{W(k)} \Rtil_e/\Itil_e(q)$. Thus indeed $\Phi^\prime_{\tilde \nu(q)}$ is of the form $\tilde g_q\Phi_{\tilde \nu(q)}$, where $\tilde g_q\in \pmb{GL}_M(\Rtil_e/\Itil_e(q))$. As $\Phi^\prime_{\tilde \nu(q)}$ and $\Phi_{\tilde \nu(q)}$ coincide modulo $\Jtil_e(q)$, $\tilde g_q$ modulo $\Jtil_e(q)$ is the identity element. We choose $g_q$ such that $g_q^{\text{left}}=\tilde g_q^{-1}$. By replacing $(N\otimes_{R_e} \Rtil_e,F^1_N\otimes_{R_e} \Rtil_e,\phi_N\otimes\Phi_k,\nabla_N)$ with its conjugate under $g_q$, we can choose $\tilde\scrI_q$ such that $\Phi^\prime_{\tilde \nu(q)}=\Phi_{\tilde \nu(q)}$. Thus $\nabla^\prime_{\tilde \nu(q)}=\nabla_{\tilde \nu(q)}$ (cf. Subsubsection 5.2.1) and therefore (24) holds.

From now on we take $\tilde\scrI_q$ such that (24) holds. As $\tilde\scrI_q$ lifts $\scrI_{q-1}$ and as no element of $\scrT(M)[{1\over p}]\otimes_{B(k)} \Jtil_e(q)[{1\over p}]$ is fixed by $\Phi_{\tilde \nu(q)}$, under the identification $\tilde\scrI_q^{-1}$ the image of $t_{H\alpha}$ in $\scrT(N\otimes_{R_e} \Rtil_e/\Itil_e(q)[{1\over p}])$ gets identified with the tensor $t_{\alpha}\in\scrT(M\otimes_{W(k)} \Rtil_e/\Itil_e(q)[{1\over p}])$ (here $\alpha\in\scrJ$). Thus we have 
$$t_{\alpha}\in F^0_q(\scrT(M\otimes_{W(k)} \Rtil_e/\Itil_e(q)))[{1\over p}],$$ 
where $(F^i_q(\scrT(M\otimes_{W(k)} \Rtil_e/\Itil_e(q))))_{i\in\dbZ}$ is the filtration of $\scrT(M\otimes_{W(k)} \Rtil_e/\Itil_e(q))$ defined by $F^1_q$. Let $\scrT(M)=\oplus_{i\in\dbZ} \tilde F^i(\scrT(M))$ be as in Subsection 2.6. We have $t_{\alpha}\in\tilde F^0(\scrT(M))[{1\over p}]$ for all $\alpha\in\scrJ$ and the filtration of $\scrT(M)$ defined by $F^1$ is $(F^i(\scrT(M)))_{i\in\dbZ}$, where $F^i(\scrT(M)):=\oplus_{s=i}^{\infty} \tilde F^s(\scrT(M))$. Let $\tilde F^{-1}(\End_{W(k)}(M))$, $U_{\text{big}}$, and $U$ be as in Subsection 2.6.

We identify $1_M$ with $1_{M\otimes_{W(k)} \Rtil_e/\Itil_e(q)}$. Let $v_q\in \tilde F^{-1}(\End_{W(k)}(M))\otimes_{W(k)} \Jtil_e(q)$ be the unique element such that $(1_M+v_q)(F^1\otimes_{W(k)} \Rtil_e/\Itil_e(q))=F^1_q$. We have 
$$(1_M-v_q)(t_{\alpha})\in F^0(\scrT(M))[{1\over p}]\otimes_{B(k)} \Rtil_e/\Itil_e(q)[{1\over p}].\leqno (25)$$ 
We view $\scrT(M)$ as a module over $\End_{W(k)}(M)$ and thus also over $\Lie(U)$. As $v_q\in\tilde F^{-1}(\scrT(M))$ and $t_{\alpha}\in\tilde F^0(\scrT(M))[{1\over p}]$, the component of $(1_M-v_q)(t_{\alpha})$ in $\tilde F^{-1}(\scrT(M))[{1\over p}]\otimes_{B(k)} \Rtil_e/\Itil_e(q)[{1\over p}]$ is $-v_q(t_{\alpha})$ as well as $0$ (cf. (25)). Thus $v_q$ annihilates $t_{\alpha}$ for all $\alpha\in\scrJ$. But $\Lie(G_{B(k)})$ is the Lie subalgebra of $\End_{B(k)}(M[{1\over p}])$ that annihilates $t_{\alpha}$ for all $\alpha\in\scrJ$. Thus, as $\Jtil_e(q)$ is a free $W(k)$-module of rank 1, we have
$$v_q\in \Lie(G_{B(k)})\otimes_{B(k)} [\Rtil_e/\Itil_e(q)[{1\over p}]]\cap [\tilde F^{-1}(\End_{W(k)}(M))\otimes_{W(k)} \Jtil_e(q)]=\Lie(U)\otimes_{W(k)} \Jtil_e(q).$$ 
\indent
We identify $\Hom_{W(k)}(F^1,F^0)$ with $\End_{W(k)}(M)/\End_{W(k)}(M)\cap F^0(\scrT(M))$ and with $\Lie(U_{\text{big}})$. As $\nabla$ respects the $G$-action (see Subsection 3.4), the image of the Kodaira--Spencer map $\Theta$ of $\nabla$ is contained in the image of $(\Lie(G_{B(k)})\cap \End_{W(k)}(M))\otimes_{W(k)} R$ in $\Hom_{W(k)}(F^1,F^0)\otimes_{W(k)} R$. Thus $\text{Im}(\Theta)$ is contained in the direct summand $\Lie(U)\otimes_{W(k)} R$ of $\Hom_{W(k)}(F^1,F^0)\otimes_{W(k)} R$. We recall that we have $R=W(k)[[z_1,\ldots,z_m]]$, $\Phi_{R}(z_j)=z_j^p$ for each $j\in\{1,\ldots,m\}$, and $\tau$ is the universal element of $U(R)$. With $\Delta$ as in Subsubsection 3.1.1 (b), $\nabla$ modulo $(p,(z_1,z_2,\ldots,z_m)^{p-1})$ is $\Delta+\tau^{-1}d\tau$ modulo $(p,(z_1,z_2,\ldots,z_m)^{p-1})$, cf. Equations (11) and (12). As $U$ is a closed subgroup scheme of $G$ and $\tau$ is the universal element of $U(R)$, we get that the $R$-submodule $\text{Im}(\Theta)$ of $\Lie(U)\otimes_{W(k)} R$ surjects onto $\Lie(U)\otimes_{W(k)} k$. Thus $\text{Im}(\Theta)=\Lie(U)\otimes_{W(k)} R$.

As $v_q\in \Lie(U)\otimes_{W(k)} \Jtil_e(q)$ we can write $v_q=\sum_{j=1}^m x_j\Theta({{\partial}\over {\partial z_j}})$, where $x_j\in\Jtil_e(q)$ and where we denote also by $\Theta({{\partial}\over {\partial z_j}})$ its reduction modulo $\Itil_e(q)$. Let $\nu(q):R\to \Rtil_e/\Itil_e(q)$ be the $W(k)$-homomorphism that takes $z_l$ to $\tilde \nu(q)(z_j)-x_j$. By replacing $\tilde \nu(q)$ with $\nu(q)$, $F^1_q$ gets replaced by $(1_M-v_q)(F^1_q)=F^1\otimes_{W(k)} \Rtil_e/\Itil_e(q)$. Thus $\scrI_q$ is defined by the identification of $(M\otimes_{W(k)} \Rtil_e/\Itil_e(q),F^1\otimes_{W(k)} \Rtil_e/\Itil_e(q),\Phi_{\nu(q)},\nabla_{\nu(q)},(t_{\alpha})_{\alpha\in\scrJ})$ with $(N\otimes_{R_e} \Rtil_e,F^1_N\otimes_{R_e} \Rtil_e,\phi_N\otimes\Phi_k,\nabla_N,(t_{H\alpha})_{\alpha\in\scrJ})$ modulo $\Itil_e(q)$. This ends the induction.

Let $\nu:R\to\tilde R_e$ be such that it lifts all $\nu(q)$'s. The extension of $(\scrM_R^0,\nabla,(t_{\alpha})_{\alpha\in\scrJ})$ through $\nu$ is isomorphic to $(N\otimes_{R_e} \Rtil_e,F^1_N\otimes_{R_e} \Rtil_e,\phi_N\otimes\Phi_k,\nabla_N,(t_{H\alpha})_{\alpha\in\scrJ})$, under an isomorphism which modulo $\Itil_e(q)$ is $\scrI_q$ for all $q\in\dbN^{\ast}$.\endproof

\medskip\smallskip\noindent 
{\bf 5.4. Theorem.} {\it We assume that $k=\bar k$ and  $p>2$. We also assume that there exists a filtration lift $F^1_N$ of the pair $(H,(t_{H\alpha})_{\alpha\in\scrJ})$ to $R_e$ (for instance, this holds if $\scrG$ is a reductive group scheme). Then there exists an isomorphism 
$$\rho_{R_e}:(N,(t_{H\alpha})_{\alpha\in\scrJ})\arrowsim (H^1(H)\otimes_{\dbZ_p} R_e,(v_{H\alpha})_{\alpha\in\scrJ}).$$ 
Thus $\scrG$ is a flat, closed subgroup scheme of $\pmb{GL}_N$ that is isomorphic to either $G\times_{\Spec\,W(k)} \Spec\,R_e$ or $\scrG^{\acute et}\times_{\Spec\,\dbZ_p} \Spec\,R_e$.}

\medskip
\proof
Let $D$ be as in the proof of Theorem 5.3. As $p>2$, let $\scrD$ be as in the proof of Theorem 5.3 (i.e., as in Theorem 3.4.1 (b)). Let $(\scrV_{\alpha})_{\alpha\in\scrJ}$ be as in Subsubsection 3.4.2; we have $\scrV_{\alpha}\in\scrT(H^1(\scrD_{K_R})[{1\over p}])$ (we recall that $Q=R$). Let $(v_{\alpha})_{\alpha\in\scrJ}$ be as in the end of Subsubsection 2.2.4. Let $\nu:R\to\Rtil_e$ be as in Theorem 5.3. Let $H_1$ be the $p$-divisible group over $\Spec\,V$ that is the pull-back of $\scrD$ via the composite morphism $\Spec\,V=\Spec\,\Rtil_e/\Itil_e(1)\hookrightarrow\Spec\,\Rtil_e\to\Spec\,R$ defined naturally by $\nu$. Let the quadruple $(N_1,\phi_{N_1},\nabla_{N_1},(t_{H_1\alpha})_{\alpha\in\scrJ})$ be the analogue of the quadruple $(N,\phi_N,\nabla_N,(t_{H\alpha})_{\alpha\in\scrJ})$ but for $H_1$ and the extension of $(t_{\alpha})_{\alpha\in\scrJ}$ via the $W(k)$-homomorphism $R\to V$ that defines the mentioned composite morphism. We can identify $(N_1,(t_{H_1\alpha})_{\alpha\in\scrJ})=(M\otimes_{W(k)} R_e,(t_{\alpha})_{\alpha\in\scrJ})$, cf. also Subsubsection 5.2.8. Let $(v_{H_1\alpha})_{\alpha\in\scrJ}$ be the analogue of $(v_{H\alpha})_{\alpha\in\scrJ}$ but for $(H_1,(t_{H_1\alpha})_{\alpha\in\scrJ})$ instead of $(H,(t_{H\alpha})_{\alpha\in\scrJ})$.

Let 
$\tilde f:(N_1\otimes_{R_e} \Rtil_e,\phi_{N_1}\otimes\Phi_k,\nabla_{N_1},(t_{H_1\alpha})_{\alpha\in\scrJ})\arrowsim (N\otimes_{R_e} \Rtil_e,\phi_{N}\otimes\Phi_k,\nabla_{N},(t_{H\alpha})_{\alpha\in\scrJ})$
be an isomorphism 
that takes the direct summand of $N_1\otimes_{R_e} \Rtil_e$ which is the Hodge filtration of $\nu^*(\scrD)$ onto $F^1_N\otimes_{R_e} \Rtil_e$, cf. Theorem 5.3. Let $\phi_{N_1N}$ be the $\Phi_k$-linear endomorphism of $\Hom_{R_e}(N_1,N)[{1\over p}]$ such that for $x\in\Hom_{R_e}(N_1,N)$ and $y\in N_1$ we have $\phi_{N_1N}(x)(\phi_{N_1}(y))=\phi_N(x(y))\in N$. As $\Phi_k(\Rtil_e)\subseteq R_e$ (see Subsection 5.1) we get that 
$$\tilde f=(\phi_{N_1N}\otimes\Phi_k)(\tilde f)\in (\Hom_{R_e}(N_1,N)\otimes_{R_e} \Rtil_e)\cap \Hom_{R_e}(N_1,N)[{1\over p}]=\Hom_{R_e}(N_1,N).$$ 
Thus $\tilde f$ is the extension to $\Rtil_e$ of an isomorphism $f$ between $(N_1,\phi_{N_1},\nabla_{N_1},(t_{H_1\alpha})_{\alpha\in\scrJ})$ and $(N,\phi_{N},\nabla_{N},(t_{H\alpha})_{\alpha\in\scrJ})$. As $f$ respects also the Verschiebung maps (of $\phi_{N_1}$ and $\phi_N$), $f$ defines naturally an isomorphism $f_{V/pV}:\dbD(H_{1,V/pV})\arrowsim\dbD(H_{V/pV})$. 

As $V/pV=k[X]/(X^e)$, $H_{V/pV}$ is uniquely determined by its Dieudonn\'e crystal $\dbD(H_{V/pV})$ (cf. [BM, Rm. 4.3.2 (i)]). Thus there exists a unique isomorphism $h_{V/pV}:H_{1,V/pV}\arrowsim H_{V/pV}$ such that $\dbD(h_{V/pV})=f_{V/pV}$. As $\tilde f$ maps the Hodge filtration of $\nu^*(\scrD)$ onto $F^1_N\otimes_{R_e} \Rtil_e$, the reduction of $\tilde f$ modulo $\Ker(\tilde e_V)$ maps the Hodge filtration of $H_1$ onto $F^1_H$ and thus (as $p>2$) the isomorphism $h_{V/pV}$ lifts naturally to an isomorphism $h_V:H_1\arrowsim H$ (cf. Grothendieck--Messing deformation theory). 

Thus we can identify $(N,(t_{H\alpha})_{\alpha\in\scrJ})$ with $(N_1,(t_{H_1\alpha})_{\alpha\in\scrJ})=(M\otimes_{W(k)} R_e,(t_{\alpha})_{\alpha\in\scrJ})$. As $(D,(v_{\alpha})_{\alpha\in\scrJ})$ and $(H_1,(v_{H_1\alpha})_{\alpha\in\scrJ})$ are pull-backs of $(\scrD,(\scrV_{\alpha})_{\alpha\in\scrJ})$, as in the proof of Lemma 3.4.3 we argue that we can identify  $(H^1(\scrD_{K_R}),(\scrV_{\alpha})_{\alpha\in\scrJ})$ with either $(H^1(D),(v_{\alpha})_{\alpha\in\scrJ})$ or $(H^1(H_1),(v_{H_1\alpha})_{\alpha\in\scrJ})$. There exist isomorphisms $(M,(t_{\alpha})_{\alpha\in\scrJ})\arrowsim (H^1(D)\otimes_{\dbZ_p} W(k),(v_{\alpha})_{\alpha\in\scrJ})$, cf. Main Theorem applied to $(D,(t_{\alpha})_{\alpha\in\scrJ})$. From the last three sentences we get the first part of the theorem. The second part of the theorem follows directly from the first part.\endproof

\medskip\smallskip\noindent
{\bf 5.5. A complement to [Ki3].} In this subsection we assume that $p>2$. Let $(\grN,\varphi_{\grN})$ be the contravariant {\it Breuil window} of $H$ relative to the $W(k)$-epimorphism $\grS\twoheadrightarrow V$ that maps $X$ to $\pi_V$. We recall that $\grN$ is a free $\grS$-module of the same rank as $N$ and that $\varphi_{\grN}:\grN\otimes_{\grS} {}_{\Phi_k} \grS\to\grN$ is a $\grS$-linear map whose cokernel is annihilated by $f_e$. To each $v_{H\alpha}$ corresponds naturally a tensor $w_{H\alpha}\in\scrT(\grN[{1\over p}])$ (cf. [Ki1]) and one can speak about the schematic closure $\grG$ in $\pmb{GL}_{\grN}$ of the closed subgroup scheme of $\pmb{GL}_{\grN[{1\over p}]}$ that fixes $w_{H\alpha}$ for all $\alpha\in\scrJ$. The canonical way of passing from Breuil windows to $F$-crystals (see [Ki1], [Ki2], [Ki3], and [Lau]) shows that one has a canonical identification
$$(\grN,(w_{H\alpha})_{\alpha\in\scrJ})\otimes_{\grS} {}_{\Phi_k} R_e=(N,(t_{H\alpha})_{\alpha\in\scrJ}).$$
In particular, if $\grG$ is a group scheme over $\Spec\,\grS$ which is flat, then we have
$$\scrG=\grG\otimes_{\Spec\,\grS} {}_{\Phi_k} \Spec\,R_e\leqno (26)$$ 
and thus $\scrG$ is a flat, closed subgroup scheme of $\pmb{GL}_N$. 

\medskip\noindent
{\bf 5.5.1. Corollary.} {\it We assume that $p>2$. Then $\scrG$ is a reductive, closed subgroup scheme of $\pmb{GL}_N$ if and only if $\scrG^{\acute et}$ is a reductive, closed subgroup scheme of $\pmb{GL}_{H^1(H)}$.}

\medskip
\proof
If $\scrG^{\acute et}$ is a reductive group scheme, then (cf. [Ki3, Cor. (1.4.3)]) there exists an isomorphism $(\grN,(w_{H\alpha})_{\alpha\in\scrJ})\arrowsim (H^1(H)\otimes_{\dbZ_p} \grS,(v_{H\alpha})_{\alpha\in\scrJ})$ and from this and (26) we get that $\scrG$ is a reductive group scheme.   

Suppose now that $\scrG$ is a reductive, closed subgroup scheme of $\pmb{GL}_N$. To prove that $\scrG^{\acute et}$ is reductive, we can assume that $k=\bar k$. Thus $\scrG$ is split. As there exists a filtration lift of the pair $(H,(t_{H\alpha})_{\alpha\in\scrJ})$ to $R_e$ (cf. Lemma 5.2.6), $\scrG^{\acute et}$ is reductive (cf. Theorem 5.4).\endproof

\medskip\smallskip\noindent
{\bf 5.6. Counterexamples.} Suppose that $e>1$, that $k=\bar k$, and that there exists an embedding $V\hookrightarrow\dbC$. Let $d\in\dbN^{\ast}$. Let $\scrA$ be an abelian scheme over $\Spec\,V$ which has complex multiplication and relative dimension $d$. Let $\scrA_k$ be the special fibre of $\scrA$. We consider a semisimple, commutative $\dbQ$--subalgebra $B$ of $\End(\scrA)\otimes_{\dbZ} \dbQ$ of dimension $2d$. There exist examples in which $d=1$, $\scrA_k$ is a supersingular elliptic curve, $B:=\dbQ(\sqrt{-p})$, and $e=2$. Let $H$ be the $p$-divisible group of $\scrA$.   

As $B$ is a product of number fields, we can speak about the ring of integers $O_B$ of $B$. Let $\bar O_B$ be the largest subring of $O_B$ such that $\scrA$ has complex multiplication by $\bar O_B$. By replacing $\scrA$ with an abelian scheme isogenous to it, we can assume that $\bar O_B=\dbZ+pO_B$ and that the $\bar O_B\otimes_{\dbZ} \dbZ_p$-module $H^1(H)$ is isomorphic to $\bar O_B\otimes_{\dbZ} \dbZ_p$ (note that this operation might enlarge both $V$ and $e$). Let $\{t_{H\alpha}|\alpha\in\scrJ\}$ be the set of tensors of $\scrT(N)[{1\over p}] $ that are crystalline realizations of Hodge cycles on the generic fibre of $\scrA$.  The group schemes $G_{B(k)}$, $G_{R_e[{1\over p}]}$, and $G^{\acute et}_{\dbQ_p}$ are forms of the Mumford--Tate group of $\scrA_{\dbC}$ (see [De2, Sect. 3]) and thus, as $\scrA_{\db C}$ has complex multiplication, are tori. The set $\{t_{H\alpha}|\alpha\in\scrJ\}$ includes the endomorphisms of $N[{1\over p}]$ that are crystalline realizations of elements of $\bar O_B$. Thus $\bar O_B\otimes_{\dbZ} \dbZ_p$ is a $\dbZ_p$-subalgebra of $\End_{\dbZ_p}(H^1(H))$ which as a $\dbZ_p$-submodule is a direct summand. 

\smallskip
{\bf (a)} We show that the assumption that either the condition 5.2.3 (i) or the condition 5.2.3 (ii) holds leads to a contradiction. We can assume that the condition 5.2.3 (ii) holds i.e., there exists an isomorphism $\rho:(M,(t_{\alpha})_{\alpha\in\scrJ})\arrowsim (H^1(H)\otimes_{\dbZ_p} W(k),(v_{H\alpha})_{\alpha\in\scrJ})$. We write $\rho\phi\rho^{-1}=t(1_{H^1(H)}\otimes\sigma)$, where $t\in G^{\acute et}_{\dbQ_p}(B(k))$. Thus $t$ is an endomorphism of $H^1(H)\otimes_{\dbZ_p} W(k)$ whose Hodge slopes are $1$ and $0$ with the same multiplicity $d$ and which centralizes $\bar O_B$.  Thus we have $t\in\bar O_B\otimes_{\dbZ} W(k)=W(k)1_M+pO_B\otimes_{\dbZ} W(k)$ and $\text{det}(t)\in p^d\dbG_m(W(k))$. But as the $\bar O_B\otimes_{\dbZ} \dbZ_p$-module $H^1(H)$ is isomorphic to $\bar O_B\otimes_{\dbZ} \dbZ_p$, the determinant of $t\in\bar O_B\otimes_{\dbZ} W(k)$ is either a unit or divisible by $p^{2d}$. Contradiction. 

\smallskip
{\bf (b)} Suppose $\scrA_k$ is a supersingular elliptic curve (so $d=1$). Then $G$ is the group scheme of invertible elements of the $W(k)$-algebra $\Lie(G_{B(k)})\cap\End_{W(k)}(M)$ and thus it is smooth over $\Spec\,W(k)$. We look at the $\bar O_B$-module $M$. The element $\varPi:=p\sqrt{-p}$ of $\bar O_B$ is an endomorphism of $(M,\phi)$ such that $\varPi^2=-p^3$. Using a $W(k)$-basis $\{e_1,e_2\}$ for $M$ such that $\phi(e_1)=e_2$ and $\phi(e_2)=pe_1$, it is easy to see that the endomorphism $\varPi$ of $(M,\phi)$ must be divisible by $p$. This implies that we have inclusions $\bar O_B\hookrightarrow O_B\hookrightarrow \End(M,\phi)$. As $O_B\otimes_{\dbZ} W(k)$ is a discrete valuation ring, the $O_B\otimes_{\dbZ} W(k)$-module $M$ is free of rank $1$. Thus the $\bar O_B$-module $M/pM$ is isomorphic to $O_B/pO_B\otimes_{\dbF_p} k$. As we have functorial identifications $M/pM=H^1_{\text{dR}}(H/V)\otimes_V k=(N/F^1(R_e)N)\otimes_V k=N\otimes_{R_e} k$, the $\bar O_B$-module $(N/F^1(R_e)N)\otimes_V k$ is isomorphic to $O_B/pO_B\otimes_{\dbF_p} k$ and thus it is not isomorphic to the $\bar O_B$-module $H^1(H)\otimes_{\dbZ_p} k$. Therefore none of the conditions 5.2.3 (i) to (iii) holds. 

\medskip
The following abstract and trivial fact captures part of the very essence of (b).

\medskip\noindent
{\bf 5.6.1. Fact.} {\it Suppose there exists an endomorphism $f$ of $H_k$ which does not lift to $H_V$ but whose multiplication by $p$ lifts to $H_V$. We also assume that there exists $\alpha\in\scrJ$ such that $t_{H\alpha}$ is the crystalline realization of $pf$. Then the conditions 5.2.3 (i) and (ii) do not hold.}

\medskip\noindent
{\bf Acknowledgment.} We would like to thank FIM of ETH--Z\"urich, University of California at Berkeley, University of Utah, University of Arizona, Binghamton University, and IAS--Princeton for providing good working conditions. We would like to thank Faltings for some comments and a suggestion that led to a shorter Subsection 4.2. We would also like to thank Deligne for many valuable comments and suggestions that led to the addition of Corollary 1.4 and to a correction of Subsection 5.6. We would also like to thank Gabber for the first proof of Theorem 2.4.1 (b) and for the example of the proof of Theorem 2.4.1 (d). We would also like to thank the referees for many valuable comments and suggestions which in particular led to a shorter and much better way to present Section 3. We are very much obliged to Milne's conjecture of 1995 and to [Fa2]: in our work on integral aspects of Shimura varieties they have been by far the most inspiring moments (things). This research was partially supported by the NSF grants DMF 97-05376 and DMS \#0900967.

\bigskip
\references{37}
{\nspace{

\Ref[An]
Y. Andr\'e,
\sl On the Shafarevich and Tate conjectures for hyperk\"ahler varieties,
\rm Math. Ann. {\bf 305} (1996), no. 2,  205--248.

\Ref[Be] P. Berthelot, 
\sl Cohomologie cristalline des sch\'emas de caract\'eristique $p>0$, 
\rm Lecture Notes in Math., Vol. {\bf 407}, Springer-Verlag, Berlin-New York, 1974.

\Ref[Bo]
A. Borel,
\sl Linear algebraic groups,
\rm Grad. Texts in Math., Vol. {\bf 126}, Springer-Verlag, New York, 1991.

\Ref[BBM]
P. Berthelot, L. Breen, and W. Messing,
\sl Th\'eorie de Dieudonn\'e cristalline II,
\rm Lecture Notes in Math., Vol. {\bf 930}, Springer-Verlag, Berlin, 1982.

\Ref[BLR]
S. Bosch, W. L\"utkebohmert, and M. Raynaud,
\sl N\'eron models,
\rm Ergebnisse der Mathematik und ihrer Grenzgebiete (3), Vol. {\bf 21}, Springer-Verlag, Berlin, 1990.

\Ref[BM]  P. Berthelot and W. Messing, 
\sl Th\'eorie de Dieudonn\'e cristalline III. Th\'eor\`emes d'\'equivalence et de pleine fid\'elit\'e,
\rm The Grothendieck Festschrift, Vol. {\bf I},   173--247, Progr. Math., Vol. {\bf 86}, Birkh\"auser Boston, Boston, MA, 1990.

\Ref[dJ] J. de Jong, 
\sl Crystalline Dieudonn\'e module theory via formal and rigid geometry, 
\rm Inst. Hautes \'Etudes Sci. Publ. Math., Vol. {\bf 82},  5--96, 1995.

\Ref[De1] P. Deligne, 
\sl Cristaux ordinaires et coordon\'ees canoniques, 
\rm Algebraic surfaces (Orsay, 1976--78), Lecture Notes in Math., Vol. {\bf 868},  80--137, Springer-Verlag, Berlin-New York, 1981. 

\Ref[De2]
P. Deligne,
\sl Hodge cycles on abelian varieties,
\rm Hodge cycles, motives, and Shimura varieties, Lecture Notes in Math., Vol. {\bf 900},  9--100, Springer-Verlag, Berlin-New York, 1982.

\Ref[Dem] M. Demazure, 
\sl Lectures on $p$-divisible groups, 
\rm Lecture Notes in Math., Vol. {\bf 302}, Springer-Verlag, Berlin-New York, 1972.

\Ref[Fa1]
G. Faltings,
\sl Crystalline cohomology and $p$-adic Galois representations,
\rm Algebraic analysis, geometry, and number
theory,  25--79, Johns Hopkins Univ. Press, Baltimore, MD, 1989.

\Ref[Fa2]
G. Faltings,
\sl Integral crystalline cohomology over very ramified
valuation rings,
\rm J. Amer. Math. Soc. {\bf 12} (1999), no. 1,  117--144.

\Ref[Fa3]
G. Faltings,
\sl Almost \'etale extensions,
\rm Ast\'erisque {\bf 279},  185--270, Soc. Math. France, Paris, 2002.

\Ref[Fo1]
J.-M. Fontaine,
\sl Groupes $p$-divisibles sur les corps locaux, 
\rm Ast\'erisque {\bf 47/48}, Soc. Math. France, Paris, 1977.

\Ref[Fo2] 
J.-M. Fontaine,
\sl Cohomologie de de Rham, cohomologie cristalline et repr\'esentations $p$-adiques,
\rm Algebraic geometry (Tokyo/Kyoto, 1982),  86--108, Lecture Notes in Math., Vol. {\bf 1016}, Springer-Verlag, Berlin-New York, 1983.

\Ref[Fo3] 
J.-M. Fontaine, 
\sl Le corps des p\'eriodes $p$-adiques, 
\rm P\'eriodes $p$-adiques (Bures-sur-Yvette, 1988), Ast\'erisque {\bf 223},  59--101, Soc. Math. France, Paris, 1994.

\Ref[Fo4] 
J.-M. Fontaine, 
\sl Repr\'esentations $p$-adiques semi-stables, 
\rm P\'eriodes $p$-adiques (Bures-sur-Yvette, 1988), Ast\'erisque {\bf 223},  113--185, Soc. Math. France, Paris, 1994.

\Ref[FC]
G. Faltings and C.-L. Chai,
\sl Degeneration of abelian varieties,
\rm Ergebnisse der Mathematik und ihrer Grenzgebiete (3), Vol. {\bf 22}, Springer-Verlag, Heidelberg, 1990.

\Ref[FL] J.-M. Fontaine and G. Laffaille, 
\sl Construction de repr\'esentations $p$-adiques, 
\rm Ann. Sci. \'Ecole Norm. Sup. {\bf 15} (1982), no. 4,  547--608.

\Ref[Fu]
W. Fulton, 
\sl Intersection theory,
\rm  Ergebnisse der Mathematik und ihrer Grenzgebiete (3), Vol. {2}, Springer-Verlag, Berlin, 1984.

\Ref[Go]
D. Goss,
\sl Basic structures of function field arithmetic,
\rm Second corrected printing, Springer-Verlag, Berlin Heidelberg, 1998. 

\Ref[Gr] A. Grothendieck et al.,
\sl Rev\^etements \'etales et groupe fondamental, 
\rm S\'eminaire de g\'eom\'etrie alg\'ebrique du Bois Marie 1960-61 (SGA 1), Lecture Notes in Math., Vol. {\bf 224}, Springer-Verlag, Berlin-New York, 1971.

\Ref[GHKR] U. G\"ortz, T. Haines, R. Kottwitz, and D. Reuman,
\sl Affine Deligne--Lusztig varieties in affine flag varieties,
\rm Compos. Math. {\bf 146}  (2010),  no. 5, 1339--1382. 

\Ref[Ki1] M. Kisin,
\sl Crystalline representations and $F$-crystals, 
\rm Algebraic geometry and number theory,  459--496, Progr. Math., Vol. {\bf 253}, Birkh\"auser,
  Boston, MA, 2006.

\Ref[Ki2] M. Kisin,
\sl Modularity of 2-adic Barsotti--Tate representations, \rm Invent. Math. {\bf 178} (2009), no. 3,  587--634.

\Ref[Ki3] M. Kisin,
\sl Integral canonical models of Shimura varieties of abelian type, 
\rm J. Amer. Math. Soc. {\bf 23} (2010), no. 4,  967--1012.

\Ref[Ko]
R. E. Kottwitz,
\sl Points on some Shimura varieties over finite fields,
\rm J. Amer. Math. Soc. {\bf 5} (1992), no. 2,  373--444.

\Ref[Il] 
L. Illusie, 
\sl D\'eformations des groupes de Barsotti--Tate (d'apr\`es A. Grothendieck), 
\rm Seminar on arithmetic bundles: the Mordell conjecture (Paris, 1983/84),  151--198, Ast\'erisque {\bf 127}, Soc. Math. France, Paris, 1985.

\Ref[La]
G. Laffaille, 
\sl Groupes $p$-divisibles et modules filtr\'es: le cas peu ramifi\'e,
\rm Bull. Soc. Math. France  {\bf 108}  (1980), no. 2,  187--206.

\Ref[Lau]
E. Lau,
\sl A relation between Dieudonn\'e displays and crystalline Dieudonn\'e theory,
\rm manuscript, Dec. 2010, http://arxiv.org/abs/1006.2720.

\Ref[Ly] G. Lyubeznik, 
\sl $F$-modules: applications to local cohomology and $\scrD$-modules in characteristic $p>0$,
\rm J. Reine Angew. Math. {\bf 491} (1997),  65--130.

\Ref[LR] R. Langlands and M. Rapoport,
\sl Shimuravariet\"aten und Gerben, 
\rm J. Reine Angew. Math. {\bf 378} (1987),  113--220. 

\Ref[Ma] H. Matsumura, 
\sl Commutative algebra. Second edition, 
\rm The Benjamin/Cummings Publ. Co., 1980.

\Ref[Me] W. Messing, 
\sl The crystals associated to Barsotti--Tate groups: with applications to abelian schemes, 
\rm Lecture Notes in Math., Vol. {\bf 264}, Springer-Verlag, Berlin-New York, 1972.

\Ref[Mi1]
J. S. Milne,
\sl \'Etale cohomology,
\rm Princeton Mathematical Series, Vol. {\bf 33}, Princeton Univ. Press, Princeton, NJ, 1980.

\Ref[Mi2]
J. S. Milne,
\sl The points on a Shimura variety modulo a prime of good
reduction,
\rm The Zeta functions of Picard modular surfaces,  153--255, Univ. Montr\'eal Press, Montreal, QC, 1992.

\Ref[Mi3]
J. S. Milne,
\sl Shimura varieties and motives,
\rm Motives (Seattle, WA, 1991),  447--523, Proc. Sympos. Pure Math., Vol. {\bf 55}, Part 2, Amer. Math. Soc., Providence, RI, 1994.

\Ref[Mi4]
J. S. Milne,
\sl Points on Shimura varieties over finite fields: the conjecture of Langlands and Rapoport,
\rm manuscript, Nov. 2009, http://arxiv.org/abs/0707.3173.

\Ref[NV]
M.-H. Nicole and A.~Vasiu,
\sl Minimal truncations of supersingular $p$-divisible groups,
\rm Indiana Univ. Math. J. {\bf 56} (2007), no. 6, 2887--2898.

\Ref[Ra1] M. Raynaud, 
\sl Anneaux locaux hens\'eliens, 
\rm Lecture Notes in Math., Vol. {\bf 169}, Springer-Verlag, Berlin-New York, 1970.

\Ref[Ra2]
M. Raynaud,
\sl Sch\'emas en groupes de type $(p,\ldots,p)$,
\rm Bull. Soc. Math. France {\bf 102} (1974),  241--280.

\Ref[Va1]
A. Vasiu,
\sl Integral canonical models for Shimura varieties of preabelian type,
\rm Asian J. Math. {\bf 3} (1999), no. 2,  401--518.

\Ref[Va2]
A. Vasiu,
\sl A purity theorem for abelian schemes,
\rm Michigan Math. J. {\bf 52} (2004), no. 1,  71--81.

\Ref[Va3]
A. Vasiu,
\sl Geometry of Shimura varieties of Hodge type over finite fields,
\rm Proceedings of the NATO Advanced Study Institute on Higher dimensional geometry over finite fields, G\"ottingen, Germany, June 25 - July 06 2007, 197--243, IOS Press. 

\Ref[Va4]
A. Vasiu,
\sl Level $m$ stratifications of versal deformations of $p$-divisible groups,
\rm J. Alg. Geom. {\bf 17} (2008), no. 4, 599--641.

\Ref[Va5]
A. Vasiu,
\sl Manin problems for Shimura varieties of Hodge type,
\rm J. Ramanujan Math. Soc. {\bf 26} (2011), no. 1, 31--84.

\Ref[Va6]
A. Vasiu,
\sl Generalized Serre--Tate ordinary theory,
\rm manuscript, Dec. 2007, 

http://arxiv.org/abs/math/0208216.

\Ref[Va7]
A. Vasiu,
\sl CM-lifts of isogeny classes of Shimura $F$-crystals over finite fields,
\rm manuscript, July 2007, http://arxiv.org/abs/math/0304128.

\Ref[Va8]
A. Vasiu,
\sl Good reductions of Shimura varieties of Hodge type in arbitrary unramified mixed characteristic, Parts I and II,
\rm manuscripts, 2007, http://arxiv.org/abs/0707.1668 and http://arxiv.org/pdf/0712.1572.

\Ref[Va9]
A. Vasiu,
\sl On the Tate and Langlands--Rapoport conjecture for special fibres of integral canonical models of Shimura varieties of abelian type,
\rm manuscript, Nov. 2008, http://www.math.binghamton.edu/adrian.

\Ref[Wi]
J.-P. Wintenberger,
\sl Un scindage de la filtration de Hodge pour certaines vari\'et\'es alg\'ebriques sur les corps locaux,
\rm Ann. of Math. (2) {\bf 119} (1984), no. 3,  511--548.

}}

\bigskip\medskip
\hbox{Adrian Vasiu}
\hbox{Department of Mathematical Sciences, Binghamton University,}
\hbox{Binghamton, P. O. Box 6000, New York 13902-6000, U.S.A.}
\hbox{e-mail: adrian\@math.binghamton.edu\;\;\;fax: 1-607-777-2450} 
\enddocument